\input amstex\documentstyle {amsppt}  
\pagewidth{12.5 cm}\pageheight{19 cm}\magnification\magstep1
\topmatter
\title Lectures on Hecke algebras with unequal parameters\endtitle
\author G. Lusztig\endauthor
\endtopmatter
\document     
\define\lan{\langle}
\define\ran{\rangle}

\define\mto{\mapsto}
\define\lra{\leftrightarrow}
\define\lras{\leftrightarrows}
\define\Lra{\Leftrightarrow}

\define\sm{\smallmatrix}
\define\esm{\endsmallmatrix}
\define\sub{\subset}
\define\boxt{\boxtimes}

\define\tim{\times}
\define\ti{\tilde}
\define\nl{\newline}
\redefine\i{^{-1}}
\define\fra{\frac}
\define\und{\underline}
\define\ov{\overline}
\define\ot{\otimes}

\define\Hom{\text{\rm Hom}}
\define\End{\text{\rm End}}

\define\sgn{\text{\rm sgn}}
\define\tr{\text{\rm tr}}
\define\Tr{\text{\rm Tr}}

\define\bst{{}^\bigstar\!}

\define\spa{\spadesuit}

\redefine\i{^{-1}}

\define\al{\alpha}

\define\ga{\gamma}
\define\de{\delta}
\define\ep{\epsilon}
\define\io{\iota}
\define\om{\omega}
\define\si{\sigma}

\define\ka{\kappa}
\define\la{\lambda}
\define\ze{\zeta}

\define\Ga{\Gamma}
\define\De{\Delta}

\define\La{\Lambda}

\define\boc{\bold c}

\define\boj{\bold j}

\define\boo{\bold o}

\define\bos{\bold s}

\define\bc{\bold C}

\define\be{\bold E}

\define\bn{\bold N}

\define\bq{\bold Q}
\define\br{\bold R}

\define\bz{\bold Z}

\define\ca{\Cal A}

\define\cd{\Cal D}
\define\ce{\Cal E}

\define\cg{\Cal G}
\define\ch{\Cal H}
\define\ci{\Cal I}

\define\cl{\Cal L}
\define\cm{\Cal M}

\define\cs{\Cal S}

\define\fA{\frak A}

\define\fS{\frak S}

\define\opl{\oplus}
\define\sha{\sharp}
\redefine\aa{\bold a}
\define\Mod{\text{\rm Mod}}
\define\Irr{\text{\rm Irr}}

\define\Inv{\text{\rm Inv}}  
\define\Sy{\text{\rm Sy}}
\define\BO{B}
\define\BR{Br}
\define\HO{Ho}
\define\KL{KL}
\define\LSp{L1}
\define\LU{L2}
\define\LC{L3}
\define\LCC{L4}
\define\LL{L5}

These are notes for lectures given at MIT during the Fall of 1999.

\head 1. Coxeter groups\endhead
\subhead 1.1\endsubhead
Let $S$ be a finite 
set and let $M=(m_{s,s'})_{(s,s')\in S\tim S}$ be a matrix with 
entries in $\bn\cup\{\infty\}$ such that $m_{s,s}=1$ for all $s$ and
$m_{s,s'}=m_{s',s}\ge 2$ for all $s\ne s'$. (A {\it "Coxeter matrix"}.) Let $W$
be the group defined by the generators $s (s\in S)$ and relations

$(ss')^{m_{s,s'}}=1$
\nl
for any $s,s'$ in $S$ such that $m_{s,s'}<\infty$. (A {\it "Coxeter group"}.) 
In $W$ we have $s^2=1$ for all $s$. 

Clearly, there is a unique homomorphism $\sgn:W@>>>\{1,-1\}$ such that
$\sgn(s)=-1$ for all $s$. ("{\it Sign} representation".)

For $w\in W$ let $l(w)$ be the smallest integer $q\ge 0$ such that
$w=s_1s_2\dots s_q$ with $s_1,s_2,\dots, s_q$ in $S$. (We then say that 
$w=s_1s_2\dots s_q$ is a {\it reduced expression} and $l(w)$ is the 
{\it length} of $w$.) Note that $l(1)=0,l(s)=1$ for $s\in S$. (Indeed, $s\ne 1$
in $W$ since $\sgn(s)=-1,\sgn(1)=1$.)

\proclaim{Lemma 1.2} Let $w\in W,s\in S$.

(a) We have either $l(sw)=l(w)+1$ or $l(sw)=l(w)-1$.

(b) We have either $l(ws)=l(w)+1$ or $l(ws)=l(w)-1$.
\endproclaim
Clearly, $\sgn(w)=(-1)^{l(w)}$. Since $\sgn(sw)=-\sgn(w)$, we have 
$(-1)^{l(sw)}=-(-1)^{l(w)}$. Hence $l(sw)\ne l(w)$. This, together with the 
obvious inequalities $l(w)-1\le l(sw)\le l(w)+1$ gives (a). The proof of (b) is
similar.

\proclaim{Proposition 1.3} Let $E$ be an $\br$-vector space with basis 
$(e_s)_{s\in S}$. For $s\in S$ let $\si_s:E@>>>E$ be the linear map defined by 

$\si_s(e_{s'})=e_{s'}+2\cos\fra{\pi}{m_{s,s'}}e_s$ for all $s'\in S$.

(a) There is a unique homomorphism $\si:W@>>>GL(E)$ such that $\si(s)=\si_s$ 
for all $s\in S$.

(b) If $s\ne s'$ in $S$, then $ss'$ has order $m_{s,s'}$ in $W$. In particular,
$s\ne s'$ in $W$.
\endproclaim
We have $\si_s(e_s)=-e_s$ and $\si$ induces the identity map on $E/\br e_s$. It
follows that $\si_s^2=1$. Now let $s\ne s'$ in $S$. Let $m=m_{s,s'}$ and let
$\Phi=\si_s\si_{s'}$. We have

$\Phi(e_s)=(4\cos^2\fra{\pi}{m}-1)e_s+2\cos\fra{\pi}{m}e_{s'}$,

$\Phi(e_{s'})=-2\cos\fra{\pi}{m}e_s-e_{s'}$.
\nl
Hence $\Phi$ restricts to an endomorphism $\phi$ of $\br e_s\opl\br e_{s'}$ 
whose characteristic polynomial is

$X^2-2\cos\fra{2\pi}{m}X+1=(X-e^{2\pi\sqrt{-1}/m})(X-e^{-2\pi\sqrt{-1}/m})$.
\nl
It follows that, if $2<m<\infty$, then $1+\phi+\phi^2+\dots+\phi^{m-1}=0$. The
same is true if $m=2$ (in this case we see directly that $\phi=-1$). Since 
$\Phi$ induces the identity map on $E/(\br e_s\opl\br e_{s'})$, it follows that
$\Phi:E@>>>E$ has order $m$ (if $m<\infty$). If $m=\infty$, we have $\phi\ne 1$
and $(\phi-1)^2=0$, hence $\phi$ has infinite order and $\Phi$ has also 
infinite order. Now both (a), (b) follow.

\proclaim{Corollary 1.4} Let $s_1\ne s_2$ in $S$. Let $\lan s_1,s_2\rangle$ be
the subgroup of $W$ generated by $s_1,s_2$. For $k\ge 0$ let
$1_k=s_1s_2s_1\dots$ ($k$ factors), $2_k=s_2s_1s_2\dots$ ($k$ factors).

(a) Assume that $m=m_{s_1,s_2}<\infty$. Then $\lan s_1,s_2\rangle$ consists of 
the elements $1_k,2_k$ ($k=0,1,\dots,m$); these elements are distinct except 
for the equalities $1_0=2_0, 1_m=2_m$. For $k\in[0,m]$ we have 
$l(1_k)=l(2_k)=k$.

(b) Assume that $m_{s_1,s_2}=\infty$. Then $\lan s_1,s_2\rangle$ consists of 
the elements $1_k,2_k$ ($k=0,1,\dots$); these elements are distinct except for
the equality $1_0=2_0$. For all $k\ge 0$ we have $l(1_k)=l(2_k)=k$.
\endproclaim
This follows immediately from 1.3(b).

We identify $S$ with a subset of $W$ (see 1.3(b)). Let
$T=\cup_{w\in W}wSw\i\sub W$.

\proclaim{Proposition 1.5} Let $R=\{1,-1\}\tim T$. For $s\in S$ let 
$U_s:R@>>>R$ be the map defined by $U_s(\ep,t)=(\ep(-1)^{\de_{s,t}},sts)$ where
$\de$ is the Kronecker symbol. There is a unique homomorphism $U$ of $W$ into
the group of permutations of $R$ such that $U(s)=U_s$ for all $s\in S$.
\endproclaim
We have $U_s^2(\ep,t)=(\ep(-1)^{\de_{s,t}+\de_{s,sts}},t)=(\ep,t)$ since the
conditions $s=t$, $s=sts$ are equivalent. Thus, $U_s^2=1$. For $s\ne s'$ in
$S$ with $m=m_{s,s'}<\infty$ we have

$U_sU_{s'}(\ep,t)=(\ep(-1)^{\de_{s',t}+\de_{s,s'ts'}},ss'ts's)$
\nl
hence
$$\align&(U_sU_{s'})^m(\ep,t)=
(\ep(-1)^{\de_{s',t}+\de_{s,s'ts'}+\de_{s',ss'ts's}+\de_{s,s'ss'ts'ss'}+\dots},
(ss')^mt(s's)^m)\\&
=(\ep(-1)^{\de_{s',t}+\de_{s'ss',t}+\de_{s'ss'ss',t}+\dots},t)\endalign$$
(both sums have exactly $2m$ terms). It is enough to show that 
$\de_{s',t}+\de_{s'ss',t}+\de_{s'ss'ss',t}+\dots$ is even, or that $t$ appears
an even number of times in the $2m$-term sequence $s',s'ss',s'ss'ss',\dots$.
This follows from the fact that in this sequence the $k$-th terms is equal to 
the $(k+m)$-th term for $k=1,2,\dots,m$.

\proclaim{Proposition 1.6} Let $w\in W$. 

(a) If $w=s_1s_2\dots s_q$ is a reduced expression, then the elements

$s_1,s_1s_2s_1,s_1s_2s_3s_2s_1,\dots,s_1s_2\dots s_q\dots s_2s_1$ 
\nl 
are distinct. 

(b) These elements form a subset of $T$ that depends only on $w$, not on the 
choice of reduced expression for it.
\endproclaim
Assume that
$s_1s_2\dots s_i\dots s_2s_1=s_1s_2\dots s_j\dots s_2s_1$ for some
$1\le i<j\le q$. Then $s_i=s_{i+1}s_{i+2}\dots s_j\dots s_{i+2}s_{i+1}$ hence
$$\align s_1s_2\dots s_q&=s_1s_2\dots s_{i-1}(s_{i+1}s_{i+2}\dots s_j\dots 
s_{i+2}s_{i+1})s_{i+1}\dots s_js_{j+1}\dots s_q\\&=
s_1s_2\dots s_{i-1}s_{i+1}s_{i+2}\dots s_{j-1}s_{j+1}\dots s_q,\endalign$$
which shows that $l(w)\le q-2$, contradiction. This proves (a). 

For $(\ep,t)\in R$ we have (see 1.5)
$U(w\i)(\ep,t)=(\ep\eta(w,t),w\i tw)$ where $\eta(w,t)=\pm 1$ depends only on 
$w,t$. On the other hand, 
$$\align U(w\i)(\ep,t)&=U_{s_q}\dots U_{s_1}(\ep,t)\\&=(\ep(-1)^{\de_{s_1,t}+
\de_{s_2,s_1ts_1}+\dots+\de_{s_q,s_{q-1}\dots s_1ts_1\dots s_{q-1}}},w\i tw)\\&
=(\ep(-1)^{\de_{s_1,t}+\de_{s_1s_2s_1,t}+\dots+
\de_{s_1\dots s_q\dots s_1,t}},w\i tw).\endalign$$
Thus, $\eta(w,t)=(-1)^{\de_{s_1,t}+\de_{s_1s_2s_1,t}+\dots+
\de_{s_1\dots s_q\dots s_1,t}}$. Using (a), we see that for $t\in T$, the sum 
$\de_{s_1,t}+\de_{s_1s_2s_1,t}+\dots+\de_{s_1\dots s_q\dots s_1,t}$ is $1$ if
$t$ belongs to the subset in (b) and is $0$, otherwise. Hence the subset in (b)
is just $\{t\in T|\eta(w,t)=-1\}$. This completes the proof.

\proclaim{Proposition 1.7} Let $w\in W,s\in S$ be such that $l(sw)=l(w)-1$. Let
$w=s_1s_2\dots s_q$ be a reduced expression. Then there exists $j\in[1,q]$ such
that 

$ss_1s_2\dots s_{j-1}=s_1s_2\dots s_j$.
\endproclaim
Let $w'=sw$. Let $w'=s'_1s'_2\dots s'_{q-1}$ be a reduced expression. Then
$w=ss'_1s'_2\dots s'_{q-1}$ is another reduced expression. By 1.6(b), the 
$q$-term sequences 

$s_1,s_1s_2s_1,s_1s_2s_3s_2s_1,\dots$ and $s,ss'_1s,ss'_1s'_2s'_1s,\dots$ 
\nl
coincide up to rearranging terms. In particular, 
$s=s_1s_2\dots s_j\dots s_2s_1$ for some $j\in[1,q]$. The proposition follows.

\subhead 1.8 \endsubhead
Let $X$ be the set of all sequences $(s_1,s_2,\dots,s_q)$ in $W$ such that
$s_1s_2\dots s_q$ is a reduced expression in $W$. We regard $X$ as the vertices
of a graph in which $(s_1,s_2,\dots,s_q),(s'_1,s'_2,\dots,s'_{q'})$ are joined
if one is obtained from the other by replacing $m$ consecutive entries of form
$s,s',s,s',\dots$ by the $m$ entries $s',s,s',s,\dots$; here $s\ne s'$ in $S$ 
are such that $m=m_{s,s'}<\infty$. We use the notation 

$(s_1,s_2,\dots,s_q)\sim(s'_1,s'_2,\dots,s'_{q'})$
\nl
for
"$(s_1,s_2,\dots,s_q),(s'_1,s'_2,\dots,s'_{q'})$ are in the same connected
component of $X$".
\nl
(When this holds we have necessarily $q=q'$ and
$s_1s_2\dots s_q=s'_1s'_2\dots s'_q$ in $W$.)

\proclaim{Theorem 1.9} Let $(s_1,s_2,\dots,s_q),(s'_1,s'_2,\dots,s'_q)$ in $X$
be such that $s_1s_2\dots s_q=s'_1s'_2\dots s'_q=w\in W$. Then 
$(s_1,s_2,\dots,s_q)\sim(s'_1,s'_2,\dots,s'_q)$.
\endproclaim
We shall use the following notation. 

Let $\bos=(s_1,s_2,\dots,s_q),\bos'=(s'_1,s'_2,\dots,s'_q)$. Let $C$ (resp. 
$C'$) be the connected component of $X$ that contains $\bos$ (resp. $\bos'$). 
For $i\in[1,q]$ we set

$\bos(i)=(\dots,s'_1,s_1,s'_1,s'_1,s_1,s_2,s_3,\dots,s_i)$ (a $q$-elements
sequence in $S$),
 
$\und{\bos}(i)=\dots s'_1s_1s'_1s'_1s_1s_2s_3\dots s_i\in W$ (the product of 
this sequence).
\nl
Let $C(i)$ be the connected component of $X$ that contains $\bos(i)$. Then 
$\bos=\bos(q)$. Hence $C=C(q)$.

We argue by induction on $q$. The theorem is obvious for $q\le 1$. We now 
assume that $q\ge 2$ and that the theorem is known for $q-1$. We first prove 
the following weaker statement.

$(A)$ {\it In the setup of the theorem we have either 

$(s_1,s_2,\dots,s_q)\sim(s'_1,s'_2,\dots,s'_q)$, or
\nl
(a) $s_1s_2\dots s_q=s'_1s_1s_2\dots s_{q-1}$ and
$(s'_1,s_1,s_2,\dots,s_{q-1})\sim (s'_1,s'_2,\dots,s'_q)$.}
\nl
We have $l(s'_1w)=l(w)-1$. By 1.7 we have 
$s'_1s_1s_2\dots s_{i-1}=s_1s_2\dots s_i$
for some $i\in[1,q]$, so that $w=s'_1s_1s_2\dots s_{i-1}s_{i+1}\dots s_q$. In 
particular, 

$(s'_1,s_1,s_2,\dots,s_{i-1},s_{i+1},\dots,s_q)\in X$. 
\nl
By the induction hypothesis, we have 
$(s_1,s_2,\dots,s_{i-1},s_{i+1},\dots,s_q)\sim(s'_2,\dots,s'_q)$. Hence 

(b) $(s'_1,s_1,s_2,\dots,s_{i-1},s_{i+1},\dots,s_q)\sim(s'_1,s'_2,\dots,s'_q)$.
\nl
Assume first that $i<q$. Then from 
$s'_1s_1s_2\dots s_{i-1}s_{i+1}\dots s_{q-1}=s_1s_2\dots s_{q-1}$ 
and the induction hypothesis we deduce that 

$(s'_1,s_1,s_2,\dots,s_{i-1},s_{i+1},\dots,s_{q-1})
\sim(s_1,s_2,\dots,s_{q-1})$, hence

$(s'_1,s_1,s_2,\dots,s_{i-1},s_{i+1},\dots,s_{q-1},s_q)\in C$.
\nl
Combining this with (b) we deduce that $C=C'$.

Assume next that $i=q$ so that $s_1s_2\dots s_q=s'_1s_1s_2\dots s_{q-1}$. Then
(b) shows that (a) holds. Thus, $(A)$ is proved.

Next we prove for $p\in[0,q-2]$ the following generalization of $(A)$.

$(A'_p)$ {\it In the setup of the theorem we have either $C=C'$ or:

for $i\in[q-p-1,q]$ we have $\bos(i)\in X,\und\bos(i)=w$, $C_i=C$ if 
$i=q\mod 2$ and $C_i=C'$ if $i=q+1\mod 2$.}
\nl
For $p=0$ this reduces to $(A)$. Assume now that $p>0$ and that $(A'_{p-1})$ is
already known. We prove that $(A'_p)$ holds.

If $C=C'$, then we are done. Hence by $(A'_{p-1})$ we may assume that:
for $i\in[q-p,q]$ we have $\bos(i)\in X,\und\bos(i)=w$, $C_i=C$ if $i=q\mod 2$
and $C_i=C'$ if $i=q+1\mod 2$.

Applying $(A)$ to $\bos(q-p),\bos(q-p+1)$ (instead of $\bos,\bos'$), we see 
that either $C_{q-p}=C_{q-p+1}$ or:

$\bos(q-p),\bos(q-p-1)$ are in $X$,$\und\bos(q-p)=\und\bos(q-p-1)$ and
$C_{q-p-1}=C_{q-p+1}$.
\nl
In both cases, we see that $(A'_p)$ holds.

This completes the inductive proof of $(A'_p)$. In particular, $(A'_{q-2})$ 
holds. In other words, in the setup of the theorem, either $C=C'$ holds or:

(c) {\it for $i\in[1,q]$ we have $\bos(i)\in X,\und\bos(i)=w$, $C_i=C$ if 
$i=q\mod 2$ and $C_i=C'$ if $i=q+1\mod 2$.}
\nl
If $C=C'$, then we are done. Hence we may assume that (c) holds. In particular,

(d) $\bos(2)\in X,\bos(1)\in X,\und\bos(2)=\und\bos(1)$. 
\nl
From $\bos(1)\in X$ and $q\ge 2$ we see that $s'\ne s$ and that 
$q\le m=m_{s,s'}$. From $\und\bos(2)=\und\bos(1)$ we see that 
$s_2\in\lan s_1,s'_1\rangle$, hence $s_2$ is either $s_1$ or $s'_1$. In fact we
cannot have $s_2=s_1$ since this would contradict $\bos(2)\in X$. Hence 
$s_2=s'_1$. We see that $\bos(2)=(\dots,s'_1,s_1,s'_1,s_1,s'_1)$ (the number of
terms is $q, q\le m$). Since $\und\bos(2)=\und\bos(1)$, it follows that $q=m$,
so that $\bos(2),\bos(1)$ are joined in $X$. It follows that $C_2=C_1$. By (c),
for some permutation $a,b$ of $1,2$ we have $C_a=C, C_b=C'$. Since $C_a=C_b$ it
follows that $C=C'$.It follows that $C=C'$. The theorem is proved.

\proclaim{Proposition 1.10} Let $w\in W$ and let $s,t\in S$ be such that 
$l(swt)=l(w), l(sw)=l(wt)$. Then $sw=wt$.
\endproclaim
Let $w=s_1s_2\dots s_q$ be a reduced expression.

Assume first that $l(wt)=q+1$. Then $s_1s_2\dots s_qt$ is a reduced expression
for $wt$. Now $l(swt)=l(wt)-1$ hence by 1.7 there exists $i\in[1,q]$ such that
$ss_1s_2\dots s_{i-1}=s_1s_2\dots s_i$ or else 
$ss_1s_2\dots s_q=s_1s_2\dots s_qt$. If the second alternative occurs, we are
done. If the first alternative occurs, we have
$sw=s_1s_2\dots s_{i-1}s_{i+1}\dots s_q$ hence $l(sw)\le q-1$. This contradicts
$l(sw)=l(wt)$.

Assume next that $l(wt)=q-1$. Let $w'=wt$. Then $l(sw't)=l(w'),l(sw')=l(w't)$.
We have $l(w't)=l(w')+1$ hence the first part of the proof applies and gives
$sw'=w't$. Hence $sw=wt$. The proposition is proved.

\subhead 1.11\endsubhead
We can regard $S$ as the set of vertices of a graph in which $s,s'$ are joined
if $m_{s,s'}>2$. We say that $W$ is {\it irreducible} if this graph is 
connected. It is easy to see that in general, $W$ is naturally a product of 
irreducible Coxeter groups, corresponding to the connected components of $S$.

In the setup of 1.3, let $(,):E\tim E\to\br$ be the symmetric $\br$-bilinear
form given by $(e_s,e_{s'})=-\cos\fra{\pi}{m_{s,s'}}$. Then $\si(w):E@>>>E$
preserves $(,)$ for any $w\in W$. We say that $W$ is {\it tame} if $(e,e)\ge 0$
for any $e\in E$. It is easy to see that, if $W$ is finite then $W$ is tame. 

We say that $W$ is {\it integral} if, for any $s\ne s'$ in $S$, we have
$4\cos^2\fra{\pi}{m_{s,s'}}\in\bz$ (or equivalently 
$m_{s,s'}\in\{2,3,4,6,\infty\}$).

We will be mainly interested in the case where $W$ is tame. The tame, 
irreducible $W$ are of three kinds:

(a) finite, integral;

(b) finite, non-integral;

(c) tame, infinite (and automatically integral).
\nl
The $W$ of type (c) are called {\it affine Weyl groups}.

\head 2. Partial order on $W$\endhead
\subhead 2.1\endsubhead
Let $y,w$ be two elements of $W$. We say that $y\le w$ if there exists a 
sequence $y=y_0,y_1,y_2,\dots,y_n=w$ such that $l(y_k)-l(y_{k-1})=1$ for 
$k\in[1,n]$ and $y_ky_{k-1}\i\in T$ (or equivalently $y_{k-1}y_k\i\in T$, or 
$y_k\i y_{k-1}\in T$, or $y_{k-1}\i y_k\in T$) for $k\in[1,n]$.

This is clearly a partial order on $W$. Note that $y\le w$ implies 
$l(y)\le l(w)$. Also, $y\le w$ implies $y\i\le w\i$. If $w\in W,s\in S$ then, 
clearly:

$sw<w$ if and only if $l(sw)=l(w)-1$;

$sw>w$ if and only if $l(sw)=l(w)+1$.

\proclaim{Lemma 2.2} Let $w=s_1s_2\dots s_q$ be a reduced expression and let 
$t\in T$. The following are equivalent:

(i) $U(w\i)(\ep,t)=(-\ep,w\i tw)$ for $\ep=\pm 1$;

(ii) $t=s_1s_2\dots s_i\dots s_2s_1$ for some $i\in[1,q]$;

(iii) $l(tw)<l(w)$.
\endproclaim
The equivalence of (i),(ii) has been proved earlier.

Proof of (ii)$\implies$(iii). Assume that (ii) holds. Then 
$tw=s_1\dots s_{i-1}s_{i+1}\dots s_q$ hence $l(tw)<q$ and (iii) holds.

Proof of (iii)$\implies$(i). First we check that 

(a) $U(t)(\ep,t)=(-\ep,t)$. 
\nl
If $t\in S$, (a) is clear. If (a) is true for $t$ then it is also true for 
$sts$ where $s\in S$. Indeed, 
$$\align U(sts)(\ep,sts)&=U_sU(t)U_s(\ep,sts)=U_sU(t)(\ep(-1)^{\de_{s,sts}},t)
=U_s(-\ep(-1)^{\de_{s,sts}},t)\\&=(-\ep(-1)^{\de_{s,sts}+\de_{s,t}},t)=
(-\ep,t);\endalign$$
(a) follows. Assume now that (i) does not hold; thus, 
$U(w\i)(\ep,t)=(\ep,w\i tw)$. Then
$$\align U((tw)\i)(\ep,t)&=U(w\i)U(t)(\ep,t)=U(w\i)(-\ep,t)=(-\ep,w\i tw)\\&=
(-\ep,(tw)\i t(tw)).\endalign$$
Since (i)$\implies$(iii) we deduce that $l(w)<l(tw)$; thus, (iii) does not 
hold. The lemma is proved.

\proclaim{Lemma 2.3} Let $y,z\in W$ and let $s\in S$. If $sy\le z<sz$, then
$y\le sz$.
\endproclaim
We argue by induction on $l(z)-l(sy)$. If $l(z)-l(sy)=0$ then $z=sy$ and the
result is clear. Now assume that $l(z)>l(sy)$. Then $sy<z$. We can assume that
$sy<y$ (otherwise the result is trivial). We can find $t\in T$ such that 
$sy<tsy\le z$ and $l(tsy)=l(sy)+1$. If $t=s$, then $y\le z$ and we are done. 
Hence we may assume that $t\ne s$. We show that 

(a) $y<stsy$.
\nl
Assume that (a) does not hold. Then $y,tsy,sy,stsy$ have lengths $q+1,q+1,q,q$.
We can find a reduced expression $y=ss_1s_2\dots s_q$. Since $l(stsy)<l(y)$, we
see from 2.2 that either $sts=ss_1\dots s_i\dots s_1s$ for some $i\in[1,q]$ or
$sts=s$. (This last case has been excluded.) It follows that

$tsy=s_1\dots s_i\dots s_1sss_1s_2\dots s_q=s_1\dots s_{i-1}s_{i+1}\dots s_q$.
\nl
Thus, $l(tsy)\le q-1$, a contradiction. Thus, (a) holds. Let $y'=stsy$. We have
$sy'\le z\le sz$ and $l(z)-l(sy')<l(z)-l(sy)$. By the induction hypothesis, we
have $y'\le sz$. We have $y<y'$ by (a), hence $y\le sz$. The lemma is proved.

\proclaim{Proposition 2.4} The following three conditions on $y,w\in W$ are
equivalent:

(i) $y\le w$;

(ii) for any reduced expression $w=s_1s_2\dots s_q$ there exists a subsequence
$i_1<i_2<\dots<i_r$ of $1,2,\dots,q$ such that $y=s_{i_1}s_{i_2}\dots s_{i_r}$,
$r=l(y)$;

(iii) there exists a reduced expression $w=s_1s_2\dots s_q$ and a subsequence
$i_1<i_2<\dots<i_r$ of $1,2,\dots,q$ such that $y=s_{i_1}s_{i_2}\dots s_{i_r}$.
\endproclaim
Proof of (i)$\implies$(ii). We may assume that $y<w$. Let
$y=y_0,y_1,y_2,\dots,y_n=w$ be as in 2.1.  Let $w=s_1s_2\dots s_q$ be a reduced
expression. Since $y_{n-1}y_n\i\in T$, $l(y_{n-1})=l(y_n)-1$, we see from 2.2 
that there exists $i\in[1,q]$ such that 
$y_{n-1}y_n\i=s_1s_2\dots s_i\dots s_2s_1$ hence
$y_{n-1}=s_1s_2\dots s_{i-1}s_{i+1}\dots s_q$. This is a reduced expression. 
Similarly, since $y_{n-2}y_{n-1}\i\in T$, $l(y_{n-2})=l(y_{n-1})-1$, we see 
from 2.2 (applied to $y_{n-1}$) that there exists $j\in[1,q]-\{i\}$ such that 
$y_{n-2}$ equals

$s_1s_2\dots s_{i-1}s_{i+1}\dots s_{j-1}s_{j+1}\dots s_q$ or
$s_1s_2\dots s_{j-1}s_{j+1}\dots s_{i-1}s_{i+1}\dots s_q$
\nl
(depending on whether $i<j$ or $i>j$). Continuing in this way we see that $y$ 
is of the required form.

Proof of (ii)$\implies$(iii). This is trivial.

Proof of (iii)$\implies$(i). Assume that $w=s_1s_2\dots s_q$ (reduced 
expression) and $y=s_{i_1}s_{i_2}\dots s_{i_r}$ where $i_1<i_2<\dots<i_r$ is a
subsequence of $1,2,\dots,q$. We argue by induction on $q$. If $q=0$ there is
nothing to prove. Now assume $q>0$. 

If $i_1>1$, then the induction hypothesis is applicable to $y,w'=s_2\dots s_q$
and yields $y\le w'$. But $w'\le w$ hence $y\le w$. If $i_1=1$ then the 
induction hypothesis is applicable to $y'=s_{i_2}\dots s_{i_r},w'=s_2\dots s_q$
and yields $y'\le w'$. Thus, $s_1y\le s_1w<w$. By 2.3 we then have $y\le w$.
The proposition is proved.

\proclaim{Corollary 2.5} Let $y,z\in W$ and let $s\in S$. 

(a) Assume that $sz<z$. Then $y\le z\lra sy\le z$.

(b) Assume that $y<sy$. Then $y\le z \lra y\le sz$.
\endproclaim
We prove (a). We can find a reduced expression of $z$ of form
$z=ss_1s_2\dots s_q$. Assume that $y\le z$. By 2.4 we can find a subsequence
$i_1<i_2<\dots<i_r$ of $1,2,\dots,q$ such that either
$y=s_{i_1}s_{i_2}\dots s_{i_r}$ or $y=ss_{i_1}s_{i_2}\dots s_{i_r}$. In the 
first case we have $sy=ss_{i_1}s_{i_2}\dots s_{i_r}$ and in the second case we
have $sy=s_{i_1}s_{i_2}\dots s_{i_r}$. In both cases we have $sy\le z$ by 2.4.
The same argument shows that, if $sy\le z$ then $y\le z$. This proves (a).

We prove (b). Assume that $y\le z$. We must prove that $y\le sz$. If $z<sz$,
this is clear. Thus we may assume that $sz<z$. We can find a reduced expression
of $z$ of form $z=ss_1s_2\dots s_q$. By 2.4 we can find a subsequence
$i_1<i_2<\dots<i_r$ of $1,2,\dots,q$ such that either
$y=s_{i_1}s_{i_2}\dots s_{i_r}, l(y)=r$ or 
$y=ss_{i_1}s_{i_2}\dots s_{i_r},l(y)=r+1$. In the second case we have 
$l(sy)=r<l(y)$, contradicting $y<sy$. Thus we are in the first case. Hence $y$
is the product of a subsequence of $s_1,s_2,\dots,s_q$ and using again 2.4, we
deduce that $y\le sz$ (note that $sz=s_1s_2\dots s_q$ is a reduced expression).
The lemma is proved.

\head 3. The algebra $\ch$\endhead
\subhead 3.1\endsubhead
A map $L:W@>>>\bz$ is said to be a {\it weight function} for $W$ if 
$L(ww')=L(w)+L(w')$ for any $w,w'\in W$ such that $l(ww')=l(w)+l(w')$. We will
assume that a weight function $L:W@>>>\bz$ is fixed; we then say that $W,L$ is
a {\it weighted Coxeter group}. (For example we could take $L=l$; in that case 
we say that we are in the {\it split case}.) Note that $L$ is determined by its
values $L(s)$ on $S$ which are subject only to the condition that $L(s)=L(s')$
for any $s\ne s'$ in $S$ such that $m_{s,s'}$ is finite and odd. We necessarily
have $L(1)=0$ and $L(w)=L(w\i)$ for all $w\in W$.

Let $\ca=\bz[v,v\i]$ where $v$ is an indeterminate. For $s\in S$ we set 
$v_s=v^{L(s)}\in\ca$.

\subhead 3.2\endsubhead
Let $\ch$ be the $\ca$-algebra with $1$ defined by the generators 
$T_s (s\in S)$ and the relations

(a) $(T_s-v_s)(T_s+v_s\i)=0$ for $s\in S$;

(b) $T_sT_{s'}T_s\dots=T_{s'}T_sT_{s'}\dots$ 
\nl
(both products have $m_{s,s'}$ factors) for any $s\ne s'$ in $S$ such that 
$m_{s,s'}<\infty$.  

$\ch$ is called the {\it Hecke algebra} or the {\it Iwahori-Hecke algebra}. 

For $w\in W$ we define $T_w\in\ch$ by $T_w=T_{s_1}T_{s_2}\dots T_{s_q}$, where
$w=s_1s_2\dots s_q$ is a reduced expression. By (b) and 1.9, $T_w$ is 
independent of the choice of reduced expression. From the definitions it is 
clear that for $s\in S,w\in W$ we have

$T_sT_w=T_{sw}$ if $l(sw)=l(w)+1$,

$T_sT_w=T_{sw}+(v_s-v_s\i)T_w$ if $l(sw)=l(w)-1$.
\nl
In particular, the $\ca$-submodule of $\ch$ generated by $\{T_w|w\in W\}$ is
a left ideal of $\ch$. It contains $1=T_1$ hence it is the whole of $\ch$.
Thus $\{T_w|w\in W\}$ generates the $\ca$-module $\ch$.

\proclaim{Proposition 3.3} $\{T_w|w\in W\}$ is an $\ca$-basis of $\ch$.
\endproclaim
We consider the free $\ca$-module $\ce$ with basis $(e_w)_{w\in W}$. For any 
$s\in S$ we define $\ca$-linear maps $P_s:\ce\to\ce,Q_s:\ce\to\ce$ by

$P_s(e_w)=e_{sw}$ if $l(sw)=l(w)+1$,

$P_s(e_w)=e_{sw}+(v_s-v_s\i)e_w$ if $l(sw)=l(w)-1$;

$Q_s(e_w)=e_{ws}$ if $l(ws)=l(w)+1$,

$Q_s(e_w)=e_{ws}+(v_s-v_s\i)e_w$ if $l(ws)=l(w)-1$.
\nl
We shall continue the proof assuming that

(a) $P_sQ_t=Q_tP_s$ {\it for any $s,t$ in} $S$.
\nl
Let $\fA$ be the $\ca$-subalgebra with $1$ of $\End(\ce)$ generated by 
$\{P_s|s\in S\}$. The map $\fA@>>>\ce$ given by $\pi\mto\pi(e_1)$ is 
surjective. Indeed, if $w=s_1s_2\dots s_q$ is a reduced expression, then
$e_w=P_{s_1}\dots P_{s_q}e_1$. Assume now that $\pi\in\fA$ satisfies 
$\pi(e_1)=0$. Let $\pi'=Q_{s_q}\dots Q_{s_1}$. By (a) we have $\pi\pi'=\pi'\pi$
hence

$0=\pi'\pi(e_1)=\pi\pi'(e_1)=\pi(Q_{s_q}\dots Q_{s_1}(e_1))=\pi(e_w)$.
\nl
Since $w$ is arbitrary, it follows that $\pi=0$. We see that the map 
$\fA@>>>\ce$ is injective, hence an isomorphism of $\ca$-modules. Using this 
isomorphism we transport the algebra structure of $\fA$ to an algebra structure
on $\ce$ with unit element $e_1$. For this algebra structure we have 
$P_s(e_1)\pi(e_1)=P_s(\pi(e_1))$ for $s\in S,\pi\in\fA$. Hence 
$e_se_w=P_s(e_w)$ for any $w\in W,s\in S$. It follows that

(b) $e_se_w=e_{sw}$ if $l(sw)=l(w)+1$, 

(c) $e_se_w=e_{sw}+(v_s-v_s\i)e_w$ if $l(sw)=l(w)-1$.
\nl
From (b) it follows that, if $w=s_1s_2\dots s_q$ is a reduced expression, then
$e_w=e_{s_1}e_{s_2}\dots e_{s_q}$. In particular, if $s\ne s'$ in $S$ are such
that $m=m_{s,s'}<\infty$ then $e_se_{s'}e_s\dots=e_{s'}e_se_{s'}\dots$ (both
products have $m$ factors); indeed, this follows from the equality
$e_{ss's\dots}=e_{s'ss'\dots}$ (see 1.4). From (c) we deduce that 
$e_s^2=1+(v_s-v_s\i)e_s$ for $s\in S$, or that $(e_s-v_s)(e_s+v_s\i)=0$. We see
that there is a unique algebra homomorphism $\ch@>>>\ce$ preserving $1$ such 
that $T_s\mto e_s$ for all $s\in S$. This homomorphism takes $T_w$ to 
$e_w$ for any $w\in W$. Assume now that $a_w\in\ca$ ($w\in W$) are zero for all
but finitely many $w$ and that $\sum_wa_wT_w=0$ in $\ch$. Applying 
$\ch@>>>\ce$ we obtain $\sum_wa_we_w=0$. Since $(e_w)$ is a basis of $\ce$, it
follows that $a_w=0$ for all $w$. Thus, $\{T_w|w\in W\}$ is an $\ca$-basis 
of $\ch$. This completes the proof, modulo the verification of (a).

We prove (a). Let $w\in W$. We distinguish six cases. 

{\it Case} 1. $swt,sw,wt,w$ have lengths $q+2,q+1,q+1,q$. Then 

$P_sQ_t(e_w)=Q_tP_s(e_w)=e_{swt}$.
\nl
{\it Case} 2. $w,sw,wt,swt$ have lengths $q+2,q+1,q+1,q$. Then 
$$\align&P_sQ_t(e_w)=Q_tP_s(e_w)\\&=e_{swt}+(v_t-v_t\i)e_{sw}+(v_s-v_s\i)
e_{wt}+(v_t-v_t\i)(v_s-v_s\i)e_w.\endalign$$
{\it Case} 3. $wt,swt,w,sw$ have lengths $q+2,q+1,q+1,q$. Then 

$P_sQ_t(e_w)=Q_tP_s(e_w)=e_{swt}+(v_s-v_s\i)e_{wt}$.
\nl
{\it Case} 4. $sw,swt,w,wt$ have lengths $q+2,q+1,q+1,q$. Then 

$P_sQ_t(e_w)=Q_tP_s(e_w)=e_{swt}+(v_t-v_t\i)e_{sw}$.
\nl
{\it Case} 5. $swt,w,wt,sw$ have lengths $q+1,q+1,q,q$. Then 

$P_sQ_t(e_w)=e_{swt}+(v_t-v_t\i)e_{sw}+(v_t-v_t\i)(v_s-v_s\i)e_w$,

$Q_tP_s(e_w)=e_{swt}+(v_s-v_s\i)e_{wt}+(v_t-v_t\i)(v_s-v_s\i)e_w$.
\nl
{\it Case} 6. $sw,wt,w,swt$ have lengths $q+1,q+1,q,q$. Then 

$P_sQ_t(e_w)=e_{swt}+(v_s-v_s\i)e_{wt}$,

$Q_tP_s(e_w)=e_{swt}+(v_t-v_t\i)e_{sw}$.
\nl
In case 5 we have $L(t)+L(wt)=L(w)=L(swt)=L(s)+L(wt)$ hence $L(t)=L(s)$
and $v_s=v_t$. In case 6 we have $L(t)+L(swt)=L(sw)=L(wt)=L(s)+L(swt)$, hence 
$L(t)=L(s)$ and $v_s=v_t$. In case 5 and 6 we have $sw=wt$ by 1.10. Hence 
$P_sQ_t(e_w)=Q_tP_s(e_w)$ in each case. The proposition is proved.

\subhead 3.4\endsubhead
There is a unique involutive antiautomorphism of the algebra $\ch$ which 
carries $T_s$ to $T_s$ for any $s\in S$. (This follows easily by looking at the
defining relations of $\ch$.) It carries $T_w$ to $T_{w\i}$ for any $w\in W$. 

\subhead 3.5\endsubhead
There is a unique algebra involution of $\ch$ denoted $h\mto h^\dag$ such that 
$T_s^\dag=-T_s\i$ for any $s\in S$. We have $T_w^\dag=\sgn(w)T_{w\i}\i$ for any
$w\in W$. 

\head 4. The bar operator \endhead
\subhead 4.1\endsubhead
For $s\in S$, the element $T_s\in\ch$ is invertible: we have
$T_s\i=T_s-(v_s-v_s\i)$. It follows that $T_w$ is invertible for each $w\in W$;
if $w=s_1s_2\dots s_q$ is a reduced expression, then
$T_w\i=T_{s_q}\i\dots T_{s_2}\i T_{s_1}\i$.

Let $\bar{}:\ca@>>>\ca$ be the ring involution which takes $v^n$ to $v^{-n}$
for any $n\in\bz$.

\proclaim{Lemma 4.2}(a) There is a unique ring homomorphism $\bar{}:\ch@>>>\ch$
which is $\ca$-semilinear with respect to $\bar{}:\ca@>>>\ca$ and satisfies
$\ov{T}_s=T_s\i$ for all $s\in S$.

(b) This homomorphism is involutive. It takes $T_w$ to $T_{w\i}\i$ for
any $w\in W$.
\endproclaim
The following two identities can be deduced easily from 3.2(a),(b):

$(T_s\i-v_s\i)(T_s\i+v_s)=0$ for $s\in S$,

$T_s\i T_{s'}\i T_s\i\dots=T_{s'}\i T_s\i T_{s'}\i\dots$ 
\nl
(both products have $m_{s,s'}$ factors) for any $s\ne s'$ in $S$ such that 
$m_{s,s'}<\infty$; (a) follows.

We prove (b).  Let $s\in S$. Applying $\bar{}$ to $T_s\ov{T_s}=1$ gives
$\bar{T_s}\bar{\bar{T_s}}=1$. We have also $\ov{T_s}T_s=1$ hence 
$\bar{\bar{T_s}}=T_s$. It follows that the square of $\bar{}$ is $1$. The
second assertion of (b) is immediate. The lemma is proved.

\subhead 4.3\endsubhead
For any $w\in W$ we can write uniquely

$\ov{T}_w=\sum_{y\in W}\ov{r}_{y,w}T_y$
\nl
where $r_{y,w}\in\ca$ are zero for all but finitely many $y$.

\proclaim{Lemma 4.4} Let $w\in W$ and $s\in S$ be such that $w>sw$. For 
$y\in W$ we have

$r_{y,w}=r_{sy,sw}$ if $sy<y$,

$r_{y,w}=r_{sy,sw}+(v_s-v_s\i)r_{y,sw}$ if $sy>y$.
\endproclaim
We have 
$$\align\ov{T}_w&=T_s\i\ov{T}_{sw}=(T_s-(v_s-v_s\i))\sum_y\bar r_{y,sw}T_y\\&=
\sum_y\bar r_{y,sw}T_{sy}-\sum_y(v_s-v_s\i)\bar r_{y,sw}T_y+\sum_{y;sy<y}
(v_s-v_s\i)\bar r_{y,sw}T_y\\&=\sum_y\bar r_{sy,sw}T_y-\sum_{y;sy>y}(v_s-v_s\i)
\bar r_{y,sw}T_y.\endalign$$
The lemma follows.

\proclaim{Lemma 4.5} For any $y,w$ we have 
$\bar r_{y,w}=\sgn(yw)r_{y,w}$.
\endproclaim
We argue by induction on $l(w)$. If $w=1$ the result is obvious. Assume now 
that $l(w)\ge 1$. We can find $s\in S$ such that $w>sw$. Assume first that 
$sy<y$. From 4.4 we see, using the induction hypothesis, that

$\bar r_{y,w}=\bar r_{sy,sw}=\sgn(sysw)r_{sy,sw}=\sgn(yw)r_{y,w}$.
\nl
Assume next that $sy>y$. From 4.4 we see, using the induction hypothesis, that
$$\align\bar r_{y,w}&=\bar r_{sy,sw}+(v_s\i-v_s)\bar r_{y,sw}
=\sgn(sysw)r_{sy,sw}+(v_s\i-v_s)\sgn(ysw)r_{y,sw}\\&=\sgn(yw)
(r_{sy,sw}+(v_s-v_s\i)r_{y,sw})=\sgn(yw)r_{y,w}.\endalign$$
The lemma is proved.

\proclaim{Lemma 4.6} For any $x,z\in W$ we have 
$\sum_y\ov{r}_{x,y}r_{y,z}=\de_{x,z}$.
\endproclaim
Using the fact that $\bar{}$ is an involution, we have

$T_z=\ov{\ov{T}}_z=\ov{\sum_y\ov{r}_{y,z}T_y}=
\sum_yr_{y,z}\ov{T}_y=\sum_y\sum_xr_{y,z}\ov{r}_{x,y}T_x$.
\nl
We now compare the coefficients of $T_x$ on both sides. The lemma follows.

\proclaim{Proposition 4.7} Let $y,w\in W$.

(a) If $r_{y,w}\ne 0$, then $y\le w$.

(b)  Assume that $L(s)>0$ for all $s\in S$. If $y\le w$, then 

$r_{y,w}=v^{L(w)-L(y)} \mod v^{L(w)-L(y)-1}\bz[v\i]$,

$r_{y,w}=\sgn(yw)v^{-L(w)+L(y)} \mod v^{-L(w)+L(y)+1}\bz[v]$.

(c) Without assumption on $L$, $r_{y,w}\in v^{L(w)-L(y)}\bz[v^2,v^{-2}]$.
\endproclaim
We prove (a) by induction on $l(w)$. If $w=1$ the result is obvious. Assume now
that $l(w)\ge 1$. We can find $s\in S$ such that $w>sw$. Assume first that 
$sy<y$. From 4.4 we see that $r_{sy,sw}\ne 0$ hence, by the induction 
hypothesis, $sy\le sw$. Thus $sy\le sw<w$ and, by 2.3, we deduce $y\le w$. 
Assume next that $sy>y$. From 4.4 we see that either $r_{sy,sw}\ne 0$ or 
$r_{y,sw}\ne 0$ hence, by the induction hypothesis, $sy\le sw$ or $y\le sw$.
Combining this with $y<sy$ and $sw<w$ we see that $y\le w$. This proves (a).

We prove the first assertion of (b) by induction on $l(w)$. If $w=1$ the result
is obvious. Assume now that $l(w)\ge 1$. We can find $s\in S$ such that $w>sw$.
Assume first that $sy<y$. Then we have also $sy<w$ and, using 2.5(b), we deduce
$sy\le sw$. By the induction hypothesis, we have
$$\align r_{sy,sw}&=v^{L(sw)-L(sy)}+\text{strictly lower powers}\\&
=v^{L(w)-L(y)}+\text{strictly lower powers}.\endalign$$
But $r_{y,w}=r_{sy,sw}$ and the result follows. Assume next that $sy>y$. From 
$y<sy,y\le w$ we deduce using 2.5(b) that $y\le sw$. By the induction 
hypothesis, we have

$r_{y,sw}=v^{L(sw)-L(y)}+\text{strictly lower powers}$.
\nl
Hence
$$\align&(v_s-v_s\i)r_{y,sw}=v^{L(s)}v^{L(sw)-L(y)}+\text{strictly lower 
powers }\\&=v^{L(w)-L(y)}+\text{strictly lower powers}.\endalign$$
On the other hand, if $sy\le sw$, then by the induction hypothesis,
$$\align&r_{sy,sw}=v^{L(sw)-L(sy)}+\text{strictly lower powers}\\&
=v^{L(w)-L(y)-2L(s)}+\text{strictly lower powers}\endalign$$
while if $sy\not\le sw$ then $r_{sy,sw}=0$ by (a). Thus, in 

$r_{y,w}=r_{sy,sw}+(v_s-v_s\i)r_{y,sw}$,
\nl
the term $r_{sy,sw}$ contributes only powers of $v$ which are strictly smaller
than $L(w)-L(y)$ and thus, 
$r_{y,w}=v^{L(w)-L(y)}+\text{strictly lower powers}$. This proves the first
assertion of (b). The second assertion of (b) follows from the first using 4.5.

We prove (c) by induction on $l(w)$. If $w=1$ the result is obvious. Assume now
that $l(w)\ge 1$. We can find $s\in S$ such that $w>sw$. Assume first that 
$sy<y$. By the induction hypothesis, we have

$r_{y,w}=
r_{sy,sw}\in v^{L(sw)-L(sy)}\bz[v^2,v^{-2}]=v^{L(w)-L(y)}\bz[v^2,v^{-2}]$
\nl
as required. Assume next that $sy>y$. By the induction hypothesis, we have
$$\align&r_{y,w}=r_{sy,sw}+(v_s-v_s\i)r_{y,sw}\\&\in v^{L(sw)-L(sy)}
\bz[v^2,v^{-2}]+v^{L(s)}v^{L(sw)-L(y)}\bz[v^2,v^{-2}]=v^{L(w)-L(y)}
\bz[v^2,v^{-2}],\endalign$$
as required. The proposition is proved. 

\proclaim{Proposition 4.8} For any $x<z$ in $W$ we have
$\sum_{y;x\le y\le z}\sgn(y)=0$.
\endproclaim
Using 4.5 we can rewrite 4.6 (in our case) in the form

(a) $\sum_y\sgn(xy)r_{x,y}r_{y,z}=0$.
\nl
Here we may restrict the summation to $y$ such that $x\le y\le z$. In the rest 
of the proof we shall take $L=l$. Then 4.7(b) holds and we see that if 
$x\le y\le z$, then 

$r_{x,y}r_{y,z}$ is $v^{l(y)-l(x)}v^{l(z)-l(y)}+$ strictly lower powers of $v$.
\nl 
Hence (a) states that

$\sum_{y;x\le y\le z}\sgn(xy)v^{l(z)-l(x)}+$ strictly lower powers of $v$ is 
$0$.
\nl
In particular $\sum_{y;x\le y\le z}\sgn(xy)=0$. The proposition is proved.

\subhead 4.9\endsubhead
The involution $\bar{}:\ch@>>>\ch$ commutes with the involution in 3.4. (This
is clear on the generators of $\ch$.) It follows that

(a) $r_{y\i,w\i}=r_{y,w}$
\nl
for any $y,w\in W$.

On the other hand, it is clear that $\bar{}:\ch@>>>\ch$ and 
${}^\dag:\ch@>>>\ch$ commute.

\head 5. The elements $c_w$\endhead
\subhead 5.1\endsubhead
For any $n\in\bz$ let 

$\ca_{\le n}=\opl_{m;m\le n}\bz v^m, \ca_{\ge n}=\opl_{m;m\ge n}\bz v^m$,
$\ca_{<n}=\opl_{m;m<n}\bz v^m, \ca_{>n}=\opl_{m;m>n}\bz v^m$.
\nl
Note that $\ca_{\le 0}=\bz[v\i]$.

Let $\ch_{\le 0}=\opl_w\ca_{\le 0}T_w$, $\ch_{<0}=\opl_w\ca_{<0}T_w$. We have 
$\ch_{<0}\sub\ch_{\le 0}\sub\ch$. 

\proclaim{Theorem 5.2} (a) Let $w\in W$. There exists a unique element 
$c_w\in\ch_{\le 0}$ such that $\ov{c}_w=c_w$ and $c_w=T_w\mod\ch_{<0}$.

(b) $\{c_w|w\in W\}$ is an $\ca_{\le 0}$-basis of $\ch_{\le 0}$ and an 
$\ca$-basis of $\ch$.
\endproclaim
We prove the existence part of (a). We will construct, for any $x$ such that
$x\le w$, an element $u_x\in\ca_{\le 0}$ such that

(c) $u_w=1$,

(d) $u_x\in \ca_{<0}, \bar u_x-u_x=\sum_{y; x<y\le w}r_{x,y}u_y$
for any $x<w$.
\nl
We argue by induction on $l(w)-l(x)$. If $l(w)-l(x)=0$ then $x=w$ and we define
$u_x$ by (c). Assume now that $l(w)-l(x)>0$ and that $u_z$ is already defined
whenever $z\le w, l(w)-l(z)<l(w)-l(x)$ so that (c) holds and (d) holds if $x$ 
is replaced by any such $z$. Then the right hand side of the equality in (d) is
defined. We denote it by $a_x\in\ca$. We have
$$\align  a_x+\bar a_x &=\sum_{y;x<y\le w}r_{x,y}u_y+\sum_{y;x<y\le w}
\bar r_{x,y}\bar u_y\\&=\sum_{y;x<y\le w}r_{x,y}u_y+\sum_{y;x<y\le w}
\bar r_{x,y}(u_y+\sum_{z;y<z\le w}r_{y,z}u_z)\\&=\sum_{z;z<y\le w}r_{z,y}u_y
+\sum_{z;x<z\le w}\bar r_{x,z}u_z+\sum_{z;x<z\le w}\sum_{y; x<y<z}\bar r_{x,y}
r_{y,z}u_z\\&=\sum_{z;x<z\le w}\sum_{y; x\le y\le z}\bar r_{x,y}r_{y,z}u_z
=\sum_{z;x<z\le w}\de_{x,z}u_z=0.\endalign$$
(We have used 4.6 and the equality $r_{y,y}=1$.) Since $a_x+\bar a_x=0$, we 
have $a_x=\sum_{n\in\bz}\ga_nv^n$ (finite sum) where $\ga_n\in\bz$ satisfy 
$\ga_n+\ga_{-n}=0$ for all $n$ and in particular, $\ga_0=0$. Then 
$u_x=\sum_{n<0}\ga_nv^n\in \ca_{<0}$ satisfies $\bar u_x-u_x=a_x$. This 
completes the inductive construction of the elements $u_x$. We now define 
$c_w=\sum_{y;y\le w}u_yT_y\in\ch_{\le 0}$. It is clear that 
$c_w=T_w\mod\ch_{<0}$.
We have
$$\align&\ov{c}_w=\sum_{y;y\le w}\bar u_y\ov{T}_y=\sum_{y;y\le w}\bar u_y
\sum_{x;x\le y}\bar r_{x,y}T_x=\sum_{x;x\le w}(\sum_{y;x\le y\le w}
\bar r_{x,y}\bar u_y)T_x\\&=\sum_{x;x\le w}u_xT_x=c_w.\endalign$$
(We have used the fact that $r_{x,y}\ne 0$ implies $x\le y$, see 4.7, and (d).)
Thus, the existence of the element $c_w$ is established.

To prove uniqueness, it suffices to verify the following statement:

(e) {\it If $h\in\ch_{<0}$ satisfies $\bar h=h$ then} $h=0$.
\nl
We can write uniquely $h=\sum_{y\in W}f_yT_y$ where $f_y\in\ca_{<0}$ are zero
for all but finitely many $y$. Assume that not all $f_y$ are $0$. Then we can 
find $l_0\in\bn$ such that 

$Y_0=\{y\in W| f_y\ne 0, l(y)=l_0\}\ne\emptyset$ and
$\{y\in W| f_y\ne 0, l(y)>l_0\}=\emptyset$.
\nl
The equality $\sum_yf_yT_y=\ov{\sum_yf_yT_y}$ implies then

$\sum_{y\in Y_0}f_yT_y=\sum_{y\in Y_0}\bar f_yT_y\mod
\sum_{y; l(y)<l(y_0)}\ca T_y$
\nl
hence $\bar f_y=f_y$ for any $y\in Y_0$.  Since $f_y\in\ca_{<0}$, it follows 
that $f_y=0$ for any $y\in Y_0$, a contradiction. We have proved that $f_y=0$ 
for all $y$; (e) is verified and (a) is proved.

The elements $c_w$ constructed in (a) (for various $w$) are related to the
basis $T_w$ by a triangular matrix (with respect to $\le$) with $1$ on the
diagonal. Hence these elements satisfy (b). The theorem is proved.

\subhead 5.3\endsubhead
For any $w\in W$ we set $c_w=\sum_{y\in W}p_{y,w}T_y$ where 
$p_{y,w}\in\ca_{\le 0}$. By the proof of 5.2 we have 

$p_{y,w}=0$ unless $y\le w$, 

$p_{w,w}=1$,

$p_{y,w}\in\ca_{<0}$ if $y<w$.
\nl
Moreover, for any $x<w$ in $W$ we have

$\bar p_{x,w}=\sum_{y;x\le y\le w}r_{x,y}p_{y,w}$.

\proclaim{Proposition 5.4}(a) Assume that $L(s)>0$ for all $s\in S$. If 
$x\le w$, then 

$p_{x,w}=v^{-L(w)+L(x)} \mod v^{-L(w)+L(x)+1}\bz[v]$.

(b) Without assumption on $L$, for $x\le w$ we have
$p_{x,w}=v^{L(w)-L(x)}\bz[v^2,v^{-2}]$.
\endproclaim
We prove (a) by induction on $l(w)-l(x)$. If $l(w)-l(x)=0$ then $x=w$, 
$p_{x,w}=1$ and the result is obvious. Assume now that $l(w)-l(x)>0$. Using 
4.7(b) and the induction hypothesis, we see that
$\sum_{y;x<y\le w}r_{x,y}p_{y,w}$ is equal to
$$\sum_{y;x<y\le w}\sgn(x)\sgn(y)v^{-L(y)+L(x)}v^{-L(w)+L(y)}=
\sum_{y;x<y\le w}\sgn(x)\sgn(y)v^{-L(w)+L(x)}$$
plus strictly higher powers of $v$. Using 4.8, we see that this is 
$-v^{-L(w)+L(x)}$ plus strictly higher powers of $v$. Thus, 

$\bar p_{x,w}-p_{x,w}=-v^{-L(w)+L(x)}$ plus strictly higher powers of $v$.
\nl
Since $\bar p_{x,w}\in v\bz[v]$, it is in particular a $\bz$-linear combination
of powers of $v$ strictly higher than $-L(w)+L(y)$. Hence 

$-p_{x,w}=-v^{-L(w)+L(x)}$ plus strictly higher powers of $v$. 
\nl
This proves (a).

We prove (b) by induction on $l(w)-l(x)$. If $l(w)-l(x)=0$, then $x=w$, 
$p_{x,w}=1$ and the result is obvious. Assume now that $l(w)-l(x)>0$. Using 
4.7(c) and the induction hypothesis, we see that
$$\sum_{y;x<y\le w}r_{x,y}p_{y,w}\in\sum_{y;x<y\le w}v^{L(y)-L(x)}
v^{L(w)-L(y)}\bz[v^2,v^{-2}]=v^{L(w)-L(x)}\bz[v^2,v^{-2}].$$
Thus, $\bar p_{x,w}-p_{x,w}\in v^{L(w)-L(x)}\bz[v^2,v^{-2}]$. Hence
$p_{x,w}\in v^{L(w)-L(x)}\bz[v^2,v^{-2}]$. The proposition is proved.

\subhead 5.5\endsubhead
Let $s\in S$. From $T_s\i=T_s-(v_s-v_s\i)$ we see that
$r_{1,s}=v_s-v_s\i$. We also see that

$\ov{T_s+v_s\i}=T_s-(v_s-v_s\i)+v_s=T_s+v_s\i$,

$\ov{T_s-v_s}=T_s-(v_s-v_s\i)-v_s\i=T_s-v_s$.
\nl
If $L(s)=0$ we have $T_s\i=T_s$. Hence, 

$c_s=T_s+v_s\i$ if $L(s)>0$, 

$c_s=T_s-v_s$ if $L(s)<0$, 

$c_s=T_s$ if $L(s)=0$.

\subhead 5.6\endsubhead
The involution in 3.4 carries $T_w$ to $T_{w\i}$ hence it carries $\ch_{\le 0}$
into itself; moreover, it commutes with $\bar{}:\ch@>>>\ch$ (as pointed out in 
4.9). Hence it carries $c_w$ to $c_{w\i}$ for any $w\in W$. It follows that

(a) $p_{y\i,w\i}=p_{y,w}$
\nl
for any $y,w\in W$.

\head 6. Left or right multiplication by $c_s$\endhead
\subhead 6.1\endsubhead 
In this section we fix $s\in S$. Assume first that $L(s)=0$. In this case we 
have $c_s=T_s$; moreover, $T_sT_y=T_{sy}$. Hence for $w\in W$ we have

$c_sc_w=\sum_yp_{y,w}T_sT_y=\sum_yp_{y,w}T_{sy}=\sum_yp_{sy,w}T_y$.
\nl
We see that $c_sc_w\in\ch_{\le 0}$ and $c_sc_w=T_{sw}\mod\ch_{<0}$. Since
$\ov{c_sc_w}=c_sc_w$, it follows that, in this case, $c_sc_w=c_{sw}$. Similarly
we have $c_wc_s=c_{ws}$.

\subhead 6.2\endsubhead 
{\it In the remainder of this section (except in 6.8) we assume that} $L(s)>0$.

\proclaim{Proposition 6.3} To any $y,w\in W$ such that $sy<y<w<sw$ one can 
assign uniquely an element $\mu_{y,w}^s\in\ca$ so that

(i) $\ov{\mu}_{y,w}^s=\mu_{y,w}^s$ and

(ii) $\sum_{z;y\le z<w;sz<z}p_{y,z}\mu_{z,w}^s-v_sp_{y,w}\in\ca_{<0}$
\nl
for any $y,w\in W$ such that $sy<y<w<sw$.
\endproclaim
Let $y,w$ be as above. We may assume that $\mu_{z,w}^s$ are already defined for
all $z$ such that $y<z<w;sz<z$. Then condition (ii) is of the form:

$\mu_{y,w}^s$ {\it equals a known element of $\ca$ modulo} $\ca_{<0}$.
\nl
This condition determines uniquely the coefficients of $v^n$ with $n\ge 0$ in
$\mu_{y,w}^s$. Then condition (i) determines uniquely the coefficients of $v^n$
with $n<0$ in $\mu_{y,w}^s$. The proposition is proved.

\proclaim{Proposition 6.4} Let $y,w\in W$ be such that $sy<y<w<sw$. Then 
$\mu_{y,w}^s$ is a $\bz$-linear combination of powers $v^n$ with
$-L(s)+1\le n\le L(s)-1$ and $n=L(w)-L(y)-L(s)\mod 2$.
\endproclaim
We may assume that this is already known for all $\mu_{z,w}^s$ with $z$ such 
that $y<z<w;sz<z$. Using 6.3(ii) and 5.4, we see that $\mu_{y,w}^s$ is a 
$\bz$-linear combination of powers $v^n$ such that, whenever $n\ge 0$, we have
$n\le\L(s)-1$ and $n=L(w)-L(y)-L(s)\mod 2$. Using now 6.3(i), we deduce the
remaining assertions of the proposition.

\proclaim{Corollary 6.5} Assume that $L(s)=1$. Let $y,w\in W$ be such that 
$sy<y<w<sw$. Then $\mu_{y,w}^s$ is an integer, equal to the coefficient of 
$v\i$ in $p_{y,w}$. In particular, it is $0$ unless $L(w)-L(y)$ is odd.
\endproclaim
In this case, the inequalities of 6.4 become $0\le n\le 0$. They imply $n=0$.
Thus, $\mu_{y,w}^s$ is an integer. Picking up the coefficient of $v^0$ in the 
two sides of 6.3(ii), we see that $\mu_{y,w}^s$ is equal to the coefficient
of $v\i$ in $p_{y,w}$. The last assertion follows from 5.4.

\proclaim{Theorem 6.6} Let $w\in W$.

(a) If $w<sw$, then $c_sc_w=c_{sw}+\sum_{z;sz<z<w}\mu_{z,w}^sc_z$.

(b) If $sw<w$, then $c_sc_w=(v_s+v_s\i)c_w$.
\endproclaim
Since $c_s=T_s+v_s\i$ (see 5.5), we see that (b) is equivalent to
$(T_s-v_s)c_w=0$, or to

(c) $p_{x,w}=v_s\i p_{sx,w}$
\nl
(where  $sw<w,x<sx$). We prove the theorem by induction on $l(w)$. If $w=1$, 
the result is obvious. Assume now that $l(w)\ge 1$ and that the result holds 
when $w$ is replaced by $w'$ with $l(w')<l(w)$. 

{\it Case} 1. Assume that $w<sw$. Using $c_s=T_s+v_s\i$, we see that the 
coefficient of $T_y$ in the left hand side minus the right hand side of (a)
is
$$f_y=v_s^\si p_{y,w}+p_{sy,w}-p_{y,sw}-\sum_{z;y\le z<w;sz<z}p_{y,z}
\mu_{z,w}^s$$
where $\si=1$ if $sy<y$ and $\si=-1$ if $sy>y$. We must show that $f_y=0$. We 
first show that 

(d) $f_y\in\ca_{<0}$. 
\nl
If $sy<y$ this follows from 6.3(ii). (The contribution of $p_{sy,w}-p_{y,sw}$ 
is in $\ca_{<0}$ if $sy\ne w$ and is $1-1=0$ if $sy=w$.) 

If $sy>y$ then, by (c) (applied to $z$ in the sum, instead of $w$), we have
$$\align f_y&=v_s\i p_{y,w}+p_{sy,w}-p_{y,sw}-\sum_{z;y\le z<w;sz<z}v_s\i 
p_{sy,z}\mu_{z,w}^s\\&=v_s\i f_{sy}+v_s\i p_{sy,sw}-p_{y,sw}\endalign$$
(the second equality holds by 2.5(a)) and this is in $\ca_{<0}$ since 
$f_{sy}\in\ca_{<0}$ (by the previous paragraph), $v_s\i\in\ca_{<0}$ and since 
$y\ne sw$. Thus, (d) is proved.

Since both sides of (a) are fixed by $\bar{}$, the sum $\sum_yf_yT_y$ is fixed
by $\bar{}$. From (d) and 5.2(e) we see that $f_y=0$ for all $y$, as required.

{\it Case} 2. Assume that $w>sw$. Then case 1 is aplicable to $sw$ (by the 
induction hypothesis). We see that

$c_w=(T_s+v_s\i)c_{sw}-\sum_{z;sz<z<sw}\mu_{z,sw}^sc_z$.
\nl
Now $(T_s-v_s)(T_s+v_s\i)=0$ and $(T_s-v_s)c_z=0$ for each $z$ in the sum (by 
the induction hypothesis). Hence $(T_s-v_s)c_w=0$. The theorem is proved.

\proclaim{Corollary 6.7} Let $w\in W$.

(a) If $w<ws$, then $c_wc_s=c_{ws}+\sum_{z;zs<z<w}\mu_{z\i,w\i}^sc_z$.

(b) If $ws<w$, then $c_wc_s=(v_s+v_s\i)c_w$.
\endproclaim
We write the equalities in 6.6(a),(b) for $w\i$ instead of $w$ and we apply to
these equalities the involution 3.4. Since this involution carries $c_w$ to 
$c_{w\i}$, the corollary follows.

\subhead 6.8\endsubhead
6.3, 6.6, 6.7 remain valid when $L(s)<0$ provided that we replace in their 
statements and proofs $v_s$ by $-v_s\i$.

\head 7. Dihedral groups \endhead
\subhead 7.1\endsubhead
In this section we assume that $S$ consists of two elements $s_1,s_2$. For
$i=1,2$, let $L_i=L(s_i),T_i=T_{s_i},c_i=c_{s_i}$. We assume that 
$L_1>0,L_2>0$. Let $m=m_{s_1,s_2}$. Let $1_k,2_k$ be as in 1.4. For $w\in W$ 
we set

$\Ga_w=\sum_{y;y\le w}v^{-L(w)+L(y)}T_y$.

\proclaim{Lemma 7.2} We have

$c_1\Ga_{2_k}=\Ga_{1_{k+1}}+v^{L_1-L_2}\Ga_{1_{k-1}}$ if $k\in[2,m)$,

$c_2\Ga_{1_k}=\Ga_{2_{k+1}}+v^{-L_1+L_2}\Ga_{2_{k-1}}$ if $k\in[2,m)$,

$c_1\Ga_{2_k}=\Ga_{1_{k+1}}$ if $k=0,1$,

$c_2\Ga_{1_k}=\Ga_{2_{k+1}}$ if $k=0,1$.
\endproclaim
Since $c_i=T_i+v^{-L_i}$, the proof is an easy exercise.

\proclaim{Proposition 7.3} Assume that $L_1=L_2$. For any $w\in W$ we have 
$c_w=\Ga_w$.
\endproclaim
This is clear when $l(w)\le 1$. In the present case Lemma 7.2 gives
$$\Ga_{1_{k+1}}=c_1\Ga_{2_k}-\Ga_{1_{k-1}},\quad 
\Ga_{2_{k+1}}=c_2\Ga_{1_k}-\Ga_{2_{k-1}}\tag c$$
for $k\in[1,m)$. This shows by induction on $k$ that $\ov{\Ga_w}=\Ga_w$ for all
$w\in W$. Clearly, $\Ga_w=T_w\mod\ch_{<0}$. The lemma follows.

\subhead 7.4\endsubhead
{\it In 7.4-7.6 we assume that $L_2>L_1$.} In this case, if $m<\infty$, then 
$m$ is even. (See 3.1.) For $2k+1\in[1,m)$ we set 
$$\align&\Ga'_{2_{2k+1}}=\sum_{s\in[0,k-1]}(1-v^{2L_1}+v^{4L_1}-\dots+(-1)^s
v^{2sL_1})v^{-sL_1-sL_2}\\&\tim(T_{2_{2k-2s+1}}+v^{-L_2}T_{2_{2k-2s}}+v^{-L_2}
T_{1_{2k-2s}}+v^{-2L_2}T_{1_{2k-2s-1}})\\&+(1-v^{2L_1}+v^{4L_1}-\dots+(-1)^k
v^{2kL_1})v^{-kL_1-kL_2}(T_{2_1}+v^{-L_2}T_{2_0}).\endalign$$
For $2k+1\in[3,m)$ we set
$$\align&\Ga'_{1_{2k+1}}\\&=T_{1_{2k+1}}+v^{-L_1}T_{1_{2k}}+v^{-L_1}T_{2_{2k}}+
v^{-2L_1}T_{2_{2k-1}}+\sum\Sb y\\y\le 1_{2k-1}\endSb
v^{-L(w)+L(y)}(1+v^{2L_1})T_y.\endalign$$
For $w$ such that $l(w)$ is even and for $w=1_1$ we set $\Ga'_w=\Ga_w$.

\proclaim{Lemma 7.5} Let $\ze=v^{L_1-L_2}+v^{L_2-L_1}\in\ca$. We have

(a) $c_1\Ga'_{2_{k'}}=\Ga'_{1_{k'+1}}$, if $k'\in[0,m)$;

(b) $c_2\Ga'_{1_{k'}}=\Ga'_{2_{k'+1}}+\ze\Ga'_{2_{k'-1}}+\Ga'_{2_{k'-3}}$,
if $k'\in[4,m)$;

(c) $c_2\Ga'_{1_{k'}}=\Ga'_{2_{k'+1}}+\ze\Ga'_{2_{k'-1}}$, if $k'=2,3, k'<m$;

(d) $c_2\Ga'_{1_{k'}}=\Ga'_{2_{k'+1}}$ if $k'=0,1$.
\endproclaim
From the definitions we have

(e) $\Ga'_{2_{2k+1}}=\sum_{s\in[0,k]}(-1)^sv^{s(L_1-L_2)}\Ga_{2_{2k-2s+1}}$ if
$2k+1\in[1,m)$, 

(f) $\Ga'_{1_{2k+1}}=\Ga_{1_{2k+1}}+v^{L_1-L_2}\Ga_{1_{2k-1}}$ if 
$2k+1\in[3,m)$.
\nl
We prove (a) for $k'=2k+1$. The left hand side can be computed using (e) and 
7.2:
$$\align&c_1\Ga'_{2_{2k+1}}=c_1(\Ga_{2_{2k+1}}-v^{L_1-L_2}\Ga_{2_{2k-1}}+
v^{2L_1-2L_2}\Ga_{2_{2k-3}}+\dots)\\&=\Ga_{1_{2k+2}}+v^{L_1-L_2}\Ga_{1_{2k}}
-v^{L_1-L_2}\Ga_{1_{2k}}-v^{2L_1-2L_2}\Ga_{1_{2k-2}}\\&+v^{2L_1-2L_2}
\Ga_{1_{2k-2}}-v^{3L_1-3L_2}\Ga_{1_{2k-4}}+\dots=\Ga_{1_{2k+2}}=
\Ga'_{1_{2k+2}}.\endalign$$
This proves (a) for $k'=2k+1$. Now (a) for $k'=0$ is trivial. We prove (a) for
$k'=2k\ge 2$. The left hand side can be  computed using 7.2 and (f):

$c_1\Ga'_{2_{2k}}=c_1\Ga_{2_{2k}}=\Ga_{1_{2k+1}}+v^{L_1-L_2}\Ga_{1_{2k-1}}
=\Ga'_{1_{2k+1}}$.
\nl
This proves (a) for $k'=2k$. We prove (b) for $k'=2k$. The left hand side can 
be computed using 7.2:

$c_2\Ga'_{1_{2k}}=c_2\Ga_{1_{2k}}=\Ga_{2_{2k+1}}+v^{-L_1+L_2}\Ga_{2_{2k-1}}$.
\nl
The right hand side of (b) is (using (e)):
$$\align&\Ga_{2_{2k+1}}-v^{L_1-L_2}\Ga_{2_{2k-1}}+v^{2L_1-2L_2}\Ga_{2_{2k-3}}+
\dots\\&+\ze\Ga_{2_{2k-1}}-v^{L_1-L_2}\ze\Ga_{2_{2k-3}}+v^{2L_1-2L_2}\ze
\Ga_{2_{2k-5}}+\dots\\&+\Ga_{2_{2k-3}}-v^{L_1-L_2}\Ga_{2_{2k-5}}+v^{2L_1-2L_2}
\Ga_{2_{2k-7}}+\dots=\Ga_{2_{2k+1}}+v^{-L_1+L_2}\Ga_{2_{2k-1}}.\endalign$$
This proves (b) for $k'=2k$. We prove (b) for $k'=2k+1$. The left hand side can
be computed using (f) and 7.2:
$$\align& c_2\Ga'_{1_{2k+1}}=c_2(\Ga_{1_{2k+1}}+v^{L_1-L_2}\Ga_{1_{2k-1}})\\&
=\Ga_{2_{2k+2}}+v^{-L_1+L_2}\Ga_{2_{2k}}+v^{L_1-L_2}\Ga_{2_{2k}}+\Ga_{2_{2k-2}}
=\Ga'_{2_{2k+2}}+\ze\Ga'_{2_{2k}}+\Ga'_{2_{2k-2}}.
\endalign$$
This proves (b) for $k'=2k+1$. The proof of (c),(d) is similar to that of (b).
This completes the proof.

\proclaim{Proposition 7.6} For any $w\in W$ we have $c_w=\Ga'_w$.
\endproclaim
Clearly, $\Ga'_w=T_w\mod\ch_{<0}$. From the formulas in 7.5 we see by induction
on $l(w)$ that $\ov{\Ga'_w}=\Ga'_w$ for all $w$. The proposition is proved.

\proclaim{Proposition 7.7} Assume that $m=\infty$. For  $a\in\{1,2\}$, let
$f_a=v^{L(a)}+v^{-L(a)}$.

(a) Assume that $L_1=L_2$. For $k,k'\ge 0$, we have 

$c_{a_{2k+1}}c_{a_{2k'+1}}=f_a\sum_{u\in[0,\min(2k,2k')]}c_{a_{2k+2k'+1-2u}}$.

(b) Assume that $L_2>L_1$. For $k,k'\ge 0$, we have 

$c_{2_{2k+1}}c_{2_{2k'+1}}=f_2\sum_{u\in[0,\min(k,k')]}c_{2_{2k+2k'+1-4u}}$.

(c) Assume that $L_2>L_1$. For $k,k'\ge 1$, we have 

$$c_{1_{2k+1}}c_{1_{2k'+1}}=
f_1\sum_{u\in[0,\min(k-1,k'-1)]}p_uc_{1_{2k+2k'+1-2u}}$$
where $p_u=\ze$ for $u$ odd, $p_u\in\bz$ for $u$ even.
\endproclaim
We prove (a). For $k=k'=0$ the equality in (a) is clear. Assume now that 
$k=0,k'\ge 1$. Using 7.2, 7.3, we have
$$\align&c_2c_{2_{2k'+1}}=c_2(c_2c_{1_{2k'}}-c_{2_{2k'-1}})=
f_2c_2c_{1_{2k'}}-f_2c_{2_{2k'-1}}\\&=
f_2c_{2_{2k'+1}}+f_2c_{2_{2k'-1}}-f_2c_{2_{2k'-1}}=f_2c_{2_{2k'+1}},\endalign$$
as required. We now prove the equality in (a) for fixed $k'$, by induction on 
$k$. The case $k=0$ is already known. Assume now that $k=1$. From 7.2,7.3 we
have $c_{2_3}=c_2c_1c_2-c_2$. Using this and 7.2, 7.3, we have
$$\align&c_{2_3}c_{2_{2k'+1}}=c_2c_1c_2c_{2_{2k'+1}}-c_2c_{2_{2k'+1}}=f_2c_2
c_{1_{2k'+2}}+f_2c_2c_{1_{2k'}}-f_2c_{2_{2k'+1}}\\&=f_2c_{2_{2k'+3}}+f_2
c_{2_{2k'+1}}+f_2c_{2_{2k'+1}}+(1-\de_{k',0})f_2c_{2_{2k'-1}}-f_2c_{2_{2k'+1}}
\\&=f_2c_{2_{2k'+3}}+f_2c_{2_{2k'+1}}+f_2(1-\de_{k',0})c_{2_{2k'-1}},\endalign
$$
as required. Assume now that $k\ge 2$. From 7.2,7.3 we have

$c_{2_{2k+1}}=c_2c_1c_{2_{2k-1}}-2c_{2_{2k-1}}-c_{2_{2k-3}}$.
\nl
Using this and the induction hypothesis we have
$$\align&c_{2_{2k+1}}c_{2_{2k'+1}}=c_2c_1c_{2_{2k-1}}c_{2_{2k'+1}}
-2c_{2_{2k-1}}c_{2_{2k'+1}}-c_{2_{2k-3}}c_{2_{2k'+1}}\\&=f_2c_1c_2
\sum_{u\in[0,\min(2k-2,k')]}c_{2_{2k+2k'-1-2u}}-f_2\sum_{u\in[0,\min(2k-2,k')]}
c_{2_{2k+2k'-1-2u}}\\&-f_2\sum_{u\in[0,\min(2k-4,k')]}c_{2_{2k+2k'-3-2u}}.
\endalign$$
We now use 7.2,7.3 and (a) follows (for $a=2$). The case $a=1$ is similar.

We prove (b). For $k=k'=0$ the equality in (b) is clear. Assume now that 
$k=0,k'=1$. Using 7.5, 7.6, we have
$$c_2c_{2_3}=c_2(c_2c_{1_2}-\ze c_{2_1})=f_2c_2c_{1_2}-f_2\ze c_{2_1}
=f_2c_{2_3}+f_2\ze c_{2_1}-f_2\ze c_{2_1}=f_2c_{2_3},$$
as required. Assume next that $k=0,k'\ge 2$. Using 7.5, 7.6, we have
$$\align&c_2c_{2_{2k'+1}}=c_2(c_2c_{1_{2k'}}-\ze c_{2_{2k'-1}}-c_{2_{2k'-3}})
=f_2c_2c_{1_{2k'}}-f_2\ze c_{2_{2k'-1}}-f_2c_{2_{2k'-3}}\\&
=f_2c_{2_{2k'+1}}+f_2\ze c_{2_{2k'-1}}+f_2c_{2_{2k'-3}}
-f_2\ze c_{2_{2k'-1}}-f_2\ze c_{2_{2k'-3}}=f_2c_{2_{2k'+1}},\endalign$$
as required. We now prove the equality in (a) for fixed $k'$, by induction on 
$k$. The case $k=0$ is already known. Assume now that $k=1$. From 7.5,7.6 we
have $c_{2_3}=c_2c_1c_2-\ze c_2$. Using this and 7.5,7.6, we have
$$\align&c_{2_3}c_{2_{2k'+1}}=c_2c_1c_2c_{2_{2k'+1}}-\ze c_2c_{2_{2k'+1}}=
f_2c_2c_{1_{2k'+2}}-f_2\ze c_{2_{2k'+1}}\\&=f_2c_{2_{2k'+3}}+f_2\ze 
c_{2_{2k'+1}}+(1-\de_{k',0})f_2c_{2_{2k'-1}}-f_2\ze c_{2_{2k'+1}}\\&=f_2
c_{2_{2k'+3}}+(1-\de_{k',0})f_2c_{2_{2k'-1}}\endalign$$
as required. Assume now that $k\ge 2$. From 7.5,7.6 we have

$c_{2_{2k+1}}=c_2c_1c_{2_{2k-1}}-\ze c_{2_{2k-1}}-c_{2_{2k-3}}$.
\nl
Using this and the induction hypothesis we have
$$\align&c_{2_{2k+1}}c_{2_{2k'+1}}=c_2c_1c_{2_{2k-1}}c_{2_{2k'+1}}
-\ze c_{2_{2k-1}}c_{2_{2k'+1}}-c_{2_{2k-3}}c_{2_{2k'+1}}\\&
=f_2c_2c_1\sum_{u\in[0,\min(k-1,k')]}c_{2_{2k+2k'-1-4u}}
-f_2\ze\sum_{u\in[0,\min(k-1,k')]}c_{2_{2k+2k'-1-4u}}\\&
-f_2\sum_{u\in[0,\min(k-2,k')]}c_{2_{2k+2k'-3-4u}}.\endalign$$
We now use 7.5,7.6 and (b) follows.

The proof of (c) is similar to that of (b). This completes the proof.

\proclaim{Proposition 7.8} Assume that $4\le m<\infty$ and $L_2>L_1$. Then 
$m=2k+2$ with $k\ge 1$. Let

$p=(-1)^k(v^{L_2}+v^{-L_2})
(v^{k(L_2-L_1)}+v^{(k-2)(L_2-L_1)}+\dots+v^{-k(L_2-L_1)})$.
\nl
Then, for some $q\in\ca$, we have 

(a) $c_{2_{m-1}}c_{2_{m-1}}=pc_{2_{m-1}}+qc_{2_m}$.
\endproclaim
From 7.5,7.6, we see that $\ca c_{2_{m-1}}+\ca c_{2_m}$ is a two-sided ideal of
$\ch$. Hence (a) holds for some (unknown) $p,q\in\ca$. It remains to compute 
$p$. Let $\chi:\ch@>>>\ca$ be the algebra homomorphism defined by 
$\chi(T_1)=-v^{-L_1},\chi(T_2)=v^{L_2}$. Since $c_{2_m}=(T_1+v^{-L_1})h$ for 
some $h\in\ch$ (see 7.5,7.6) we see that $\chi(c_{2_m})=0$. Hence applying 
$\chi$ to (a) gives $\chi(c_{2_{m-1}})^2=p\chi(c_{2_{m-1}})$. It is thus enough
to show that 

(b) $\chi(c_{2_{m-1}})=p$. 
\nl
We verify (b) for $m=4$:
$$\align&\chi(T_2T_1T_2+v^{-L_2}T_2T_1+v^{-L_2}T_1T_2+v^{-2L_2}T_1\\&+
(v^{-L_1-L_2}-v^{L_1-L_2})T_2+(v^{-L_1-2L_2}-v^{L_1-2L_2}))\\&=-v^{-L_1+2L_2}
-2v^{-L_1}-v^{-L_1-2L_2}+(v^{-L_1-L_2}-v^{L_1-L_2})v^{L_2}\\&+v^{-L_1-2L_2}-
v^{L_1-2L_2}=-v^{-L_1+2L_2}-v^{-L_1}-v^{L_1}-v^{L_1-2L_2})\endalign$$
and for $m=6$:
$$\align&\chi(T_2T_1T_2T_1T_2+v^{-L_2}T_2T_1T_2T_1+v^{-L_2}T_1T_2T_1T_2+
v^{-2L_2}T_1T_2T_1\\&+(v^{-L_1-L_2}-v^{L_1-L_2})T_2T_1T_2+(v^{-L_1-2L_2}-
v^{L_1-2L_2})T_1T_2\\&+(v^{-L_1-2L_2}-v^{L_1-2L_2})T_2T_1+(v^{-L_1-3L_2}-
v^{L_1-3L_2})T_1\\&+(v^{-2L_1-2L_2}-v^{-2L_2}+v^{2L_1-2L_2})T_2
+(v^{-2L_1-3L_2}-v^{-3L_2}+v^{2L_1-3L_2}))\\&=v^{-2L_1+3L_2}+2v^{-2L_1+L_2}+
v^{-2L_1-L_2}-v^{-2L_1+L_2}-2v^{-2L_1-L_2}-v^{-2L_1-3L_2}+v^{L_2}\\&+2v^{-L_2}
+v^{-3L_2}+v^{-2L_1-L_2}-v^{-L_2}+v^{2L_1-L_2}+v^{-2L_1-3L_2}-v^{-3L_2}+
v^{2L_1-3L_2}\\&=v^{-2L_1+3L_2}+v^{-2L_1+L_2}+v^{L_2}+v^{-L_2}+v^{2L_1-L_2}+
v^{2L_1-3L_2}.\endalign$$
Analogous computations can be carried out for any even $m$. The proposition is
proved.

\head 8. Cells\endhead
\subhead 8.1\endsubhead
For $z\in W$ define $D_z\in\Hom_\ca(\ch,\ca)$ by

$D_z(c_w)=\de_{z,w}$ for all $w\in W$.
\nl
If $w,w'\in W$ we write $w\gets_\cl w'$ or $w'\to_\cl w$ if there exists 
$s\in S$ such that $D_w(c_sc_{w'})\ne 0$; we write $w\gets_{\Cal R} w'$ or 
$w'\to_{\Cal R}w$ if there exists $s\in S$ such that $D_w(c_{w'}c_s)\ne 0$;
we write $w\gets_{\cl\Cal R} w'$ or $w'\to_{\cl\Cal R}w$ if there exists 
$s\in S$ such that $D_w(c_sc_{w'})\ne 0$ or $D_w(c_{w'}c_s)\ne 0$.

If $w,w'\in W$, we say that $w\le_\cl w'$ (resp. $w\le_{\Cal R}w'$) if there 
exists a sequence $w=w_0,w_1,\dots,w_n=w'$ in $W$ such that

$w_0\gets_\cl w_1, w_1\gets_\cl w_2,\dots, w_{n-1}\gets_\cl w_n$ 

(resp. 
$w_0\gets_{\Cal R}w_1, w_1\gets_{\Cal R}w_2,\dots,w_{n-1}\gets_{\Cal R}w_n$).
\nl
If $w,w'\in W$, we say that $w\le_{\cl\Cal R}w'$ if there exists a sequence 

$w=w_0,w_1,\dots,w_n=w'$
\nl
in $W$ such that

$w_0\gets_{\cl\Cal R}w_1, w_1\gets_{\cl\Cal R}w_2,\dots,
w_{n-1}\gets_{\cl\Cal R}w_n$.
\nl
Clearly $\le_\cl,\le_{\Cal R},\le_{\cl\Cal R}$ are preorders on $W$. Let
$\sim_\cl,\sim_{\Cal R},\sim_{\cl\Cal R}$ be the associated equivalence
relations. (For example, we have $w\sim_\cl w'$ if and only if $w\le_\cl w'$
and $w'\le_\cl w$.) The equivalence classes on $W$ for
$\sim_\cl,\sim_{\Cal R},\sim_{\cl\Cal R}$ are called  respectively {\it left 
cells, right cells, two-sided cells} of $W$. They depend on $L:W@>>>\bn$.

If $w,w'\in W$, we say that $w<_\cl w'$ (resp. $w<_{\Cal R}w'$;
$w<_{\cl\Cal R}w'$) if $w\le_\cl w'$ and $w\not\sim_\cl w'$ (resp. 
$w<_{\Cal R}w'$ and $w\not\sim_{\Cal R}w'$; $w<_{\cl\Cal R}w'$ and 
$w\not\sim_{\cl\Cal R}w'$). 

Let $w,w'\in W$. It is clear that 

$w\le_\cl w'$ if and only if $w\i\le_{\Cal R}w'{}\i$,

$w\le_{\cl\Cal R}w'$ if and only if $w\i\le_{\cl\Cal R}w'{}\i$.
\nl
It follows that $w\mto w\i$ carries left cells to right cells, right cells to
left cells and two-sided cells to two-sided cells.

\proclaim{Lemma 8.2} Let $w'\in W$.

(a) $\ch_{\le_\cl w'}=\opl_{w;w\le_\cl w'}\ca c_w$ is a left ideal of $\ch$.

(b) $\ch_{\le_{\Cal R}w'}=\opl_{w;w\le_{\Cal R}w'}\ca c_w$ is a right ideal of
$\ch$.

(c) $\ch_{\le_{\cl\Cal R}w'}=\opl_{w;w\le_{\cl\Cal R}w'}\ca c_w$ is a 
two-sided ideal of $\ch$.
\endproclaim
We prove (a). Since $c_s (s\in S)$ generate $\ch$ as an $\ca$-algebra, it is
enough to verify the following statement:

{\it if $w\in W$ is such that $w\le_\cl w'$ and $s\in S$ then}
$c_sc_w\in\ch_{\le_\cl w'}$.
\nl
From the definition, $c_sc_w$ is an $\ca$-linear combination of elements
$c_{w''}$ with $w''\gets_\cl w$. For such $w''$ we clearly have 
$w''\le_\cl w'$. This proves (a). The proof of (b),(c) is entirely similar.

\subhead 8.3\endsubhead
Let $Y$ be a left cell of $W$. From 8.3(a) we see that for $y\in Y$,

$\opl_{w;w\le_\cl y}\ca c_w/\opl_{w;w<_\cl y}\ca c_w$
\nl
is a quotient of two left ideals of $\ch$ (independent of the choice of $y$)
hence it is naturally a left $\ch$-module; it has an $\ca$-basis consisting of
the images of $c_w (w\in Y)$.

Similarly, if $Y'$ is a right cell of $W$ then, for $y'\in Y'$,

$\opl_{w;w\le_{\Cal R}y'}\ca c_w/\opl_{w;w<_{\Cal R}y'}\ca c_w$
\nl
is a quotient of two right ideals of $\ch$ (independent of the choice of $y'$)
hence it is naturally a right $\ch$-module; it has an $\ca$-basis consisting of
the images of $c_w (w\in Y')$.

If $Y''$ is a two-sided cell of $W$ then, for $y''\in Y''$,

$\opl_{w;w\le_{\cl\Cal R}y''}\ca c_w/\opl_{w;w<_{\cl\Cal R}y''}\ca c_w$
\nl
is a quotient of two two-sided ideals of $\ch$ (independent of the choice of
$y''$) hence it is naturally a $\ch$-bimodule; it has an $\ca$-basis consisting
of the images of $c_w (w\in Y'')$.

\proclaim{Lemma 8.4} Let $s\in S$. Assume that $L(s)>0$. Let 
$\ch^s=\opl_{w;sw<w}\ca c_w$,${}^s\ch^s=\opl_{w;ws<w}\ca c_w$.

(a) $\{h\in\ch|(c_s-v_s-v_s\i)h=0\}=\ch^s$. Hence $\ch^s$ is a right ideal of
$\ch$.

(b) $\{h\in\ch|h(c_s-v_s-v_s\i)=0\}={}^s\ch$. Hence ${}^s\ch$ is a left ideal
of $\ch$.
\endproclaim
We prove the equality in (a). The right hand side is contained in the left hand
side by 6.6(b). Conversely, by 6.6, we have $c_sh\in\ch^s$ for any $h\in\ch$. 
Hence, if $h\in\ch$ is such that $c_sh=(v_s+v_s\i)h$, then 
$(v_s+v_s\i)h\in\ch^s$ so that $h\in\ch^s$ (since $\ch/\ch^s$ is a free 
$\ca$-module). This proves (a). The proof of (b) is entirely similar. The lemma
is proved.

\subhead 8.5\endsubhead
For $w\in W$ we set $\cl(w)=\{s\in S|sw<w\}, \Cal R(w)=\{s\in S|ws<w\}$.

\proclaim{Lemma 8.6} Let $w,w'\in W$. Assume that $L(s)>0$ for all $s\in S$.

(a) If $w\le_\cl w'$, then $\Cal R(w')\sub\Cal R(w)$. If $w\sim_\cl w'$, then 
$\Cal R(w')=\Cal R(w)$. 

(b) If $w\le_{\Cal R}w'$, then $\cl(w')\sub\cl(w)$. If $w\sim_{\Cal R}w'$, then
$\cl(w')=\cl(w)$.
\endproclaim
To prove the first assertion of (a), we may assume that $D_w(c_sc_{w'})\ne 0$ 
for some $s\in S$. Let $t\in\Cal R(w')$. We must prove that $t\in\Cal R(w)$. We
have $c_{w'}\in {}^t\ch$. By 8.4, ${}^t\ch$ is a left ideal of $\ch$. Hence 
$c_sc_{w'}\in{}^t\ch$. By the definition of ${}^t\ch$, for $h\in{}^t\ch$ we 
have $D_w(h)=0$ unless $wt<w$. Hence from $D_w(c_sc_{w'})\ne 0$ we deduce 
$wt<w$, as required. This proves the first assertion of (a). The second 
assertion of (a) follows immediately from the first. The proof of (b) is 
entirely similar to that of (a). The lemma is proved.

\subhead 8.7\endsubhead
We describe the left cells of $W$ in the setup of 7.3. From 7.2 and 7.3 we can
determine all pairs $y\ne w$ such that $y\gets w$ (we write $\gets,\to$ instead
of $\gets_\cl,\to_\cl$):

$1_0\to 2_1\lras 1_2\lras 2_3\lras \dots$,

$2_0\to 1_1\lras 2_2\lras 1_3\lras\dots$,
\nl
if $m=\infty$,

$1_0\to 2_1\lras 1_2\lras 2_3\lras \dots\lras 2_{m-1}\to 1_m$,

$2_0\to 1_1\lras 2_2\lras 1_3\lras\dots\lras 1_{m-1}\to 2_m$,
\nl
if $m<\infty, m$ even,

$1_0\to 2_1\lras 1_2\lras 2_3\lras\dots\lras 1_{m-1}\to 2_m$,

$2_0\to 1_1\lras 2_2\lras 1_3\lras\dots\lras 2_{m-1}\to 1_m$,
\nl
if $m<\infty, m$ odd. Hence the left cells are 

$\{1_0\}, \{2_1,1_2,2_3,\dots\}, \{1_1,2_2,1_3,\dots\}$,
\nl
if $m=\infty$, 

$\{1_0\},\{2_1,1_2,2_3,\dots,2_{m-1}\},\{1_1,2_2,1_3,\dots,1_{m-1}\},
\{2_m\}$,
\nl
if $m<\infty, m$ even,

$\{1_0\},\{2_1,1_2,2_3,\dots,1_{m-1}\},\{1_1,2_2,1_3,\dots,2_{m-1}\},
\{2_m\}$,
\nl
if $m<\infty, m$ odd.

The two-sided cells are $\{1_0\},W-\{1_0\}$ if $m=\infty$ and
$\{1_0\},\{2_m\},W-\{1_0,2_m\}$ if $m<\infty$.

\subhead 8.8\endsubhead
We describe the left cells of $W$ in the setup of 7.4. From 7.5 and 7.6 we can
determine all pairs $y\ne w$ such that $y\gets w$ (we write $\gets,\to$ instead
of $\gets_\cl,\to_\cl$). If $m=\infty$, these pairs are:

$1_0\to 2_1\lras 1_2\to 2_3\lras 1_4\to \dots,\quad
2_0\to 1_1\to 2_2\lras 1_3\to 2_4\lras\dots$,
\nl
and $2_1\gets 1_4, 2_2\gets 1_5, 2_3\gets 1_6,\dots$.

If $m=4$, these pairs are:

$1_0\to 2_1\lras 1_2\to 2_3\to 1_4,\quad 2_0\to 1_1\to 2_2\lras 1_3\to 2_4$.
\nl
If $m=6$, these pairs are:

$1_0\to 2_1\lras 1_2\to 2_3\lras 1_4\to 2_5\to 1_6,\quad
2_0\to 1_1\to 2_2\lras 1_3\to 2_4\lras 1_5\to 2_6$,
\nl
and $2_1\gets 1_4, 2_2\gets 1_5$. An analogous pattern holds for any even $m$.

Hence the left cells are 

$\{1_0\},\{2_1,1_2,2_3,1_4,\dots\},\{1_1\},\{2_2,1_3,2_4,1_5,\dots\}$,
\nl
if $m=\infty$,

$\{1_0\},\{2_1,1_2,2_3,\dots, 1_{m-2}\},\{2_{m-1}\},\{1_1\}, 
\{2_2,1_3,2_4,\dots,1_{m-1}\},\{2_m\}$,
\nl
if $m<\infty$.

The two-sided cells are 

$\{1_0\},\{1_1\}, W-\{1_0,1_1\}$, if $m=\infty$ and

$\{1_0\},\{1_1\},\{2_{m-1}\},\{2_m\},W-\{1_0,1_1,2_{m-1},2_m\}$, if $m<\infty$.

\head 9. Cosets of parabolic subgroups\endhead
\proclaim{Lemma 9.1}Assume that $w=s_1s_2\dots s_q$ with $s_i\in S$. We can 
find a subsequence $i_1<i_2<\dots<i_r$ of $1,2,\dots,r$ such that
$w=s_{i_1}s_{i_2}\dots s_{i_r}$ is a reduced expression.
\endproclaim
We argue by induction on $q$. If $q=0$ the result is obvious. Assume that 
$q>0$. Using the induction hypothesis we can assume that $s_2\dots s_q$ is a
reduced expression. If $s_1s_2\dots s_q$ is a reduced expression, we are done.
Hence we may assume that $s_1s_2\dots s_q$ is not a reduced expression. Then
$l(w)=q-1$. By 1.7, we can find $j\in [2,q]$ such that
$s_1s_2\dots s_{j-1}=s_2s_3\dots s_j$. Then
$w=s_2s_3\dots s_{j-1}s_{j+1}\dots s_q$ is a reduced expression. The 
lemma is proved.

\subhead 9.2\endsubhead
Let $w\in W$. Let $w=s_1s_2\dots s_q$ be a reduced expression of $w$. Using
1.9, we see that the set $\{s\in S| s=s_i \text{ for some } i\in[1,q]\}$ is 
independent of the choice of reduced expression. We denote it by $S_w$.

\subhead 9.3\endsubhead
In the remainder of this section we fix $I\sub S$. Let $W_I$ be the subgroup of
$W$ generated by $I$. 

If $w\in W_I$ then we can find a reduced expression $w=s_1s_2\dots s_q$ in $W$
with all $s_i\in I$ (we first write $w=s_1s_2\dots s_q$ a not necessarily 
reduced expression with all $s_i\in I$ and then we apply 9.1). Thus, 
$S_w\sub I$. Conversely, it is clear that if $w'\in W$ satisfies $S_{w'}\sub I$
then $w'\in W_I$. It follows that

$W_I=\{w\in W| S_w\sub I\}$.

\subhead 9.4\endsubhead
Replacing $S,(m_{s,s'})_{(s,s')\in S\tim S}$ by 
$I,(m_{s,s'})_{(s,s')\in I\tim I}$ in the definition of $W$ we obtain a Coxeter
group denoted by $W^*_I$. We have an obvious homomorphism $f:W^*_I@>>>W_I$ 
which takes $s$ to $s$ for $s\in I$.

\proclaim{Proposition 9.5}$f:W^*_I@>>>W_I$ is an isomorphism.
\endproclaim
We define $f':W_I@>>>W_I^*$ as follows: for $w\in W_I$ we choose a reduced
expression $w=s_1s_2\dots s_q$ in $W$; then $s_i\in I$ for all $i$ (see 9.3)
and we set $f'(w)=s_1s_2\dots s_q$ (product in $W_I^*$). This map is well 
defined. Indeed, if $s'_1s'_2\dots s'_q$ is another reduced expression for $w$
with all $s_i\in I$, then we can pass from $(s_1,s_2,\dots, s_q)$ to 
$(s'_1,s'_2,\dots,s'_q)$ by moving along edges of the graph $X$ (see 1.9); but
each edge involved in this move will necessarily involve only pairs $(s,s')$ in
$I$, hence the equation $s_1s_2\dots s_q=s'_1s'_2\dots s'_q$ must hold in 
$W_I^*$. It is clear that $ff'(w)=w$ for all $w\in W_I$. Hence $f'$ is 
injective. 

We show that $f'$ is a group homomorphism. It suffices to show that
$f'(sw)=f'(s)f'(w)$ for any $w\in W_I,s\in I$. This is clear if $l(sw)=l(w)+1$
(in $W$). Assume now that $l(sw)=l(w)-1$ (in $W$). Let $w=s_1s_2\dots s_r$ be a
reduced expression in $W$. Then $s_i\in I$ for all $i$. By 1.7 we have (in $W$)
$sw=s_1s_2\dots s_{i-1}s_{i+1}\dots s_q$ for some $i\in[1,q]$. Since 
$ss_1s_2\dots s_{i-1}s_{i+1}\dots s_q$ is a reduced expression for $w$ in $W$,
we have $f'(w)=ss_1s_2\dots s_{i-1}s_{i+1}\dots s_q$ (product in $W_I^*$). We 
also have $f'(w)=s_1s_2\dots s_q$ (product in $W_I^*$). Hence
$ss_1s_2\dots s_{i-1}s_{i+1}\dots s_q=s_1s_2\dots s_q$ (in $W_I^*$). Hence
$s_1s_2\dots s_{i-1}s_{i+1}\dots s_q=ss_1s_2\dots s_q$ (in $W_I^*$). Hence
$f'(sw)=f'(s)f'(w)$, as required.

Since the image of $f'$ contains the generators $s\in I$ of $W_I^*$ and $f'$ is
a group homomorphism, it follows that $f'$ is surjective. Hence $f'$ is 
bijective. Since $ff'=1$ it follows that $f$ is bijective. The proposition is
proved.

\subhead 9.6\endsubhead
We identify $W_I^*$ and $W_I$ via $f$. Thus, $W_I$ is naturally a Coxeter 
group. Let $l_I:W_I@>>>\bn$ be the length function of this Coxeter group. Let 
$w\in W_I$. Let $w=s_1s_2\dots s_q$ be a reduced expression of $w$ (in $W$).
Then $s_i\in I$ for all $i$ (see 9.3). Hence $l_I(w)\le l(w)$. The reverse 
inequality $l(w)\le l_I(w)$ is obvious. Hence $l_I(w)=l(w)$.

From 2.4 we see that the partial order on $W_I$ defined in the same way as 
$\le$ on $W$ is just the restriction of $\le$ from $W$ to $W_I$.

\proclaim{Lemma 9.7} Let $W_Ia$ be a coset in $W$. 

(a) This coset has a unique element $w$ of minimal length. 

(b) If $y\in W_I$ then $l(yw)=l(y)+l(w)$.

(c) $w$ is characterized by the property that $l(sw)>l(w)$ for all $s\in I$.
\endproclaim
Let $w$ be an element of minimal length in the coset. Let $w=s_1s_2\dots s_q$ 
be a reduced expression. Let $y\in W_I$ and let $y=s'_1s'_2\dots s'_p$ be a
reduced expression in $W_I$. Then $yw=s'_1s'_2\dots s'_ps_1s_2\dots s_q$. By 
9.1 we can drop some of the factors in the last product so that we are left
with a reduced expression for $yw$. The factors dropped cannot contain any
among the last $q$ since we would find an element in $W_Ia$ of strictly smaller
length than $w$. Thus, we can find a subsequence $i_1<i_2<\dots<i_r$ of 
$1,2,\dots,p$ such that $yw=s'_{i_1}s'_{i_2}\dots s'_{i_r}s_1s_2\dots s_q$ is a
reduced expression. It follows that 
$y=s'_1s'_2\dots s'_p=s'_{i_1}s'_{i_2}\dots s'_{i_r}$. Since $p=l(y)$, we must
have $r=p$ so that $s'_1s'_2\dots s'_ps_1s_2\dots s_q$ is a reduced expression
and $l(yw)=p+q=l(y)+l(w)$.

If now $w'$ is another element of minimal length in $W_Ia$ then $w'=yw$ for
some $y\in W_I$. We have $l(w)=l(w')=l(y)+l(w)$ hence $l(y)=0$ hence $y=1$ and
$w'=w$. This proves (a). Now (b) is already proved. Note that by (b), $w$ has 
the property in (c). Conversely, let $w'\in aW_I$ be an element such that
$l(sw')>l(w')$ for all $s\in I$. We have $w'=yw$ for some $y\in W_I$. If 
$y\ne 1$ then for some $s\in I$ we have $l(y)=l(sy)+1$. By (b) we have 
$l(w')=l(y)+l(w)$, $l(sw')=l(sy)+l(w)$. Thus $l(w')-l(sw')=l(y)-l(sy)=1$, a 
contradiction. Thus $y=1$ and $w'=w$. The lemma is proved.

\proclaim{Lemma 9.8} Let $W_Ia$ be a coset in $W$. 

(a) If $W_I$ is finite, this coset has a unique element $w$ of maximal length. 
If $W_I$ is infinite, this coset has no element of maximal length.

(b) Assume that $W_I$ is finite. If $y\in W_I$ then $l(yw)=l(w)-l(y)$.

(c) Assume that $W_I$ is finite. Then $w$ is characterized by the property that
$l(sw)<l(w)$ for all $s\in I$.
\endproclaim
Assume that $w$ has maximal length in $W_Ia$. We show that for any $y\in W_I$
we have 

(d) $l(yw)=l(w)-l(y)$. 
\nl
We argue by induction on $l(y)$. If $l(y)=0$, the result is clear. Assume now 
that $l(y)=p+1\ge 1$. Let $y=s_1\dots s_ps_{p+1}$ be a reduced expression. By 
the induction hypothesis, $l(w)=l(s_1s_2\dots s_pw)+p$. Hence we can find a 
reduced expression of $w$ of the form $s_p\dots s_2s_1s'_1s'_2\dots s'_q$. 
Since $s_{p+1}\in I$, by our assumption on $w$ we have $l(s_{p+1}w)=l(w)-1$. 
Using 1.7, we deduce that either

(1) $s_{p+1}s_p\dots s_{j+1}=s_p\dots  s_{j+1}s_j$ for some $j\in[1,p]$ or

(2) $s_{p+1}s_p\dots s_2s_1s'_1s'_2\dots s'_{i+1}= 
s_p\dots s_2s_1s'_1s'_2\dots s'_{i+1}s'_i$ for some $i\in[1,q]$.
\nl
In case (1) it follows that 
$y=s_1\dots s_ps_{p+1}=s_1s_2\dots s_{j-1}s_{j+1}\dots s_p$
\nl
contradicting $l(y)=p+1$. Thus, we must be in case (2). We have

$yw=s'_1s'_2\dots s'_{i-1}s'_{i+1}\dots s'_q$ and 
$l(yw)\le q-1=l(w)-p-1=l(w)-l(y)$. 
\nl
Thus, $l(w)\ge l(yw)+l(y)$. The reverse inequality is obvious. Hence 
$l(w)=l(yw)+l(y)$. This completes the induction.

From (d) we see that $l(y)\le l(w)$. Thus $l:W_I@>>>\bn$ is bounded above.
Hence there exists $y\in W_I$ of maximal length in $W_I$. Applying (d) to
$y,W_I$ instead of $w,W_Ia$ we see that 

$l(y)=l(y'{}\i)+l(y'{}\i y)=l(y')+l(y'{}\i y)$ 
\nl
for any $y'\in W_I$. Hence a reduced expression of $y'$ followed by a reduced
expression of $y'{}\i y$ gives a reduced expression of $y$. In particular
$y'\le y$. Since the set $\{y'\in W|y'\le y\}$ is finite, we see that $W_I$ is
finite. Conversely, if $W_I$ is finite then $W_Ia$ clearly has some element of
maximal length.

If $w'$ is another element of maximal length in $aW_I$ then $w'=yw$ for some 
$y\in W_I$. We have $l(w)=l(w')=l(w)-l(y)$ hence $l(y)=0$ hence $y=1$ and
$w'=w$. This proves (a) and (b). The proof of (c) is entirely similar to that 
of 9.7(c). The lemma is proved.

\subhead 9.9\endsubhead 
Replacing $W,L$ by $W_I,L|_{W_I}$ in the definition of $\ch$ we obtain an
$\ca$-algebra $\ch_I$ (naturally a subalgebra of $\ch)$; instead of
$r_{x,y},p_{x,y},c_y,\mu^s_{x,y}$ we obtain for $x,y\in W_I$ elements 
$r_{x,y}^I\in\ca,p_{x,y}^I\in\ca_{\le 0},c_y^I\in\ch_I,\mu^{s,I}_{x,y}\in\ca$.

\proclaim{Lemma 9.10} Let $z\in W$ be such that $z$ is the element of minimal 
length of $W_Iz$. Let $x,y\in W_I$. We have

(a) $\{u'\in W|xz\le u'\le yz\}=\{u\in W_I|x\le u\le y\}z$;

(b) $r_{xz,yz}=r^I_{x,y}$;

(c) $p_{xz,yz}=p^I_{x,y}$;

(d) $c^I_y=c_y$.

(e) If in addition, $s\in I$ and $sx<x<y<sy$, then $sxz<xz<yz<syz$ and
$\mu^{s,I}_{x,y}=\mu^s_{xz,yz}$.
\endproclaim
We first prove the following statement.

{\it Assume that $z_1,z_2$ have minimal length in $W_Iz_1,W_Iz_2$ 
respectively, that $u_1,u_2\in W_I$ and that $u_1z_1\le u_2z_2$. Then

(f) $z_1\le z_2$; if in addition, $z_1=z_2$ then $u_1\le u_2$. }
\nl
Indeed, using 2.4 we see that there exist $u'_1,z'_1$ such that

$u_1z_1=u'_1z'_1, u'_1\le u_2,z'_1\le z_2$.
\nl
Then $u'_1\in W_I$ and $z'_1\in W_Iz_1$ hence $z'_1=wz_1$ where $w\in W_I$, 
$l(z'_1)=l(w)+l(z_1)$. Hence $z_1\le z'_1$. Since $z'_1\le z_2$, we see that
$z_1\le z_2$. If we know that $z_1=z_2$, then $z'_1=z_1$ hence $u_1=u'_1$. 
Since $u'_1\le u_2$, it follows that $u_1\le u_2$ and (f) is proved.

We prove (a). If $u\in W_I$ and $x\le u\le y$, then $xz\le uz\le yz$ by 2.4 and
9.7(b). Conversely, assume that $u'\in W$ satisfies $xz\le u'\le yz$. Then 
$u'=uz_1$ where $z_1$ has minimal length in $W_Iu'$ and $u\in W_I$. Applying
(f) to $xz\le uz_1$ and to $uz_1\le yz$ we deduce $z\le z_1\le z$. Hence 
$z=z_1$. Applying the second part of (f) to $xz\le uz$ and to $uz\le yz$ we 
deduce $x\le u\le y$. This proves (a).

We prove (b) by induction on $l(y)$. Assume first that $y=1$. Then
$r^I_{x,y}=\de_{x,1}$. Now $r_{xz,z}=0$ unless $xz\le z$ (see 4.7(a)) in which 
case $x=1$ and $r_{z,z}=1$. Thus, (b) holds for $y=1$. Assume now that 
$l(y)\ge 1$. We can find $s\in I$ such that $l(sy)=l(y)-1$. We have 

$l(syz)=l(sy)+l(z)=l(y)-1+l(z)=l(yz)-1$.
\nl
If $sx<x$ then we have (as in the previous line) $sxz<xz$. Using 4.4 and the
induction hypothesis, we have

$r_{xz,yz}=r_{sxz,syz}=r^I_{sx,sy}=r^I_{x,y}$.
\nl
If $sx>x$ then we have (as above) $sxz>xz$. Using 4.4 and the induction 
hypothesis, we have

$r_{xz,yz}=r_{sxz,syz}+(v_s-v_s\i)r_{xz,syz}=
r^I_{sx,sy}+(v_s-v_s\i)r^I_{x,sy}=r^I_{x,y}$.
\nl
This completes the proof of (b).

We prove (c). Using (a), we may assume that $x\le y$ (otherwise, both sides are
zero.) We argue by induction on $l(y)-l(x)\ge 0$. If $y=x$, the result is clear
(both sides are $1$). Assume now that $l(y)-l(x)\ge 1$. Using 5.3, then (a),(b)
and the induction hypothesis, we have
$$\align &\bar p_{xz,yz}=\sum_{u';xz\le u'\le yz}r_{xz,u'}p_{u',yz}=\sum_{u\in
W_I;x\le u\le y}r_{xz,uz}p_{uz,yz}\\&=\sum_{u\in W_I;x\le u\le y}r^I_{x,u}
p_{uz,yz}=\sum_{u\in W_I;x<u\le y}r^I_{x,u}p^I_{u,y}+p_{xz,yz}.\endalign$$
Using 5.3 for $W_I$ we have
$\bar p^I_{x,y}=\sum_{y;x\le u\le y}r^I_{x,u}p^I_{u,y}$. Comparing with the 
previous equality we deduce

$\bar p_{xz,yz}-\bar p^I_{x,y}=p_{xz,yz}-p^I_{x,y}$.
\nl
The right hand side of this equality is in $\ca_{<0}$. Since it is fixed by
$\bar{}$, it must be $0$. This proves (c). Now (d) is an immediate consequence
of (c) (with $z=1$). 

We prove (e). By 6.3(ii) we have 

$\sum_{u';xz\le u'<yz;su'<u'}p_{xz,u'}\mu_{u',yz}^s-v_sp_{xz,yz}\in\ca_{<0}$.
\nl
We rewrite this using (a): 

$\sum_{u\in W_I;x\le u<y;su<u}p_{xz,uz}\mu_{uz,yz}^s-v_sp_{xz,yz}\in\ca_{<0}$.
\nl
We may assume that for all $u$ in the sum, other than $u=x$ we have
$\mu_{uz,yz}^s=\mu_{u,y}^{s,I}$. Using this and (d), we obtain

$\mu^s_{xz,yz}+\sum_{u\in W_I;x<u<y;su<u}p^I_{x,u}\mu_{u,y}^{s,I}-v_sp^I_{x,y}
\in\ca_{<0}$.
\nl
By 6.3(ii) for $W_I$ we have

$\mu^{s,I}_{x,y}+\sum_{u\in W_I;x<u<y;su<u}p^I_{x,u}\mu_{u,y}^{s,I}
-v_sp^I_{x,y}\in\ca_{<0}$.
\nl
It follows that $\mu^s_{xz,yz}-\mu^{s,I}_{x,y}\in\ca_{<0}$. On the other hand,
$\mu^s_{xz,yz}-\mu^{s,I}_{x,y}$ is fixed by $\bar{}$ (see 6.3(ii)) hence it is
$0$. This proves (e). The lemma is proved.

\proclaim{Proposition 9.11} Assume that $L(s)>0$ for all $s\in I$.

(a) Let $z\in W$ be such that $z$ is the element of
minimal length of $W_Iz$. If $x,y$ in $W_I$ satisfy $x\le_\cl y$ (relative to 
$W_I$), then $xz\le_\cl yz$ (in $W$). If $x,y$ in $W_I$ satisfy $x\sim_\cl y$ 
(relative to $W_I$), then $xz\sim_\cl yz$ (in $W$). 

(b)  Let $z\in W$ be such that $z$ is the element of minimal length of $zW_I$.
If $x,y$ in $W_I$ satisfy $x\le_{\Cal R}y$ (relative to $W_I$), then 
$zx\le_{\Cal R}zy$ (in $W$). If $x,y$ in $W_I$ satisfy $x\sim_{\Cal R}y$ 
(relative to $W_I$), then $zx\sim_{\Cal R}zy$ (in $W$). 
\endproclaim
We prove the first assertion of (a). We may assume that $x\gets_\cl y$ 
(relative to $W_I$) and $x\ne y$. Thus, there exists $s\in I$ such that $sy>y$,
$sx<x$ and we have either $x=sy$ or $x<y$ and $\mu^s_{x,y}\ne 0$. If $x=sy$, 
then $sxz<xz=syz>yz$, hence $xz\gets_\cl yz$ (in $W$). Thus, we may assume that
$x<y$ and $\mu^s_{x,y}\ne 0$. By 9.10(e) we then have $\mu^s_{xz,yz}\ne 0$, 
hence $xz\gets_\cl yz$ (in $W$). The first assertion of (a) is proved. The 
second assertion of (a) follows from the first. (b) follows by applying (a) to
$z\i, x\i,y\i$ instead of $z,x,y$.

\subhead 9.12\endsubhead 
Assume that $z\in W$ is such that $W_Iz=zW_I$ and $z$ is the element of minimal
length of $W_Iz=zW_I$. Then $y\mto z\i yz$ is an automorphism of $W_I$. If 
$s\in I$ then, by 9.7, we have $l(sz)=l(s)+l(z)=1+l(z)$; by 9.7 applied to 
$W_Iz\i$ instead of $W_Iz$ we have $l((z\i sz)z\i)=l(z\i sz)+l(z\i)$ hence 
$l(z\i s)=l(z\i sz)+l(z\i)$; since $l(z\i s)=l(sz)$ and $l(z\i)=l(z)$, it 
follows that $l(z\i sz)+l(z\i)=1+l(z)$, hence $l(z\i sz)=1$. We see that 
$y\mto z\i yz$ maps $I$ onto itself hence it is an automorphism of $W_I$ as a 
Coxeter group. This automorphism preserves the function $L|_{W_I}$. Indeed, if
$y\in W_I$, then

$l(zyz\i)+l(z)=l((zyz\i)z)=l(zy)=l(y\i z\i)=l(y\i)+l(z\i)=l(y)+l(z)$
\nl
(by 9.7 applied to $W_Iz$ and to $W_Iz\i$) hence 
$$\align&L(zyz\i)+L(z)=L((zyz\i)z)=L(zy)=L(y\i z\i)\\&
=L(y\i)+L(z\i)=L(y)+L(z),\endalign$$
so that $L(zyz\i)=L(z)$. In particular, this automorphism respects the 
preorders $\le_\cl,\le_{\Cal R},\le_{\cl\Cal R}$ of $W_I$ (defined in terms of 
$L|_{W_I}$) and the associated equivalence relations.

\proclaim{Proposition 9.13} Assume that $L(s)>0$ for all $s\in I$. Let $z$ be 
as in 9.12. If $x,y$ in $W_I$ satisfy $x\le_{\cl\Cal R}y$ (relative to $W_I$),
then $xz\le_{\cl\Cal R}yz$ (in $W$). If $x,y$ in $W_I$ satisfy 
$x\sim_{\cl\Cal R}y$ (relative to $W_I$), then $xz\sim_{\cl\Cal R}yz$ (in $W$).
\endproclaim
We prove the first assertion. We may assume that either $x\le_\cl y$ (in $W_I)$
or $x\le_{\Cal R}y$ (in $W_I$). In the first case then, by 9.11(a), we have 
$xz\le_\cl yz$ (in $W$) hence $xz\le_{\cl\Cal R}yz$ (in $W$). In the second 
case, by 9.12, we have $z\i xz\le_{\Cal R}z\i yz$. Applying 9.11(b) to 
$z\i xz,z\i yz$ instead of $x,y$ we see that $xz\le_{\Cal R}yz$ (in $W$) hence
$xz\le_{\cl\Cal R}yz$ (in $W$). This proves the first assertion. The second
assertion follows from the first.

\head 10. Inversion \endhead
\subhead 10.1\endsubhead
For $y,w\in W$ we define $q'_{y,w}\in\ca$ by

$q'_{y,w}=\sum (-1)^np_{z_0,z_1}p_{z_1,z_2}\dots p_{z_{n-1},z_n}$
\nl
(sum over all sequences $y=z_0<z_1<z_2<\dots<z_n=w$ in $W$) and

$q_{y,w}=\sgn(y)\sgn(w)q'_{y,w}$. 
\nl
We have

$q_{w,w}=1$,

$q_{y,w}\in\ca_{<0}$ if $y\ne w$,

$q_{y,w}=0$ unless $y\le w$.

\proclaim{Proposition 10.2} For any $y,w \in W$ we have

$\bar q_{y,w}=\sum_{z; y\le z\le w}q_{y,z}r_{z,w}$.
\endproclaim
The (triangular) matrices $Q'=(q'_{y,w}),P=(p_{y,w}),R=(r_{y,w})$ are related 
by

(a) $Q'P=PQ'=1$, $\bar P=RP$, $\bar RR=R\bar R=1$
\nl
where $\bar{}$ over a matrix is the matrix obtained by applying $\bar{}$ to 
each entry. (Although the matrices may be infinite, the products are well 
defined as each entry of a product is obtained by finitely many operations.)
The last three equations in (a) are obtained from 5.3, 4.6; the equations
involving $Q'$ follow from the definition. From (a) we deduce 

$Q'P=1=\bar Q'\bar P=\bar Q'RP$.
\nl
Hence $Q'P=\bar Q'RP$. Multiplying on the right by $Q'$ and using $PQ'=1$ we 
deduce $Q'=\bar Q'R$. Multiplying on the right by $\bar R$ we deduce 

(b) $\bar Q'=Q'\bar R$.
\nl
Let $\bos$ be the matrix whose $y,w$ entry is $\sgn(y)\de_{y,w}$. We have
$\bos^2=1$. Let $Q$ be the triangular matrix $(q_{y,w})$. Note that 
$Q=\bos Q'\bos$. By 4.5 we have $\bar R=\bos R\bos$. Hence by multiplying the 
two sides of (b) on the left and right by $\bos$ we obtain $\bar Q=QR$. The 
proposition is proved.

\subhead 10.3\endsubhead
Let $\tau:\ch@>>>\ca$ be the $\ca$-linear map defined by $\tau(T_w)=\de_{w,1}$
for $w\in W$. 

\proclaim{Lemma 10.4} (a) For $x,y\in W$ we have $\tau(T_xT_y)=\de_{xy,1}$.

(b) For $h,h'\in\ch$ we have $\tau(hh')=\tau(h'h)$.

(c) Let $x,y,z\in W$ and let $M=\min(L(x),L(y),L(z))$. We 
have $\tau(T_xT_yT_z)\in v^M\bz[v\i]$.
\endproclaim
We prove (a) by induction on $l(y)$. If $l(y)=0$, the result is clear. Assume 
now that $l(y)\ge 1$. If $l(xy)=l(x)+l(y)$ then $T_xT_y=T_{xy}$ and the result
is clear. Hence we may assume that $l(xy)\ne l(x)+l(y)$. Then 
$l(xy)<l(x)+l(y)$. Let $y=s_1s_2\dots s_q$ be a reduced expression. We can find
$i\in[1,q]$ such that 

(d) $l(x)+i-1=l(xs_1s_2\dots s_{i-1})>l(xs_1s_2\dots s_{i-1}s_i)$.
\nl
We show that

(e) $xs_1s_2\dots s_{i-1}s_{i+1}\dots s_q\ne 1$.
\nl
If (e) does not hold, then $x=s_q\dots s_{i+1}s_{i-1}\dots s_1$, so that

$l(xs_1s_2\dots s_{i-1})=l(s_q\dots s_{i+1}s_{i-1}\dots s_1s_1...s_{i-1})=
l(s_q\dots s_{i+1})=q-i$,

$l(xs_1s_2\dots s_{i-1}s_i)=l(s_q\dots s_{i+1}s_{i-1}\dots s_1s_1...s_i)=
l(s_q\dots s_{i+1}s_i)=q-i+1$,
\nl
contradicting (d). Thus (e) holds. We have 
$$\align&\tau(T_xT_y)=\tau(T_xT_{s_1}T_{s_2}\dots T_{s_q})=
\tau(T_{xs_1s_2\dots s_{i-1}}T_{s_i}T_{s_{i+1}\dots s_q})\\&=  
\tau(T_{xs_1s_2\dots s_{i-1}s_i}T_{s_{i+1}\dots s_q})+  
(v_s-v_s\i)\tau(T_{xs_1s_2\dots s_{i-1}}T_{s_{i+1}\dots s_q}).\endalign$$
Using the induction hypothesis and (d) we see that this equals

$\de_{xs_1s_2\dots s_{i-1}s_is_{i+1}\dots s_q,1}+  
(v_s-v_s\i)\de_{xs_1s_2\dots s_{i-1}s_{i+1}\dots s_q,1}=\de_{xy,1}$.  
\nl
This completes the proof of (a). To prove (b), we may assume that 
$h=T_x,h'=T_y$ for $x,y\in W$; we then use (a) and the obvious equality 
$\de_{xy,1}=\de_{yx,1}$. 

We prove (c). Using b) we see that 
$\tau(T_xT_yT_z)=\tau(T_yT_zT_x)=\tau(T_zT_xT_y)$. Hence it is enough to show 
that, for any $x,y,z$ we have

$\tau(T_xT_yT_z)\in v^{L(x)}\bz[v\i]$.
\nl
We argue by induction on $l(x)$. If $x=1$, the result follows from (a). Assume 
now that $l(x)\ge 1$. We can find $s\in S$ such that $xs<x$. If $sy>y$, then by
the induction hypothesis,

$\tau(T_xT_yT_z)=\tau(T_{xs}T_{sy}T_z)\in v^{L(x)-L(s)}\bz[v\i]\sub
v^{L(x)}\bz[v\i]$.
\nl
If $sy<y$, then by the induction hypothesis,
$$\align&\tau(T_xT_yT_z)=\tau(T_{xs}T_{sy}T_z)+(v_s-v_s\i)\tau(T_{xs}T_yT_z)\\&
\in v^{L(x)-L(s)}\bz[v\i]+v_sv^{L(x)-L(s)}\bz[v\i]\sub v^{L(x)}\bz[v\i].
\endalign$$
The lemma is proved.

\subhead 10.5\endsubhead
Let $\ch'=\Hom_\ca(\ch,\ca)$. We regard $\ch'$ as a left $\ch$-module where,
for $h\in\ch,\phi\in\ch'$ we have $(h\phi)(h_1)=\phi(h_1h)$ for all $h_1\in\ch$
and as a right $\ch$-module where, for $h\in\ch,\phi\in\ch'$, we have 
$(\phi h)(h_1)=\phi(hh_1)$ for all $h_1\in\ch$. 

\subhead 10.6\endsubhead
We sometimes identify $\ch'$ with the set of all formal sums
$\sum_{x\in W}a_xT_x$ with $a_x\in\ca$; to $\phi\in\ch'$ corresponds the formal
sum $\sum_{x\in W}\phi(T_{x\i})T_x$. Since $\ch$ is contained in the set of 
such formal sums (it is the set of sums such that $c_x=0$ for all but finitely
many $x$), we see that $\ch$ is naturally a subset of $\ch'$. Using 10.4(a) we
see that the imbedding $\ch\sub\ch'$ is an imbedding of $\ch$-bimodules; it is
an equality if $W$ is finite.

\subhead 10.7\endsubhead
Let $z\in W$. Recall that in 8.1 we have defined $D_z\in\ch'$ by
$D_z(c_w)=\de_{z,w}$ for all $w$. An equivalent definition is

(a) $D_z(T_y)=q'_{z,y}$ 
\nl
for all $y\in W$. Indeed, assuming that (a) holds, we have

$D_z(c_w)=\sum_yq'_{z,y}p_{y,w}=\de_{z,w}$.

\proclaim{Proposition 10.8} Let $z\in W,s\in S$. Assume that $L(s)>0$.

(a) If $zs<z$, then $c_sD_z=(v_s+v_s\i)D_z+
D_{zs}+\sum_{u; z<u<us} \mu^s_{z\i,u\i}D_u$.

(b) If $zs>z$, then $c_sD_z=0$.
\endproclaim
For $a,b\in W$ we define $\de_{a<b}$ to be $1$ if $a<b$ and $0$ otherwise. Let
$w\in W$. If $ws>w$, then by 6.7(a), we have
$$\align&(c_sD_z)(c_w)=D_z(c_wc_s)=D_z(c_{ws}+\sum\Sb x\\xs<x<w\endSb
\mu_{x\i,w\i}^sc_x)\\&=\de_{z,ws}+\sum\Sb x\\xs<x<w\endSb
\mu_{x\i,w\i}^s\de_{z,x}.\tag c\endalign$$
If $ws<w$, then by 6.7(b), we have
$$(c_sD_z)(c_w)=D_z(c_wc_s)=(v_s+v_s\i)D_z(c_w)=(v_s+v_s\i)\de_{z,w}.\tag d$$
If $zs<z,ws>w$, then by (c):
$$(c_sD_z)(c_w)=\de_{zs,w}+\de_{z<w}\mu_{z\i,w\i}^s
=((v_s+v_s\i)D_z+D_{zs}+\sum\Sb u\\z<u<us\endSb \mu^s_{z\i,u\i}D_u)(c_w).$$
If $zs<z, ws<w$, then by (d):
$$(c_sD_z)(c_w)=(v_s+v_s\i)\de_{z,w}=((v_s+v_s\i)D_z+D_{zs}+
\sum\Sb u\\ z<u<us\endSb\mu^s_{z\i,u\i}D_u)(c_w).$$
If $zs>z, ws>w$, then by (c), we have $(c_sD_z)(c_w)=0$. If $zs>z,ws<w$, then 
by (d), we have $(c_sD_z)(c_w)=0$. Since $(c_w)$ is an $\ca$-basis of $\ch$, 
the proposition follows.

\head 11. The longest element for a finite $W$ \endhead
\subhead 11.1\endsubhead
Let $I\sub S$ be such that $W_I$ is finite. By 9.8, there is a unique element 
of maximal length of $W_I$. We denote it by $w^I_0$. If $w_1$ has minimal
length in $W_Ia$ then $w^I_0w_1$ has maximal length in $W_Ia$.

\subhead 11.2\endsubhead
In the remainder of this section we assume that $W$ is finite. Then $w^S_0$, 
the unique element of maximal length of $W$ is well defined. Traditionally one
writes $w_0$ instead of $w^S_0$. Since $l(w_0\i)=l(w_0)$, we must have 

$w_0\i=w_0$.
\nl
By the argument in the proof of 9.8 we have $w\le w_0$ for any $w\in W$. By 9.8
we have 

(a) $l(ww_0)=l(w_0)-l(w)$
\nl
for any $w\in W$. Applying this to $w\i$ and using the equalities
$l(w\i w_0)=l(w_0\i w)=l(w_0w), l(w\i)=l(w)$, we deduce that

(b) $l(w_0w)=l(w_0)-l(w)$.
\nl
We can rewrite (a),(b) as 

$l(w_0)=l(w\i)+l(ww_0)$, $l(w_0)=l(w_0w)+l(w\i)$. 
\nl
Using this and the definition of $L$ we deduce that

$L(w\i)+L(ww_0)=L(w_0)=L(w_0w)+L(w\i)$, 
\nl
hence $L(ww_0)=L(w_0w)$. This implies $L(w_0ww_0)=L(w)$ for all $w$. Replacing
$L$ by $l$ we deduce that $l(w_0ww_0)=l(w)$. Thus, the (involutive) 
automorphism $w\mto w_0ww_0$ of $W$ maps $S$ into itself hence is a Coxeter 
group automorphism preserving the function $L$. 

\proclaim{Lemma 11.3} Let $y,w\in W$. We have

(a) $y\le w \Leftrightarrow w_0w\le w_0y \Leftrightarrow ww_0\le yw_0$.

(b) $r_{y,w}=r_{ww_0,yw_0}=r_{w_0w,w_0y}$;

(c) $\bar p_{ww_0,yw_0}=\sum_{z; y\le z\le w}p_{zw_0,yw_0}r_{z,w}$.
\endproclaim
We prove (a). To prove that $y\le w \implies w_0w\le w_0y$, we may assume that
$l(w)-l(y)=1, yw\i\in T$. Then 

$l(yw_0)-l(ww_0)=l(w_0)-l(y)-(l(w_0)-l(w))=l(w)-l(y)=1$ 
\nl
and $(ww_0)(yw_0)\i=wy\i\in T$. Hence $w_0w\le w_0y$. The opposite implication
is proved in the same way. The second equivalence in (a) follows from the last
sentence in 11.2.

We prove the first equality in (b) by induction on $l(w)$. If $l(w)=0$ then 
$w=1$. We have $r_{y,1}=\de_{y,1}$. Now $r_{w_0,w_0y}$ is zero unless 
$w_0\le w_0y$ (see 4.7). On the other hand we have $w_0y\le w_0$ (see 11.2).
Hence $r_{w_0,w_0y}$ is zero unless $w_0y=w_0$, that is unless $y=1$ in which 
case it is $1$. Thus the desired equality holds when $l(w)=0$. Assume now that 
$l(w)\ge 1$. We can find $s\in S$ such that $sw<w$. Then $sww_0>ww_0$ by (a).

Assume first that $sy<y$ (hence $syw_0>yw_0$.) Using 4.4 and the induction 
hypothesis we have

$r_{y,w}=r_{sy,sw}=r_{sww_0,syw_0}=r_{ww_0,yw_0}$.
\nl
Assume next that $sy>y$ (hence $syw_0<yw_0$.) Using 4.4 and the induction 
hypothesis we have
$$\align&r_{y,w}=r_{sy,sw}+(v_s-v_s\i)r_{y,sw}
=r_{sww_0,syw_0}+(v_s-v_s\i)r_{sww_0,yw_0}\\&
=r_{sww_0,syw_0}+(v_s-v_s\i)r_{ww_0,syw_0}=r_{ww_0,yw_0}.\endalign$$
This proves the first equality in (b). The second equality in (b) follows from
the last sentence in 11.2. 

We prove (c). We may assume that $y\le w$. By 5.3 (for $ww_0,yw_0$ instead of 
$y,w$ we have

$\bar p_{ww_0,yw_0}=\sum_{z;y\le z\le w}r_{ww_0,zw_0}p_{zw_0,yw_0}$
\nl
(we have used (a)). Here we substitute $r_{ww_0,zw_0}=r_{z,w}$ (see (b)) and 
the result follows. The lemma is proved.

\proclaim{Proposition 11.4} For any $y,w\in W$ we have $q_{y,w}=p_{ww_0,yw_0}=
p_{w_0w,w_0y}$.
\endproclaim
The second equality follows from the last sentence in 11.2. We prove the first
equality. We may assume that $y\le w$. We argue by induction on 
$l(w)-l(y)\ge 0$. If $l(w)-l(y)=0$ we have $y=w$ and both sides are $1$. Assume
now that $l(w)-l(y)\ge 1$. Substracting the identity in 11.3(c) from that in 
10.2 and using the induction hypothesis, we obtain

$\bar q_{y,w}-\bar p_{ww_0,yw_0}=q_{y,w}-p_{ww_0,yw_0}$.
\nl
The right hand side is in $\ca_{<0}$; since it is fixed by $\bar{}$, it is $0$.
The proposition is proved.

\proclaim{Proposition 11.5} We identify $\ch=\ch'$ as in 10.6. If $z\in W$, 
then $D_{z\i}\in\ch'$ (see 10.7) becomes an element of $\ch$. We have

$D_{z\i}T_{w_0}\i=\sgn(zw_0)c_{zw_0}^\dag,\quad (\text{${}^\dag$ as in 3.5})$. 
\endproclaim
By definition, $D_{z\i}\in\ch$ is characterized by

$\tau(D_{z\i}T_{y\i})=q'_{z\i,y\i}$ 
\nl
for all $y\in W$. Here $\tau$ is as in 10.3. Hence, by 10.4(a), we have 
$D_{z\i}=\sum_yq'_{z\i,y\i}T_y$. Using 11.4, we deduce

$D_{z\i}=\sum_y\sgn(yz)p_{w_0y\i,w_0z\i}T_y$.
\nl
Multiplying on the right by $T_{w_0}\i$ gives

$D_{z\i}T_{w_0}\i=\sum_y\sgn(yz)p_{w_0y\i,w_0z\i}T_{w_0y\i}\i$
\nl
since $T_{w_0y\i}T_y=T_{w_0}$. On the other hand,
$$\sgn(zw_0)c_{zw_0}^\dag=\sum_x\sgn(zw_0x)p_{x,zw_0}T_{x\i}\i
=\sum_y\sgn(zw_0yw_0)p_{yw_0,zw_0}T_{w_0y\i}\i.$$
We now use the identity $p_{yw_0,zw_0}=p_{w_0y\i,w_0z\i}$. The proposition 
follows.

\proclaim{Proposition 11.6} Let $u,z\in W, s\in S$ be such that $sz<z<u<su$. 
Assume that $L(s)>0$. Then $suw_0<uw_0<zw_0<szw_0$ and

$\mu_{uw_0,zw_0}^s=-\sgn(uz)\mu^s_{z,u}$.
\endproclaim
Let $z\in W,s\in S$ be such that $sz<z$. Using 10.8(a), we see that

$(c_s-(v_s+v_s\i))D_{z\i}T_{w_0}\i=
D_{z\i s}T_{w_0}\i+\sum_{u;z\i<u\i<u\i s}\mu^s_{z,u}D_{u\i}T_{w_0}\i$,
\nl
hence, using 11.5, we have
$$(c_s-(v_s+v_s\i))\sgn(zw_0)c_{zw_0}^\dag=\sgn(zsw_0)c_{szw_0}^\dag+
\sum_{u;z<u<su}\mu^s_{z,u}\sgn(uw_0)c_{uw_0}^\dag.$$
Applying ${}^\dag$ to both sides and using $(c_s-(v_s+v_s\i))^\dag=-c_s$ gives

(a) $-c_sc_{zw_0}=-c_{szw_0}+\sum_{u;z<u<su}\mu^s_{z,u}\sgn(uz)c_{uw_0}$.
\nl
Since $szw_0>zw_0$, we can apply 6.6(a) and we get

$c_sc_{zw_0}=c_{szw_0}+\sum_{u';su'<u'<zw_0}\mu_{u',zw_0}^sc_{u'}$
\nl
or equivalently

$c_sc_{zw_0}=c_{szw_0}+\sum_{u;z<u<su}\mu_{uw_0,zw_0}^sc_{uw_0}$.
\nl
Comparing this with (a) we get

$-\sum_{u;z<u<su}\mu^s_{z,u}\sgn(uz)c_{uw_0}=
\sum_{u;z<u<su}\mu_{uw_0,zw_0}^sc_{uw_0}$;
\nl
the proposition follows.

\proclaim{Corollary 11.7}Assume that $L(s)>0$ for all $s\in S$. Let $y,w\in W$.

(a) $y\le_\cl w\Lra ww_0\le_\cl yw_0\Lra w_0w\le_\cl w_0y$;

(b) $y\le_{\Cal R}w\Lra ww_0\le_{\Cal R}yw_0\Lra w_0w\le_{\Cal R}w_0y$;

(c) 
$y\le_{\cl\Cal R}w\Lra ww_0\le_{\cl\Cal R}yw_0\Lra w_0w\le_{\cl\Cal R}w_0y$.

(d) Left multiplication by $w_0$ carries left cells to left cells, right cells
to right cells, two-sided cells to two-sided cells. The same holds for right
multiplication by $w_0$.
\endproclaim
We prove the first equivalence in (a). It is enough to show that
$y\le_\cl w\implies ww_0\le_\cl yw_0$. We may assume that $y\gets_\cl w$ and
$y\ne w$. Then there exists $s\in S$ such that $sw>w,sy<y$ and
$D_y(c_sc_w)\ne 0$. We have $syw_0>yw_0, sww_0<ww_0$. From 6.6 we see that
either $y=ws$ or $y<w$ and $\mu^s_{y,w}\ne 0$. In the first case we have 
$ww_0=syw_0$; in the second case we have $ww_0<yw_0$ and 
$\mu^s_{ww_0,yw_0}\ne 0$ (see 11.6). In both cases, 6.6 shows that 
$D_{ww_0}(c_sc_{yw_0})\ne 0$. Hence $ww_0\le_\cl yw_0$. Thus, the first 
equivalence in (a) is established. The second equivalence in (a) follows from 
the last sentence in 11.2.

Now (b) follows by applying (a) to $y\i,w\i$ instead of (a); (c) follows from
(a) and (b); (d) follows from (a),(b),(c). The corollary is proved.

\head 12. Examples of elements $D_w$\endhead
\proclaim{Proposition 12.1} Assume that $L(s)>0$ for all $s\in S$. For any 
$y\in W$ we have $D_1(T_y)=\sgn(y)v^{-L(y)}$. Equivalently, with the
identification in 10.6, we have

$D_1=\sum_{y\in W}\sgn(y)v^{-L(y)}T_y$.
\endproclaim
An equivalent statement is that $q'_{1,y}=\sgn(y)v^{-L(y)}$. Since $q'_{1,y}$ 
are determined  by the equations $\sum_yq'_{1,y}p_{y,w}=\de_{1,w}$ (see 
10.2(a)) it is enough to show that

$\sum_y\sgn(y)v^{-L(y)}p_{y,w}=\de_{1,w}$
\nl
for all $w\in W$. If $w=1$ this is clear. Assume now that $w\ne 1$. We can find
$s\in S$ such that $sw<w$. We must prove that

$\sum_{y;y<sy}\sgn(y)v^{-L(y)}(p_{y,w}-v_s\i p_{sy,w})=0$.
\nl
Each term of the last sum is $0$, by 6.6(c). The proposition is proved.

\proclaim{Corollary 12.2} Assume that $W$ is finite and that $L(s)>0$ for all 
$s\in S$. Then 

$c_{w_0}=\sum_{y\in W}v^{-L(yw_0)}T_y$.
\endproclaim
This follows immediately from 12.1 and 11.5. Alternatively, we can argue as
follows. We prove that $p_{y,w_0}=v^{-L(yw_0)}$ for all $y$, by descending 
induction on $l(y)$. If $l(y)$ is maximal, that is $y=w_0$, then 
$p_{y,w_0}=1$. Assume now that $l(y)<l(w_0)$. We can find $s\in S$ such that 
$l(sy)=l(y)+1$. By the induction hypothesis we have $p_{sy,w_0}=v^{-L(syw_0)}$.
By 6.6(c), we have
$$\align&p_{y,w_0}\\&=v_s\i p_{sy,w_0}=v^{-L(s)-L(syw_0)}
=v^{-L(s)-L(w_0)+L(sy)}=v^{-L(w_0)+L(y)}=v^{-L(yw_0)}.\endalign$$
The corollary is proved.

\subhead 12.3\endsubhead
From 11.5 we see that $D_{z\i}$ can be explicitly computed whenever $W$ is
finite and $c_{zw_0}$ is known. In particular, in the setup of 7.4 with 
$m=2k+2<\infty$, we can compute explicitly all $D_{z\i}$ using 7.6(a). For 
example:
$$\align&D_{s_1}=\sum_{s\in[0,k-1]}(1-v^{2L_1}+v^{4L_1}-\dots+(-1)^sv^{2sL_1})
v^{-sL_1-sL_2}\\&\tim (T_{1_{2s+1}}-v^{-L_2}T_{2_{2s+2}}-v^{-L_2}T_{1_{2s+2}}
+v^{-2L_2}T_{2_{2s+3}})\\&+(1-v^{2L_1}+v^{4L_1}-\dots+(-1)^kv^{2kL_1})
v^{-kL_1-kL_2}(T_{1_{2k+1}}-v^{-L_2}T_{2_{2k+2}}).\endalign$$
Using this (for larger and larger $m$) one can deduce that an analogous formula
holds in the setup of 7.4 with $m=\infty$:
$$\align&D_{s_1}=\sum_{s\ge 0}(1-v^{2L_1}+v^{4L_1}-\dots+(-1)^sv^{2sL_1})
v^{-sL_1-sL_2}\\&\tim (T_{1_{2s+1}}-v^{-L_2}T_{2_{2s+2}}-v^{-L_2}T_{1_{2s+2}}
+v^{-2L_2}T_{2_{2s+3}})\in\ch'.\tag a\endalign$$
(We use the identification in 10.6.)

\head 13. The function $\aa$\endhead
\subhead 13.1\endsubhead
{\it In the remainder of these lectures we assume that $L(s)>0$ for all
$s\in S$.} A reference for this section is \cite{\LC}.

For $x,y,z$ in $W$ we define $f_{x,y,z}\in\ca,f'_{x,y,z}\in\ca,h_{x,y,z}\in\ca$
by

$T_xT_y=\sum_{z\in W}f_{x,y,z}T_z=\sum_{z\in W}f'_{x,y,z}c_z$,

$c_xc_y=\sum_{z\in W}h_{x,y,z}c_z$.
\nl
We have

(a) $f_{x,y,z}=\sum_{z'}p_{z,z'}f'_{x,y,z'}$

(b) $f'_{x,y,z}=\sum_uq'_{z,z'}f_{x,y,z'}$,

(c) $h_{x,y,z}=\sum_{x',y'}p_{x',x}p_{y',y}f'_{x',y',z}$,
\nl
(a),(c) follow from the definitions; (b) follows from (a) using 10.2(a). 
All sums in (a)-(c) are finite. From 8.2, 5.6, we see that

(d) $h_{x,y,z}\ne 0\implies z\le_{\Cal R}x, z\le_\cl y$,

(e) $h_{x,y,z}=h_{y\i,x\i,z\i}$. 

\subhead 13.2\endsubhead
We say that $N\in\bn$ is a {\it bound} for $W,L$ if 
$v^{-N}f_{x,y,z}\in\ca_{\le 0}$ for all $x,y,z$ in $W$. We say that $W,L$ is 
{\it bounded} if there exists $N\in\bn$ such that $N$ is a bound for $W,L$. If
$W$ is finite then $W,L$ is obviously bounded. More precisely:

\proclaim{Lemma 13.3} If $W$ is finite, then $N=L(w_0)$ is a bound for $W,L$.
\endproclaim
By 10.4(a) we have $f_{x,y,z}=\tau(T_xT_yT_{z\i})$. By 10.4(c) we have
$\tau(T_xT_yT_{z\i})\in v^{L(w_0)}\bz[v\i]$. The lemma is proved.

\subhead 13.4\endsubhead
More generally, according to \cite{\LC, 7.2}, if $W$ is tame and $L=l$ then 
$W,L$ is bounded; as a bound we can take 

(a) $N=\max_IL(w_0^I)$
\nl
where $I$ runs over the subsets of $S$ such that $W_I$ is finite. This result 
and its proof remain valid without the assumption $L=l$.

We illustrate this in the setup of 7.1 with $m=\infty$. For $a,b\in\{1,2\}$ and
$k>0,k'>0$, we have

$T_{a_k}T_{b_{k'}}=T_{a_{k+k'}}$ if $b=a+k\mod 2$,

$T_{a_k}T_{b_{k'}}
=T_{a_kb_{k'}}+\sum_{u\in[1,\min(k,k')]}\xi_{b+u-1}T_{a_{k+k'-2u+1}}$
if $b=a+k+1\mod 2$;
\nl
here, for $n\in\bz$ we set $\xi_n=v^{L_1}-v^{-L_1}$ if $n$ is odd and 
$\xi_n=v^{L_2}-v^{-L_2}$ if $n$ is even. We see that, in this case, 
$\max(L_1,L_2)$ is a bound for $W,L$.

{\it Question.} Is it true, in the general case, that $W,L$ admits a bound? If
so, can the bound be taken as in (a)?

\proclaim{Lemma 13.5} Assume that $N$ is a bound for $W,L$. Then, for any
$x,y,z$ in $W$ we have

(a) $v^{-N}f'_{x,y,z}\in\ca_{\le 0}$,

(b) $v^{-N}h_{x,y,z}\in\ca_{\le 0}$.
\endproclaim
(a) follows from 13.1(b) since $q'_{z,z'}\in\ca_{\le 0}$. (b) follows from (a)
and 13.1(c) since $p_{x',x}\in\ca_{\le 0}$, $p_{y',y}\in\ca_{\le 0}$.

\subhead 13.6\endsubhead
{\it In the remainder of this section we assume that $W,L$ is bounded.} Let $N$
be a bound for $W,L$. By 13.5(b), for any $z\in W$ there exists a unique 
integer $\aa(z)\in[0,N]$ such that

(a) $h_{x,y,z}\in v^{\aa(z)}\bz[v\i]$ for all $x,y\in W$,

(b) $h_{x,y,z}\notin v^{\aa(z)-1}\bz[v\i]$ for some $x,y\in W$.
\nl
(We use that $h_{1,z,z}=1$.) We then have for any $x,y,z$:

(c) $h_{x,y,z}=\ga_{x,y,z\i}v^{\aa(z)}\mod v^{\aa(z)-1}\bz[v\i]$
\nl
where $\ga_{x,y,z\i}\in\bz$ is well defined; moreover, for any $z\in W$ there
exists $x,y$ such that $\ga_{x,y,z\i}\ne 0$. 

For any $x,y,z$ we have 

(d) $f'_{x,y,z}=\ga_{x,y,z\i}v^{\aa(z)}\mod v^{\aa(z)-1}\bz[v\i]$.
\nl
This is proved (for fixed $z$) by induction on $l(x)+l(y)$ using (c) and 
13.1(c). (Note that $p_{x',x}p_{y',y}$ is $1$ if $x'=x,y'=y$ and is in 
$\ca_{<0}$ otherwise.)

\proclaim{Proposition 13.7} (a) $\aa(1)=0$.

(b) If $z\in W-\{1\}$, then $\aa(z)\ge\min_{s\in S}L(s)>0$.
\endproclaim
We prove (a). Let $x,y\in W$. Assume first that $y\ne 1$. We can find $s\in S$
such that $ys<y$. Then $c_y\in{}^s\ch$. Since ${}^s\ch$ is a left ideal (see 
8.4) we have $c_xc_y\in{}^s\ch$. Since $s1>1$, from the definition of 
${}^s\ch$ it then follows that $h_{x,y,1}=0$. 

Similarly, if $x\ne 1$, then $h_{x,y,1}=0$. Since $h_{1,1,1}=1$, (a) follows.

In the setup of (b) we can find $s\in S$ such that $sz<z$. By 6.6(b) we have
$h_{s,z,z}=v_s+v_s\i$. This shows that $\aa(z)\ge L(s)>0$. The proposition is 
proved.

\proclaim{Proposition 13.8}Assume that $W$ is finite. 

(a) We have $\aa(w_0)=L(w_0)$.

(b) For any $w\in W-\{w_0\}$ we have $\aa(w)<L(w_0)$. 
\endproclaim
For any $w\in W$ we have by definition $\aa(w)\le N$ where $N$ is a bound for
$W,L$; by 13.3 we can take $N=L(w_0)$, hence $\aa(w)\le L(w_0)$.

We prove (a). From 6.6(b) we see that $T_sc_{w_0}=v_sc_{w_0}$ for any $s\in S$.
Using this and 12.2, we see that 

$c_{w_0}c_{w_0}=\sum_{y\in W}v^{-L(yw_0)}v^{L(y)}c_{w_0}$,
\nl
hence 

$h_{w_0,w_0,w_0}=\sum_{y\in W}v^{-L(w_0)}v^{2L(y)}\in v^{L(w_0)}\mod
v^{L(w_0)-1}\bz[v\i]$.
\nl
It follows that $\aa(w_0)\ge L(w_0)$. Hence $\aa(w_0)=L(w_0)$. This proves (a).

We prove (b). Let $z\in W$ be such that $\aa(z)=L(w_0)$. We must prove that
$z=w_0$. By 13.6(d), we can find $x,y$ such that 

$f'_{x,y,z}=bv^{L(w_0)}+\text{strictly smaller powers of } v$
\nl
where $b\in\bz-\{0\}$. For any $z'\ne z$ we have

$f'_{x,y,z'}\in v^{L(w_0)}\bz[v\i]$
\nl
(by 13.6 and the first sentence in the proof). Since $p_{z,z'}=1$ for $z=z'$ 
and $p_{z,z'}\in\ca_{<0}$ for $z'<z$, we see that the equality 
$f_{x,y,z}=\sum_{z'}p_{z,z'}f'_{x,y,z'}$ (see 13.1(a)), implies that 

$f_{x,y,z}=bv^{L(w_0)}+\text{strictly smaller powers of } v$
\nl
with $b\ne 0$. Now $f_{x,y,z}=\tau(T_xT_yT_{z\i})$. Using now 10.4(c) we see
that 

$\min(L(x),L(y),L(z\i))=L(w_0)$.
\nl
It follows that $x=y=z\i=w_0$. The proposition is proved.

\proclaim{Proposition 13.9} (a) For any $z\in W$ we have $\aa(z)=\aa(z\i)$.

(b) For any $x,y,z\in W$ we have $\ga_{x,y,z}=\ga_{y\i,x\i,z\i}$.
\endproclaim
This follows from 13.1(e).

\subhead 13.10\endsubhead
We show that, in the setup of 7.1 with $m=\infty$ and $L_2\ge L_1$, the 
function $\aa:W@>>>\bn$ is given as follows:

(a) $\aa(1)=0$,

(b) $\aa(1_1)=L_1, \aa(2_1)=L_2$,
 
(c) $\aa(1_k)=\aa(2_k)=L_2$ if $k\ge 2$.
\nl
Now (a) is contained in 13.7(a). If $s_2z<z$ then, by the proof of 13.7(b) we 
have $\aa(z)\ge L_2$. By 13.4, $L_2$ is a bound for $W,L$ hence $\aa(z)\le L_2$
so that $\aa(z)=L_2$. If $zs_2<z$ then the previous argument is applicable to 
$z\i$. Using 13.9, we see that $\aa(z)=\aa(z\i)=L_2$. 

Assume next that $z=1_{2k+1}$ where $k\ge 1$. By 7.5, 7.6, we have 

$c_{1_2}c_{2_{2k}}=c_1c_2c_{2_{2k}}
=(v^{L_2}+v^{-L_2})c_1c_{2_{2k}}=(v^{L_2}+v^{-L_2})c_{1_{2k+1}}$,
\nl
hence $h_{1_2,2_{2k},z}=v^{L_2}+v^{-L_2}$. Thus, $\aa(z)\ge L_2$. 
By 13.4 we have
$\aa(z)\le L_2$ hence $\aa(z)=L_2$.

It remains to consider the case where $z=s_1$. Assume first that $L_1=L_2$. 
Then $\aa(s_1)\le L_1$ by 13.4 and $\aa(s_1)\ge L_1$ by 13.7(b). Hence 
$\aa(s_1)=L_1$.

Assume next that $L_1<L_2$. Then $\ci=\sum_{w\in W-\{1,s_1\}}\ca c_w$ is a 
two-sided ideal $\ci$ of $\ch$ (see 8.8). Hence if $x$ or $y$ is in 
$W-\{1,s_1\}$, then $c_xc_y\in\ci$ and $h_{x,y,s_1}=0$. Using

$h_{1,1,s_1}=0,h_{1,s_1,s_1}=h_{s_1,1,s_1}=1,h_{s_1,s_1,s_1}=v^{L_1}+v^{-L_1}$
\nl
we see that $\aa(s_1)=L_1$. Thus, (a),(b),(c) are established.

\subhead 13.11\endsubhead
In this subsection we assume that we are in the setup of 7.1 with
$4\le m<\infty$ and $L_2>L_1$. By 7.8, we have

$h_{2_{m-1},2_{m-1},2_{m-1}}
=(-1)^{(m-2)/2} v^{(mL_2-(m-2)L_1)/2}+\text{ strictly smaller powers of } v.$
\nl
Hence $\aa(2_{m-1})\ge(mL_2-(m-2)L_1)/2$.

One can show that the function $\aa:W@>>>\bn$ is given as follows:

$\aa(1)=0$,

$\aa(1_1)=L_1, \aa(2_1)=L_2$,
 
$\aa(1_{m-1})=L_2, \aa(2_{m-1})=(mL_2-(m-2)L_1)/2$,

$\aa(2_m)=m(L_1+L_2)/2$,

$\aa(1_k)=\aa(2_k)=L_2$ if $1<k<m-1$.
\nl
This remains true in the case where $L_1=L_2$.

\head  14. Conjectures \endhead
\subhead 14.1\endsubhead
{\it In this section we assume that $W,L$ is bounded.}

For an integer $n$ we denote by $\pi_n:\ca@>>>\bz$ be the group homomorphism 
defined by $\pi_n(\sum_{k\in\bz}a_kv^k)=a_n$.

For $z\in W$ we denote by $\De(z)$ the integer $\ge 0$ defined by

(a) $p_{1,z}=n_zv^{-\De(z)}+\text{strictly smaller powers of } v, 
n_z\in\bz-\{0\}$.
\nl
Note that $\De(1)=0$, $0<\De(z)\le L(z)$ for $z\ne 1$ (see ...) and
$\De(z)=\De(z\i)$ for all $z$. Let 
$$\cd=\{z\in W| \aa(z)=\De(z)\}.$$ 
Clearly, $z\in\cd\implies z\i\in\cd$.

\proclaim{Conjectures 14.2} The following properties hold.

P1. We have $\aa(z)\le\De(z)$.

P2. If $d\in\cd$ and $x,y\in W$ satisfy $\ga_{x,y,d}\ne 0$, then $x=y\i$.

P3. If $y\in W$, there exists a unique $d\in\cd$ such that 
$\ga_{y\i,y,d}\ne 0$.

P4. If $z'\le_{\cl\Cal R}z$ then $\aa(z')\ge\aa(z)$. Hence, if
$z'\sim_{\cl\Cal R}z$, then $\aa(z')=\aa(z)$. 

P5. If $d\in\cd,y\in W,\ga_{y\i,y,d}\ne 0$, then $\ga_{y\i,y,d}=n_d=\pm 1$.

P6. If $d\in\cd$, then $d^2=1$.

P7. For any $x,y,z\in W$ we have $\ga_{x,y,z}=\ga_{y,z,x}$.

P8. Let $x,y,z\in W$ be such that $\ga_{x,y,z}\ne 0$. Then $x\sim_\cl y\i$, 
$y\sim_\cl z\i$, $z\sim_\cl x\i$.

P9. If $z'\le_\cl z$ and $\aa(z')=\aa(z)$ then $z'\sim_\cl z$.

P10. If $z'\le_{\Cal R}z$ and $\aa(z')=\aa(z)$ then $z'\sim_{\Cal R}z$.

P11. If $z'\le_{\cl\Cal R}z$ and $\aa(z')=\aa(z)$ then $z'\sim_{\cl\Cal R}z$.

P12. Let $I\sub S$. If $y\in W_I$, then $\aa(y)$ computed in terms of $W_I$ is
equal to $\aa(y)$ computed in terms of $W$.

P13. Any left cell $\Ga$ of $W$ contains a unique element $d\in\cd$. We have 
$\ga_{x\i,x,d}\ne 0$ for all $x\in\Ga$.

P14. For any $z\in W$ we have $z\sim_{\cl\Cal R}z\i$.

P15. Let $v'$ be a second indeterminate and let $h'_{x,y,z}\in\bz[v',v'{}\i]$ 
be obtained from $h_{x,y,z}$ by the substitution $v\mto v'$. If $x,x',y,w\in W$
satisfy $\aa(w)=\aa(y)$, then 

$\sum_{y'}h'_{w,x',y'}h_{x,y',y}=\sum_{y'}h_{x,w,y'}h'_{y',x',y}.$
\endproclaim
In \S15-\S17 we will verify the conjecture above in a number of cases.

\proclaim{Auxiliary statement 14.3} Let $x,y,z,z'\in W$ be such that 
$\ga_{x,y,z\i}\ne 0$, $z'\gets_\cl z$. Then there exists $x'\in W$ such that 
$\pi_{\aa(z)}(h_{x',y,z'})\ne 0$. In particular, $\aa(z')\ge\aa(z)$.
\endproclaim

In this section we will show, that, if P1-P3 and 14.3 are assumed to be true, 
then P4-P14 are automatically true. The arguments follow \cite{\LC},
\cite{\LCC}.

\subhead 14.4\endsubhead 
{\it 14.3 $\implies$ P4.} Let $z',z$ be as in P4. We can assume that $z'\ne z$
and that $z'\gets_\cl z$ or $z'{}\i\gets_\cl z\i$. In the first case, from 14.3
we get $\aa(z')\ge \aa(z)$. 
(We can find $x,y$ such that $\ga_{x,y,z\i}\ne 0$.) In
the second case, from 14.3 we get $\aa(z'{}\i)\ge \aa(z\i)$ 
hence $\aa(z')\ge \aa(z)$.

\subhead 14.5\endsubhead 
{\it P1,P3 $\implies$ P5.} Let $x,y\in W$. Applying $\tau$ to 
$c_xc_y=\sum_zh_{x,y,z}c_z$ gives
$$\sum_zh_{x,y,z}p_{1,z}=\sum_{x',y'}p_{x',x}p_{y',y}\tau(T_{x'}T_{y'})=
\sum_{x',y'}p_{x',x}p_{y',y}\de_{x'y',1}=\sum_{x'}p_{x',x}p_{x'{}\i,y}$$
hence
$$\sum_zh_{x,y,z}p_{1,z}=\de_{xy,1}\mod v\i\bz[v\i].\tag a$$
We take $x=y\i$ and note that $h_{y\i,y,z}\in v^{\aa(z)}\bz[v\i]$, 
$p_{1,z}\in v^{-\De(z)}\bz[v\i]$, hence

$h_{y\i,y,z}p_{1,z}\in v^{\aa(z)-\De(z)}\bz[v\i]$.
\nl
The same argument shows that, if $z\in\cd$, then 

$h_{y\i,y,z}p_{1,z}\in\ga_{y\i,y,z\i}n_z+v\i\bz[v\i]$.
\nl
If $z\notin\cd$ then, by P1, we have $\aa(z)-\De(z)<0$ so that
$h_{y\i,y,z}p_{1,z}\in v\i\bz[v\i]$. We now see that 

$\sum_zh_{y\i,y,z}p_{1,z}=\sum_{z\in\cd}\ga_{y\i,y,z\i}n_z\mod v\i\bz[v\i]$.
\nl
Comparing with (a) we see that $\sum_{z\in\cd}\ga_{y\i,y,z\i}n_z=1$. 
Equivalently,

$\sum_{z\in\cd}\ga_{y\i,y,z}n_z=1$.
\nl
Using this and P3 we see that, in the setup of P5 we have $\ga_{y\i,y,d}n_d=1$.
Since $\ga_{y\i,y,d},n_d$ are integers, we must have $\ga_{y\i,y,d}=n_d=\pm 1$.

\subhead 14.6\endsubhead 
{\it P2,P3 $\implies$ P6.} We can find $x,y$ such that $\ga_{x,y,d}\ne 0$. By 
P2, we have $x=y\i$ so that $\ga_{y\i,y,d}\ne 0$. This implies 
$\ga_{y\i,y,d\i}\ne 0$. (See 13.9(b)). We have $d\i\in\cd$. By the uniqueness 
in P3 we have $d=d\i$.

\subhead 14.7\endsubhead 
{\it P2,P3,P4,P5 $\implies$ P7.} We first prove the following statement.

(a) {\it Let $x,y,z\in W, d\in\cd$ be such that $\ga_{x,y,z}\ne 0$,
$\ga_{z\i,z,d}\ne 0$, $\aa(d)=\aa(z)$. Then} $\ga_{x,y,z}=\ga_{y,z,x}$.
\nl
Let $n=\aa(d)$. From $\ga_{x,y,z}\ne 0$ we deduce $h_{x,y,z\i}\ne 0$ hence 
$z\i\le_{\Cal R}x$, hence $n=\aa(z)=\aa(z\i)\ge \aa(x)$ (see P4). Computing the
coefficient of $c_d$ in two ways, we obtain 

$\sum_{z'}h_{x,y,z'}h_{z',z,d}=\sum_{x'}h_{x,x',d}h_{y,z,x'}$.
\nl
Now $h_{z',z,d}\ne 0$ implies $d\le_{\Cal R}z'$ hence $\aa(z')\le\aa(d)=n$ (see
P4); similarly, $h_{x,x',d}\ne 0$ implies $d\le_\cl x'$ hence 
$\aa(x')\le\aa(d)=n$. Thus we have 

$\sum_{z';\aa(z')\le n}h_{x,y,z'}h_{z',z,d}=
\sum_{x';\aa(x')\le n}h_{x,x',d}h_{y,z,x'}$.
\nl
By P2 and our assumptions, the left hand side is 

$\ga_{x,y,z}\ga_{z\i,z,d}v^{2n}+\text{strictly smaller powers of } v$. 
\nl
Similarly, the right hand side is 

$\ga_{x,x\i,d}\pi_n(h_{y,z,x\i})v^{2n}+\text{strictly smaller powers of } v$. 
\nl
Hence $\ga_{x,x\i,d}\pi_n(h_{y,z,x\i})=\ga_{x,y,z}\ga_{z\i,z,d}\ne 0$. Thus, 

$\ga_{x,x\i,d}\ne 0, \pi_n(h_{y,z,x\i})\ne 0$.
\nl
We see that $\aa(x\i)\ge n$. But we have also $\aa(x)\le n$ hence $\aa(x)=n$ 
and 
$\pi_n(h_{y,z,x\i})=\ga_{y,z,x}$. Since $\ga_{x,x\i,d}\ne 0$, we have (by P5) 
$\ga_{x,x\i,d}=\ga_{z\i,z,d}$. Using this and 
$\ga_{x,x\i,d}\ga_{y,z,x}=\ga_{x,y,z}\ga_{z\i,z,d}$ we deduce 
$\ga_{y,z,x}=\ga_{x,y,z}$, as required. 

Next we prove the following statement.

(b) {\it Let $z\in W,d\in\cd$ be such that $\ga_{z\i,z,d}\ne 0$. Then}
$\aa(z)=\aa(d)$.
\nl
We shall 
assume that (b) holds whenever $\aa(z)>N_0$ and we shall deduce that it
also holds when $\aa(z)=N_0$. (This will prove (b) by descending induction on 
$\aa(z)$ since $\aa(z)$ is bounded above.) Assume that $\aa(z)=N_0$. From 
$\ga_{z\i,z,d}=\pm 1$ we deduce that $h_{z\i,z,d\i}\ne 0$ hence 
$d\i\le_\cl z\i$ hence $\aa(d\i)\ge \aa(z\i)$ (see P4) and $\aa(d)\ge \aa(z)$.
Assume that $\aa(d)>\aa(z)$, that is, $\aa(d)>N_0$. Let $d'\in\cd$ be such that
$\ga_{d\i,d,d'}\ne 0$ (see P3). By the induction hypothesis applied to $d,d'$ 
instead of $z,d$, we have $\aa(d)=\aa(d')$. From $\ga_{z\i,z,d}\ne 0$,
$\ga_{d\i,d,d'}\ne 0$, $\aa(d)=\aa(d')$, we deduce (using (a)) that 
$\ga_{z,d,z\i}=\ga_{z\i,z,d}$. Hence $\ga_{z,d,z\i}\ne 0$. It follows that
$h_{z,d,z}\ne 0$, hence $z\le_\cl d$, hence $\aa(z)\ge \aa(d)$ (see P4). This
contradicts the assumption $\aa(d)>\aa(z)$. Hence we must have 
$\aa(z)=\aa(d)$, as 
required.

We now prove P7. Assume first that $\ga_{x,y,z}\ne 0$. Let $d\in\cd$ be such 
that $\ga_{z\i,z,d}\ne 0$ (see P3). By (b) we have $\aa(z)=\aa(d)$. Using (a) 
we then have $\ga_{x,y,z}=\ga_{y,z,x}$. Assume next that $\ga_{x,y,z}=0$; we 
must show that $\ga_{y,z,x}=0$. We assume that $\ga_{y,z,x}\ne 0$. By the first
part of the proof, we have

$\ga_{y,z,x}\ne 0\implies\ga_{y,z,x}=\ga_{z,x,y}\ne 0\implies
\ga_{z,x,y}=\ga_{x,y,z}\ne 0$,
\nl
a contradiction.

\subhead 14.8\endsubhead 
{\it P7 $\implies$ P8.} If $\ga_{x,y,z}\ne 0$, then $h_{x,y,z\i}\ne 0$, hence 
$z\i\le_\cl y, z\le_\cl x\i$. By P7 we also have $\ga_{y,z,x}\ne 0$ (hence
$x\i\le_\cl z, x\le_\cl y\i$) and $\ga_{z,x,y}\ne 0$ (hence
$y\i\le_\cl x, y\le_\cl z\i$). Thus, we have 
$x\sim_\cl y\i,y\sim_\cl z\i,z\sim_\cl x\i$.

\subhead 14.9\endsubhead 
{\it 14.3,P4,P8 $\implies$ P9.} We can find a sequence $z'=z_0,z_1,\dots,z_n=z$
such that for any $j\in[1,n]$ we have $z_{j-1}\le_\cl z_j$. By P4 we have 
$\aa(z')=\aa(z_0)\ge \aa(z_1)\ge\dots\ge \aa(z_n)=\aa(z)$. Since 
$\aa(z)=\aa(z')$, we have $\aa(z')=\aa(z_0)=\aa(z_1)=\dots=\aa(z_n)=\aa(z)$. 
Thus, it suffices to show that, if $z'\gets_\cl z$ and $\aa(z')=\aa(z)$, then 
$z'\sim_\cl z$. Let $x,y\in W$ be such that $\ga_{x,y,z\i}\ne 0$. By 14.3, 
there exists $x'\in W$ such that $\pi_{\aa(z)}(h_{x',y,z'})\ne 0$. Since 
$\aa(z')=\aa(z)$, we have $\ga_{x',y,z'{}\i}\ne 0$. From $\ga_{x,y,z\i}\ne 0$,
$\ga_{x',y,z'{}\i}\ne 0$ we deduce, using P8, that $y\sim_\cl z, y\sim_\cl z'$,
hence $z\sim_\cl z'$.

\subhead 14.10\endsubhead 
{\it P9 $\implies$ P10.} We apply P9 to $z\i, z'{}\i$.

\subhead 14.11\endsubhead 
{\it P4,P9,P10 $\implies$ P11.} We can find a sequence $z'=z_0,z_1,\dots,z_n=z$
such that for any $j\in[1,n]$ we have $z_{j-1}\le_\cl z_j$ or 
$z_{j-1}\le_{\Cal R}z_j$. By P4, we have 
$\aa(z')=\aa(z_0)\ge \aa(z_1)\ge\dots\ge \aa(z_n)=\aa(z)$. Since 
$\aa(z)=\aa(z')$, we have 
$\aa(z')=\aa(z_0)=\aa(z_1)=\dots=\aa(z_n)=\aa(z)$. Applying P9 or P10 to 
$z_{j-1},z_j$ we obtain $z_{j-1}\sim_\cl z_j$ or $z_{j-1}\sim_{\Cal R}z_j$. 
Hence $z'\sim_{\cl\Cal R}z$.

\subhead 14.12\endsubhead 
{\it P3,P4,P8 for $W$ and $W_I$ $\implies$ P12.} We write $\aa_I:W_I@>>>\bn$ 
for the $\aa$-function defined in terms of $W_I$. For $x,y,z\in W_I$, we write 
$h^I_{x,y,z},\ga^I_{x,y,z}$ for the analogues of $h_{x,y,z},\ga_{x,y,z}$ when 
$W$ is replaced by $W_I$. Let $\ch_I\sub\ch$ be as in 9.9.

Let $d\in\cd$ be such that $\ga_{y\i,y,d}\ne 0$. (See P3.) Then
$\pi_{\aa(d)}(h_{y\i,y,d\i})\ne 0$. Now $c_{y\i}c_y\in\ch_I$ hence 
$d\in W_I$ and $\pi_{\aa(d)}(h^I_{y\i,y,d\i})\ne 0$. Thus, 
$\aa_I(d\i)\ge \aa(d\i)$. The reverse inequality is obvious hence 
$\aa_I(d)=\aa(d)$. We see that $\ga^I_{y\i,y,d}\ne 0$. By P8 we see that 
$y\sim_\cl d$ (relative to $W_I$) and $y\sim_\cl d$ (relative to $W$). From P4
we deduce that $\aa_I(y)=\aa_I(d)$ and $\aa(y)=\aa(d)$. It follows that 
$\aa(y)=\aa_I(y)$.

\subhead 14.13\endsubhead 
{\it 14.3,P2,P3,P4,P6,P8 $\implies$ P13.} If $x\in\Ga$ then, by P3, there 
exists $d\in\cd$ such that $\ga_{x\i,x,d}\ne 0$. By P8 we have $x\sim_\cl d\i$
hence $d\i\in\Ga$. By P6, we have $d=d\i$ hence $d\in\Ga$. It remains to prove
the uniqueness of $d$. Let $d',d''$ be elements of $\cd\cap\Ga$. We must prove
that $d'=d''$. We can find $x',y',x'',y''$ such that $\ga_{x',y',d'}\ne 0$,
$\ga_{x'',y'',d''}\ne 0$. By P2, we have $x'=y'{}\i,x''=y''{}\i$. By P8, we 
have $y'\sim_\cl d'{}\i=d'$ and $y''\sim_\cl d''{}\i=d''$, hence 
$y',y''\in\Ga$. By the definition of left cells, we can find a sequence 
$y'=x_0,x_1,\dots,x_n=y''$ such that $x_{j-1}\gets_\cl x_j$ for $j\in[1,n]$. 
Since $y'\sim_\cl y''$, we have $x_j\in\Ga$ for all $j$. For each $j\in[1,n-1]$
let $d_j\in\cd$ be such that $\ga_{x_j\i,x_j,d_j}\ne 0$. Let $d_0=d',d_n=d''$.
As in the beginning of the proof, we have $d_j\in\Ga$ for each $j$. Let 
$j\in[1,n]$. By P8, we have $\ga_{x_j,d_j,x_j\i}\ne 0$. Appplying 14.3 to 
$x_j,d_j,x_j,x_{j-1}$ instead of $x,y,z,z'$, we see that there exists $u$ such
that $\pi_{\aa(x_j)}(h_{u,d_j,x_{j-1}})\ne 0$. Since $x_{j-1}\sim x_j$, we 
have $\aa(x_{j-1})=\aa(x_j)$ (see P4), hence 
$\pi_{\aa(x_j)}(h_{u,d_j,x_{j-1}})=\ga_{u,d_j,x_{j-1}\i}\ne 0$. Using P8, we
deduce $\ga_{x_{j-1}\i,u,d_j}\ne 0$. Using P2 we see that $u=x_{j-1}$ and
$\ga_{x_{j-1}\i,x_{j-1},d_j}\ne 0$. We have also
$\ga_{x_{j-1}\i,x_{j-1},d_{j-1}}\ne 0$ and by the uniqueness in P3, it follows
that $d_{j-1}=d_j$. Since this holds for $j\in[1,n]$, it follows that $d'=d''$,
as required.

\subhead 14.14\endsubhead 
{\it P6,P13 $\implies$ P14.} By P13, we can find $d\in\cd$ such that 
$z\sim_\cl d$. Since $d=d\i$ (see P6), it follows that $z\i\sim_{\Cal R}d$. 
Thus, $z\sim_{\cl\Cal R}z\i$.

\subhead 14.15\endsubhead 
In this subsection we reformulate conjecture P15, assuming that P4,P9, P10 
hold. Let $\ti{\ca}=\bz[v,v\i,v',v'{}\i]$ where $v,v'$ are indeterminates. Let 
$\ti{\ch}$ be the free $\ti{\ca}$-module with basis $e_w (w\in W)$. Let 
$(\ch',c'_w,h'_{x,y,z})$ be obtained from $(\ch,c_w,h_{x,y,z})$ by changing the
variable $v$ to $v'$. In particular, $h'_{x,y,z}\in\bz[v',v'{}\i]$.

On $\ti{\ch}$ we have a left $\ch$-module structure given by
$v^nc_ye_w=\sum_xv^nh_{y,w,x}e_x$ and a right $\ch'$-module structure defined 
by $e_w(v'{}^nc'_y)=\sum_xv'{}^nh'_{w,y,x}e_x$. These module structures do not
commute in general. For each $a\ge 0$ let $\ti{\ch}_{\ge a}$ be the 
$\ti{\ca}$-submodule of $\ti{\ch}$ spanned by $\{e_w|\aa(w)\ge a\}$. By P4, 
this
is a left $\ch$-submodule and a right $\ch'$-submodule of $\ti{\ch}$. We have

$\dots\ti{\ch}_{\ge 2}\sub\ti{\ch}_{\ge 1}\sub\ti{\ch}_{\ge 0}=\ti{\ch}$
\nl
and $gr\ti{\ch}=\opl_{a\ge 0}\ti{\ch}_{\ge a}/\ti{\ch}_{\ge a+1}$ inherits a
left $\ch$-module structure and a right $\ch'$-module structure from 
$\ti{\ch}$. Clearly, P15 is equivalent to the condition that these module
stuctures on $gr\ti{\ch}$ commute. To check this last condition, it is enough 
to check that the actions of $c_s, c'_{s'}$ commute on $gr\ti{\ch}$ whenever 
$s,s'\in S$. Let $s,s'\in S, w\in W$. A computation using 6.6, 6.7, 8.2 shows 
that $(c_se_w)c'_{s'}-c_s(e_wc_{s'})$ is $0$ if $sw<w$ or $ws'<w$, while if
$sw>w,ws'>w$, it is
$$\sum_{y;sy<y,ys'<y}(h'_{w,s',y}(v_s+v_s\i)-h_{s,w,y}(v'_{s'}+v'_{s'}{}\i))e_y
+\sum_{y;sy<y,ys'<y}\al_ye_y$$
where 
$$\al_y=\sum_{y';y's'<y'<sy'}h'_{w,s',y'}h_{s,y',y}
-\sum_{y';sy'<y'<y's'}h_{s,w,y'}h'_{y',s',y}.$$
If $y$ satisfies $sy<y,ys'<y$ and either $h'_{w,s',y}$ or $h_{s,w,y}$ is 
$\ne 0$, then $\aa(y)>\aa(w)$. (We certainly have $\aa(y)\ge \aa(w)$ by P4. 
If we had 
$\aa(y)=\aa(w)$ and $h_{s,w,y}\ne 0$ then by P9 we would have $y\sim_\cl w$ 
hence
$\Cal R(y)=\Cal R(w)$ contradicting $ys'<y,ws'>w$. If we had $\aa(y)=\aa(w)$ 
and 
$h'_{s,w,y}\ne 0$ then by P10 we would have $y\sim_{\Cal R}w$ hence 
$\cl(y)=\cl(w)$ contradicting $sy<y,sw>w$.) Hence, if $sw>w,ws'>w$, we have
$$(c_se_w)c'_{s'}-c_s(e_wc'_{s'})=\sum_{y;sy<y,ys'<y, \aa(y)=\aa(w)}\al_ye_y 
\mod\ti{\ch}_{\ge \aa(w)+1}.$$
We see that P15 is equivalent to the following statement.

(a) {\it If $y,w\in W,s,s'\in S$ are such that 
$sw>w,ws'>w,sy<y,ys'<y,\aa(y)=\aa(w)$, then}

$\sum_{y';y's'<y'<sy'}h'_{w,s',y'}h_{s,y',y}
=\sum_{y';sy'<y'<y's'}h_{s,w,y'}h'_{y',s',y}$.

\head 15. Example: the split case\endhead
\subhead 15.1\endsubhead
{\it In this section we assume that we are in the split case (see 3.1), that 
$W$ is tame and that}

(a) $h_{x,y,z}\in\bn[v,v\i]$ for all $x,y,z$ in $W$,

(b) $p_{y,w}\in\bn[v\i]$ for all $y,w$ in $W$.
\nl
It is known that (a),(b) hold automatically in the split case if $W$ is tame,
integral.

Under these assumptions we will show that 14.3 and P1-P3 hold for $W,l$ hence 
all of P1-P14 hold for $W,l$; we will also show that P15 holds.

A reference for this section is \cite{\LCC}.

\subhead 15.2\endsubhead 
{\it Proof of P1.} The weaker inequality $\aa(z)\le l(z)$ was first proved in 
\cite{\LC} for finite $W$ and then for general $W$ by Springer (unpublished). 
The present argument was inspired by Springer's argument.

From 14.5(a) we see that for $x,y\in W$ we have 

(a) $\sum_zh_{x,y,z}p_{1,z}\in\bz[v\i]$.
\nl
From 15.1(a),(b) we see that $h_{x,y,z}p_{1,z}\in\bn[v,v\i]$ for any $z$. Hence
in (a) there are no cancellations, so that

(b) $h_{x,y,z}p_{1,z}\in\bn[v\i]$ for any $z$.
\nl
We now fix $z$ and choose $x,y$ so that $\ga_{x,y,z\i}\ne 0$. From the 
definitions we have

(c) $h_{x,y,z}p_{1,z}\in\ga_{x,y,z\i}n_zv^{\aa(z)-\De(z)}+
\text{strictly smaller
powers of } v$
\nl
and the coefficient of $v^{\aa(z)-\De(z)}$ is $\ne 0$. Comparing with (b) we
deduce that $\aa(z)-\De(z)\le 0$.

\subhead 15.3\endsubhead 
{\it Proof of P2.} Assume that $x\ne y\i$. From 14.5(a) we see that 

(a) $\sum_zh_{x,y,z}p_{1,z}\in v\i\bz[v\i]$.
\nl
As in 15.2, this implies (using 15.1(a),(b)) that 

(b) $h_{x,y,z}p_{1,z}\in v\i\bn[v\i]$ for any $z$.
\nl
Assume now that $z=d\i\in\cd$. The equality 15.2(c) becomes in our case

$h_{x,y,z}p_{1,z}\in\ga_{x,y,z\i}n_z+v\i\bz[v\i]$.
\nl
Comparing with (b) we deduce that $\ga_{x,y,z\i}n_z=0$. Since $n_z\ne 0$, we
have $\ga_{x,y,z\i}=0$. This proves P2.

\subhead 15.4\endsubhead 
{\it Proof of P3.} From 14.5(a) we see that 

(a) $\sum_zh_{y\i,y,z}p_{1,z}\in 1+v\i\bz[v\i]$.
\nl
As in 15.2, this implies (using 15.1(a),(b)) that there is a unique $z$, say
$z=d\i$ such that

(b) $h_{y\i,y,d\i}p_{1,d\i}\in 1+v\i\bn[v\i]$ 
\nl
and that

(c) $h_{y\i,y,z}p_{1,z}\in v\i\bn[v\i]$
\nl
for all $z\ne d\i$. For $z=d\i$, the equality 15.2(c) becomes

$h_{y\i,y,d\i}p_{1,d\i}\in\ga_{y\i,y,d}n_{d\i}v^{\aa(d)-\De(d)}
+\text{strictly smaller powers of } v$.
\nl
Here $\aa(d)-\De(d)\le 0$. Comparing with (b) we deduce that $\aa(d)-\De(d)=0$
and $\ga_{y\i,y,d}n_{d\i}=1$. Thus, $d\in\cd$ and $\ga_{y\i,y,d}\ne 0$. Thus, 
the existence part of P3 is established.

Assume that there exists $d'\ne d$ such that $d'\in\cd$ and 
$\ga_{y\i,y,d'}\ne 0$. For $z=d'{}\i$ the equality 15.2(c) becomes

$h_{y\i,y,d'{}\i}p_{1,d'{}\i}\in\ga_{y\i,y,d'}n_{d'{}\i}+v\i\bz[v\i]$.
\nl
Comparing with (c) (with $z=d'{}\i$) we deduce that 
$\ga_{y\i,y,d'}n_{d'{}\i}=0$ hence $\ga_{y\i,y,d'}=0$, a contradiction. This
proves the uniqueness part of P3.

\subhead 15.5\endsubhead
{\it Proof of 14.3.} We may assume that $z'\ne z$. Then we can find $s\in S$ 
such that $sz'<z',sz>z$ and $h_{s,z,z'}\ne 0$. Since $h_{x,y,z}\ne 0$, we have
(by 13.1(d)) $z\le_{\Cal R}x$ hence $\cl(x)\sub\cl(z)$ (by 8.6). Since 
$s\notin\cl(z)$, we have $s\notin\cl(x)$, that is, $sx>x$. We have 
$c_sc_xc_y=\sum_up_uc_u$, where

$p_u=\sum_wh_{x,y,w}h_{s,w,u}=\sum_{x'}h_{s,x,x'}h_{x',y,u}$.
\nl
In particular, 

$p_{z'}=\sum_wh_{x,y,w}h_{s,w,z'}=h_{x,y,z}h_{s,z,z'}
+\sum_{w;w\ne z}h_{x,y,w}h_{s,w,z'}$.
\nl
By 6.5, we have $h_{s,z,z'}\in\bz$ hence, taking 

(a) $\pi_n(p_{z'})=
\pi_n(h_{x,y,z})h_{s,z,z'}+\sum_{w;w\ne z}\pi_n(h_{x,y,w}h_{s,w,z'})$.
\nl
for any $n\in\bz$. In particular, this holds for $n=\aa(z)$. By assumption, we 
have $\pi_n(h_{x,y,z})\ne 0$ and $h_{s,z,z'}\ne 0$; hence, by 15.1(a), we have
$\pi_n(h_{x,y,z})>0$ and $h_{s,z,z'}>0$. Again, by 15.1(a) we have 
$\pi_n(h_{x,y,w}h_{s,w,z'})\ge 0$ for any $w\ne z$. Hence from (a) we deduce 
$\pi_n(p_{z'})>0$. Since $p_{z'}=\sum_{x'}h_{s,x,x'}h_{x',y,z'}$, there exists
$x'$ such that $\pi_n(h_{s,x,x'}h_{x',y,z'})\ne 0$. Since $sx>x$, we see from
6.5 that $h_{s,x,x'}\in\bz$ hence 

$\pi_n(h_{s,x,x'}h_{x',y,z'})=h_{s,x,x'}\pi_n(h_{x',y,z'})$.
\nl
Thus we have $\pi_n(h_{x',y,z'})\ne 0$. This proves 14.3 in our case.

\subhead 15.6\endsubhead 
Since 14.3 and P1-P3 are known, we see that P1-P11 and P13,P14 hold in our case
(see \S14). The same arguments can be applied to $W_I$ where $I\sub W$, hence 
P1-P11 and P13,P14 hold for $W_I$. By 14.12, P12 holds for $W$. Thus, P1-P14 
hold for $W$.

\subhead 15.7\endsubhead 
{\it Proof of P15.} By 14.15, we see that it is enough to prove 14.15(a). Let 
$y,w,s,s'$ be as in 14.15(a). In our case, by 6.5, the equation in 14.15(a) 
involves only integers, hence it is enough to prove it after specializing 
$v=v'$. If in 14.15 we specialize $v=v'$, then the left and right module 
structures in 14.15 clearly commute, since the left and right regular 
representations of $\ch$ commute. Hence the coefficient of $e_y$ in 
$((c_se_w)c'_{s'}-c_s(e_wc_{s'}))_{v=v'}$ is $0$. By the computation in 14.15,
this coefficient is 
$$(h_{w,s',y}-h_{s,w,y})(v+v\i)+\sum\Sb y'\\y's'<y'<sy'\endSb h_{w,s',y'}
h_{s,y',y}-\sum\Sb y'\\sy'<y'<y's'\endSb h_{s,w,y'}h_{y',s',y}=0\tag a$$
By 6.5, $h_{s,w,y}$ is the coefficient of $v\i$ in $p_{y,w}$ and
$h_{w,s',y}=h_{s',w\i,y\i}$ is the coefficient of $v\i$ in 
$p_{y\i,w\i}=p_{y,w}$. Thus, $h_{s,w,y}=h_{w,s',y}$ and (a) reduces to the 
equation in 14.15(a) (specialized at $v=v'$). This proves 14.15(a).

\head 16. Example: the quasisplit case\endhead
\subhead 16.1 \endsubhead
Let $(\ti W,\ti S)$ be a Coxeter group and let $\si:\ti W@>>>\ti W$ be an
automorphism of finite order of $\ti W$ which restricts to a permutation of
$\ti S$ such that any $\si$-orbit on $\ti S$ generates a finite subgroup of
$\ti W$. 
Assume also that the following condition is satisfied: if $s,s'\in S$ and
$k\in\bn$ are such that $m_{s,s'}\ge 4$ and $\si^k$ maps $\{s,s'\}$ into
itself then $\si^k(s)=s,\si^k(s')=s'$.

Let $\ti W^\si=\{w\in\ti W|\si(w)=w\}$. For any $\si$-orbit $\boo$ on
$\ti S$, we set $s_\boo=w_0^\boo\in\ti W^\si$. The elements $s_\boo$ for 
various $\boo$ as above form a subset $\ti S_\si$ of $\ti W^\si$. 

\proclaim{Lemma 16.2} $(\ti W^\si,\ti S_\si)$ is a Coxeter group and the 
restriction to $\ti W^\si$ of the length function $\ti l:\ti W@>>>\bn$ is a 
weight function $L:\ti W^\si@>>>\bn$.
\endproclaim
We omit the proof.

\subhead 16.3\endsubhead
{\it In this section we assume that $\ti W$ is tame, $(W,S)$ is the Coxeter 
group $(\ti W^\si,\ti S_\si)$ and that $L:W@>>>\bn$ is as in 16.2.} We then 
say that we are in the {\it quasisplit case}.

We denote $h_{x,y,z},p_{x,y},\aa(z),\ga_{x,y,z},\De(z),\cd$ defined in terms of
$\ti W,\ti l$ by 

$\ti h_{x,y,z},\ti p_{x,y},\ti\aa(z),\ti\ga_{x,y,z},\ti\De(z),\ti{\cd}$.
\nl
We shall assume that the following holds.

(a)  {\it For any $x,y,z\in\ti W$ and any integer 
$n$ there exists a $\bq$-vector space $V_{x,y,z}^n$ such that 
$\ti h_{x,y,z}=\sum_n\dim V_{x,y,z}^nv^n$; for any $x,y\in\ti W$ and any 
integer $n\le 0$ there exists a $\bq$-vector space $V_{x,y}^n$ such that}
$\ti p_{x,y}=\sum_{n\le 0}\dim V_{x,y}^nv^n$.

In other words, $(\ti W,\ti l)$ satisfies the requirements of 15.1. We shall
further assume that the following holds.

(b) {\it For any $x,y,z\in W$ and any integer $n$, $V_{x,y,z}^n$ carries a 
linear transformation $\si$ of finite order and 
$h_{x,y,z}=\sum_n\tr(\si,V_{x,y,z}^n)v^n$; for any $x,y\in W$ and any integer 
$n\le 0$, $V_{x,y}^n$ carries a linear transformation $\si$ of finite order and
$p_{x,y}=\sum_n\tr(\si,V_{x,y}^n)v^n$.}
\nl
It is known that this is automatically satisfied if $\ti W$ is tame, integral.

We deduce that (c)-(f) below hold: 

(c) If $x,y,z\in W,n\in\bz$ and $\pi_n(h_{x,y,z})\ne 0$ then 
$\pi_n(\ti h_{x,y,z})\ne 0$.

(d) If $x,y\in W,n\in\bz$ and $\pi_n(p_{x,y})\ne 0$ then 
$\pi_n(\ti p_{x,y})\ne 0$.

(e) If $x,y,z\in W,n\in\bz$ and $\pi_n(\ti h_{x,y,z})=\pm 1$ then
$\pi_n(h_{x,y,z})=\pm 1$.

(f) If $x,y\in W,n\in\bz$ and $\pi_n(\ti p_{x,y})=\pm 1$ then 
$\pi_n(p_{x,y})=\pm 1$.

\subhead 16.4\endsubhead
By our assumptions, the results of \S15 are applicable to $\ti W,\ti l$. Under
the assumptions above, we will show that 14.3 and P1-P3 hold for $W,L$ hence 
all of P1-P14 hold for $W,L$. 

\proclaim{Lemma 16.5} For $z\in W$ we have $\aa(z)=\ti\aa(z)$ and 
$\ti\De(z)\le\De(z)$.
\endproclaim
We can find $x,y\in W$ such that $\pi_{\aa(z)}(h_{x,y,z})\ne 0$. By 16.3(c) we 
have \linebreak 
$\pi_{\aa(z)}(\ti h_{x,y,z})\ne 0$. Hence $\aa(z)\le\ti\aa(z)$. By P3,P5 for 
$\ti W$, there is a unique $d\in\ti{\cd}$ such that $\ti\ga_{z\i,z,d}=\pm 1$. 
The uniqueness of $d$ implies that $d$ is fixed by $\si$. Thus $d\in W$. By P7
for $\ti W$, we have $\ti\ga_{z,d,z\i}=\pm 1$. Hence 
$\pi_{\ti\aa(z)}(\ti h_{z,d,z})=\pm 1$. By 16.3(e), we have
$\pi_{\ti\aa(z)}(h_{z,d,z})=\pm 1$. Hence $\ti\aa(z)\le\aa(z)$ so that 
$\ti\aa(z)=\aa(z)$.

By definition, we have $\pi_{-\De(z)}(p_{1,z})\ne 0$. Using 16.3(d), we deduce
that \linebreak $\pi_{-\De(z)}(\ti p_{1,z})\ne 0$. Hence 
$-\De(z)\le-\ti\De(z)$. The lemma is proved.

\proclaim{Lemma 16.6} $\cd=\ti\cd\cap W$.
\endproclaim
Let $d\in\cd$. We have $\aa(d)=\De(d)$. Using 16.5, we deduce 
$\ti\aa(d)=\De(d)$. 
By P1 for $\ti W$, we have $\ti\aa(d)\le\ti\De(d)$. Hence 
$\De(d)\le\ti\De(d)$. 
Using 16.5, we deduce $\De(d)=\ti\De(d)$ so that $\ti\De(d)=\ti\aa(d)$ and 
$d\in\ti\cd$.

Conversely, let $d\in\ti\cd\cap W$. We have $\ti\aa(d)=\ti\De(d)$. Using 16.5 
we
deduce $\aa(d)=\ti\De(d)$. By P5 for $\ti W$, we have 
$\pi_{-\ti\De(d)}(\ti p_{1,d})=\pm 1$. Using 16.3(f) we deduce 
$\pi_{-\ti\De(d)}(p_{1,d})=\pm 1$. Hence $-\ti\De(d)\le-\De(d)$. Using 16.5 we 
deduce $\De(d)=\ti\De(d)$ so that $\De(d)=\aa(d)$ and $d\in\cd$. The lemma is 
proved.

\proclaim{Lemma 16.7} (a) Let $x,y,z\in W$ be such that $\ga_{x,y,z}\ne 0$. 
Then $\ti\ga_{x,y,z}\ne 0$. 

(b) Let $x,y,z\in W$ be such that $\ti\ga_{x,y,z}=\pm 1$. Then 
$\ga_{x,y,z}=\pm 1$.
\endproclaim
In the setup of (a) we have $\pi_{\aa(z\i)}(h_{x,y,z\i})\ne 0$. Using 16.5 we 
deduce that\linebreak $\pi_{\ti\aa(z\i)}(h_{x,y,z\i})\ne 0$. Using 16.3(c), we 
deduce that $\pi_{\ti\aa(z\i)}(\ti h_{x,y,z\i})\ne 0$. Hence 
$\ti\ga_{x,y,z}\ne 0$. 

In the setup of (b) we have $\pi_{\ti\aa(z\i)}(\ti h_{x,y,z\i})=\pm 1$. Using 
16.5, we deduce $\pi_{\aa(z\i)}(\ti h_{x,y,z\i})=\pm 1$. 
Using 16.3(e), we deduce
$\pi_{\aa(z\i)}(h_{x,y,z\i})=\pm 1$. Hence $\ga_{x,y,z}=\pm 1$. 

\subhead 16.8\endsubhead 
{\it Proof of P1.} By 16.5 and P1 for $\ti W$, we have 
$\aa(z)=\ti\aa(z)\le\ti\De(z)\le\De(z)$, hence $\aa(z)\le\De(z)$.

\subhead 16.9\endsubhead 
{\it Proof of P2.} In the setup of P2, we have (by 16.7) $\ti\ga_{x,y,d}\ne 0$
and $d\in\ti\cd$ (see 16.6). Using P2 for $\ti W$, we deduce $x=y\i$.

\subhead 16.10\endsubhead 
{\it Proof of P3.} Let $y\in W$. By P3 for $\ti W$, there is a unique 
$d\in\ti{\cd}$ such that $\ti\ga_{y\i,y,d}\ne 0$. By the uniqueness of $d$, we
have $\si(d)=d$ hence $d\in W$. Using P5 for $\ti W$, we see that
$\ti\ga_{y\i,y,d}=\pm 1$. Using 16.7, we deduce $\ga_{y\i,y,d}=\pm 1$. Since 
$d\in\cd$ by 16.6, the existence part of P3 is established. Assume now that 
$d'\in\cd$ satisfies $\ga_{y\i,y,d'}\ne 0$. Using 16.7, we deduce 
$\ti\ga_{y\i,y,d'}\ne 0$. Since $d'\in\ti{\cd}$ by 16.6, we can use the 
uniqueness in P3 for $\ti W$ to deduce that $d=d'$. Thus P3 holds for $W$.

\subhead 16.11\endsubhead 
{\it Proof of P4.} We may assume that there exists $s\in S$ such that 
$h_{s,z,z'}\ne 0$ or $h_{z,s,z'}\ne 0$. In the first case, using 16.3(c), we 
deduce $\ti h_{s,z,z'}\ne 0$. Hence $z'\le_\cl z$ (in $\ti W$) and using P4 for
$\ti W$, we deduce that $\ti\aa(z')\ge\ti\aa(z)$. Using now 16.5, we see that
$\aa(z')\ge\aa(z)$. The proof in the second case is entirely similar.

\subhead 16.12\endsubhead 
Now P5 is proved as in 14.5; P6 is proved as in 14.6; P7 is proved as in 14.7;
P8 is proved as in 14.8; P12 is proved as in 14.12.

\subhead 16.13\endsubhead 
{\it Proof of P13.} If $z'\gets_\cl z$ in $W$, then there exists $s\in S$ such
that $h_{s,z,z'}\ne 0$ hence, by 16.3(c), $\ti h_{s,z,z'}\ne 0$, hence 
$z'\le_\cl z$ in $\ti W$. It follows that 

(a) $z'\le_\cl z$ (in $W$) implies $z'\le_\cl z$ (in $\ti W$). 
\nl
Hence 

(b) $z'\sim_\cl z$ (in $W$) implies $z'\sim_\cl z$ (in $\ti W$). 
\nl
Thus any left cell of $W$ is contained in a left cell of $\ti W$. 

In the setup of P13, let $\ti\Ga$ be the left cell of $\ti W$ containing $\Ga$.
Let $x\in\Ga$. By P3 for $W$, there exists $d\in\cd$ such that 
$\ga_{x\i,x,d}\ne 0$. By P8 for $W$, we have $x\sim_\cl d\i$ hence $d\i\in\Ga$.
Using P6 we have $d=d\i$, hence $d\in\Ga$. It remains to prove the uniqueness 
of $d$. Let $d',d''$ be elements of $\cd\cap\Ga$. We must prove that $d'=d''$.
Now $d',d''$ belong to $\ti\Ga$ and, by 16.6, are in $\ti{\cd}$. Using P13 for
$\ti W$, it follows that $d'=d''$. Thus P13 holds for $W$.

\proclaim{Lemma 16.14} Let $x,y\in W$. We have $x\sim_\cl y$ (in $W$) if and 
only if $x\sim_\cl y$ (in $\ti W$).
\endproclaim
If $x\sim_\cl y$ (in $W$) then $x\sim_\cl y$ (in $\ti W$), by 16.13(b). 

Assume now that $x\sim_\cl y$ (in $\ti W$). Let $d,d'\in\cd$ be such that 
$x\sim_\cl d$ (in $W$) and $y\sim_\cl d'$ (in $W$); see P13. By the first line
of the proof we have $x\sim_\cl d$ (in $\ti W$) and $y\sim_\cl d'$ (in 
$\ti W$). Hence $d\sim_\cl d'$ (in $\ti W$). Since $d,d'\in\ti{\cd}$, we deduce
(using P13 for $\ti W$) that $d=d'$. It follows that $x\sim_\cl y$ (in $W$). 
The lemma is proved.

\subhead 16.15\endsubhead 
{\it Proof of P9.} We assume that $z'\le_\cl z$ (in $W$) and $\aa(z')=\aa(z)$.
 By
16.13(a), it follows that $z'\le_\cl z$ (in $\ti W$) and, using 16.5, that 
$\ti\aa(z')=\ti\aa(z)$. Using now P9 in $\ti W$, it follows that $z'\sim_\cl z$
(in $\ti W$). Using 16.14, we deduce that $z'\sim_\cl z$ (in $W$).

\subhead 16.16\endsubhead 
Now P10 is proved as in 14.10; P11 is proved as in 14.11; P14 is proved as in 
14.14.

\subhead 16.17\endsubhead 
Assuming that $(\ti W,\ti S)$ is tame, integral, P15 can be shown to hold in 
our case. We sketch a proof which is almost (but not entirely) correct.

The proof of P15 given in 14.15,15.7, can be refined to a proof of the 
following statement:

{\it For any $w,y,x,x'$ in $\ti W$ and any $k$ there exists a natural 
isomorphism of vector spaces

$\opl_{j+j'=k}\opl_{y'\in\ti W}V^{j'}_{w,x',y'}\ot V^j_{x,y',y}@>\sim>>
\opl_{j+j'=k}\opl_{y'\in\ti W}V^j_{x,w,y'}\ot V^{j'}_{y',x',y}$.
\nl
When $\ti\aa(w)=\ti\aa(y)$, this restricts to an isomorphism 

$\opl_{y'\in\ti W}V^{j'}_{w,x',y'}\ot V^j_{x,y',y}@>\sim>>
\opl_{y'\in\ti W}V^j_{x,w,y'}\ot V^{j'}_{y',x',y}$
\nl
for any $j,j'$ such that $j+j'=k$.}
\nl
Assuming now that $w,y,x,x'\in W$ and taking traces of $\si$ in both sides, we 
deduce

$\sum_{y'\in W}\pi_{j'}(h_{w,x',y'})\pi_j(h_{x,y',y})=
\sum_{y'\in W}\pi_j(h_{x,w,y'})\pi_{j'}(h_{y',x',y})$
\nl
(the summands corresponding to $y'\in\ti W-W$ do not contibute to the trace) or
equivalently

$\sum_{y'\in W}h'_{w,x',y'}h_{x,y',y}=\sum_{y'\in W}h_{x,w,y'}h'_{y',x',y}$,
\nl
as required. (The actual proof is slightly more complicated since the direct
sum decompositions are not natural, only certain filtrations attached to them 
are.)

\head 17. Example: the infinite dihedral case\endhead
\subhead 17.1\endsubhead
In this section we preserve the setup of 7.1. We assume that $m=\infty$ and 
that $L_2>L_1$. We will show that P1-P15 hold in this case.

Let $\ze=v^{L_2-L_1}+v^{L_1-L_2}$. For $a\in\{1,2\}$, let 
$f_a=v^{L_a}+v^{-L_a}$. For $m,n\in\bz$ we define $\de_{m<n}$ to be $1$ if
$m<n$ and to be $0$ otherwise.

\subhead 17.2\endsubhead
From 7.5, 7.6 we have for all $k'\in\bn$:

$c_1c_{2_{k'}}=c_{1_{k'+1}}$,

$c_2c_{1_{k'}}=c_{2_{k'+1}}+\de_{k'>1}\ze c_{2_{k'-1}}+\de_{k'>3}c_{2_{k'-3}}$.

\proclaim{Proposition 17.3} For $k\ge 0,k'\ge 1$ we have

(a) $c_{2_{2k+1}}c_{2_{k'}}=f_2\sum_{u\in[0,k];2u\le k'-1}c_{2_{2k+k'-4u}}$,

(b) $c_{1_{2k+2}}c_{2_{k'}}=f_2\sum_{u\in[0,k];2u\le k'-1}c_{1_{2k+k'+1-4u}}$.
\endproclaim
Assume that $k=0$. Using 17.2 we have $c_2c_{2_{k'}}=f_2c_{2_{k'}}$.

Assume now that $k=1$. Using 17.2, we have $c_{2_3}=c_2c_1c_2-\ze c_2$. Using 
this and 17.2, we have
$$\align&c_{2_3}c_{2_{k'}}=c_2c_1c_2c_{2_{k'}}-\ze c_2c_{2_{k'}}=
f_2c_2c_{1_{k'+1}}-f_2\ze c_{2_{k'}}\\&=
f_2c_{2_{k'+2}}+f_2\ze c_{2_{k'}}+\de_{k'>2}f_2c_{2_{k'-2}}
-f_2\ze c_{2_{k'}}=f_2c_{2_{k'+2}}+\de_{k'>2}f_2c_{2_{k'-2}},\endalign$$
as required. We prove the equality in (a) for fixed $k'$, by induction on $k$.
The cases $k=0,1$ are already known. If $k=2$ then using 17.2, we have
$c_{2_5}=c_2c_1c_{2_3}-\ze c_{2_3}-c_{2_1}$. Using this, 17.2, and the 
induction hypothesis, we have
$$\align&c_{2_5}c_{2_{k'}}=c_2c_1c_{2_3}c_{2_{k'}}-\ze c_{2_3}c_{2_{k'}}
-c_{2_1}c_{2_{k'}}\\&=f_2c_2c_1c_{2_{k'+2}}+\de_{k'>2}f_2c_2c_1c_{2_{k'-2}}
-\ze f_2c_{2_{k'+2}}-\de_{k'>2}\ze f_2c_{2_{k'-2}}-f_2c_{2_{k'}}\\&
=f_2c_2c_{1_{k'+3}}+\de_{k'>2}f_2c_2c_{1_{k'-1}}-\ze f_2c_{2_{k'+2}}-\de_{k'>2}
\ze f_2c_{2_{k'-2}}-f_2c_{2_{k'}}\\&=f_2c_{2_{k'+4}}+f_2\ze c_{2_{k'+2}}+f_2
c_{2_{k'}}+\de_{k'>2}f_2c_{2_{k'}}+\de_{k'>2}f_2\ze c_{2_{k'-2}}+\de_{k'>4}f_2
c_{2_{k'-4}}\\&-\ze f_2c_{2_{k'+2}}-\de_{k'>2}\ze f_2c_{2_{k'-2}}-f_2c_{2_{k'}}
=f_2c_{2_{k'+4}}+\de_{k'>2}f_2c_{2_{k'}}+\de_{k'>4}f_2c_{2_{k'-4}},\endalign$$
as required. A similar argument applies for $k\ge 3$. This proves (a).

(b) is obtained by multiplying both sides of (a) by $c_1$ on the left. The 
proposition is proved.

\proclaim{Proposition 17.4} For $k\ge 0,k'\ge 1$, we have 

(a) $c_{2_{2k+1}}c_{1_{k'}}=\sum_{u\in[0,2k+2]}p_uc_{2_{k'+2k+1-2u}}$,

(b) $c_{1_{2k+2}}c_{1_{k'}}=\sum_{u\in[0,2k+2]}p_uc_{1_{k'+2k+2-2u}}$,

(c) $c_{1_{k'}\i}c_{2_{2k+1}}=\sum_{u\in[0,2k+2]}p_uc_{2_{k'+2k+1-2u}\i}$,

(d) $c_{1_{k'}\i}c_{1_{2k+2}\i}=\sum_{u\in[0,2k+2]}p_uc_{1_{k'+2k+1-2u}\i}$,

(e) $c_{2_{2k+2}}c_{1_{k'}}=\sum_{u\in[0,2k+2]}f_1p_uc_{2_{k'+2k+1-2u}}$,    

(f) $c_{1_{2k+3}}c_{1_{k'}}=\sum_{u\in[0,2k+2]}f_1p_uc_{1_{k'+2k+2-2u}}$,

(g) $c_{1_1}c_{1_{k'}}=f_1c_{1_{k'}}$,
\nl
where 

$p_0=1$, $p_{2k+2}=\de_{k'>2k+3}$,

$p_u=\de_{k'>u}\ze$ for $u=1,3,5,\dots,2k+1$,

$p_u=\de_{k'>u-1}+\de_{k'>u+1}$ for $u=2,4,6,\dots,2k$.
\endproclaim
We prove (a). For $k=0$ the equality in (a) is
$c_2c_{1_{k'}}=c_{2_{k'+1}}+\de_{k'>1}\ze c_{2_{k'-1}}+\de_{k'>3}c_{2_{k'-3}}$
which is contained in 17.2. Assume now that $k=1$. Using 
$c_{2_3}=c_2c_1c_2-\ze c_2$ and 17.2, we have
$$\align&c_{2_3}c_{1_{k'}}=c_2c_1c_2c_{1_{k'}}-\ze c_2c_{1_{k'}}=c_2c_1
c_{2_{k'+1}}+\de_{k'>1}\ze c_2c_1c_{2_{k'-1}}+\de_{k'>3}c_2c_1c_{2_{k'-3}}\\&-
\ze c_{2_{k'+1}}-\de_{k'>1}\ze^2 c_{2_{k'-1}}-\de_{k'>3}\ze c_{2_{k'-3}}=c_2
c_{1_{k'+2}}+\de_{k'>1}\ze c_2c_{1_{k'}}\\&+\de_{k'>3}c_2c_{1_{k'-2}}-\ze 
c_{2_{k'+1}}-\de_{k'>1}\ze^2 c_{2_{k'-1}}-\de_{k'>3}\ze c_{2_{k'-3}}\\&=
c_{2_{k'+3}}+\ze c_{2_{k'+1}}+\de_{k'>1}c_{2_{k'-1}}+\de_{k'>1}\ze c_{2_{k'+1}}
+\de_{k'>1}\ze^2c_{2_{k'-1}}+\de_{k'>3}\ze c_{2_{k'-3}}\\&+\de_{k'>3}
c_{2_{k'-1}}+\de_{k'>3}\ze c_{2_{k'-3}}+\de_{k'>5}c_{2_{k'-5}}-\ze c_{2_{k'+1}}
-\de_{k'>1}\ze^2c_{2_{k'-1}}-\de_{k'>3}\ze c_{2_{k'-3}}\\&=c_{2_{k'+3}}+
\de_{k'>1}\ze c_{2_{k'+1}}+(\de_{k'>1}+\de_{k'>3})c_{2_{k'-1}}+\de_{k'>3}\ze 
c_{2_{k'-3}}+\de_{k'>5}c_{2_{k'-5}},\endalign$$
as required.

We prove the equality in (a) for fixed $k'$, by induction on $k$. The cases
$k=0,1$ are already known. Assume now that $k=2$. Using 
$c_{2_5}=c_2c_1c_{2_3}-\ze c_{2_3}-c_{2_1}$, 17.2, and the case $k=1$, we have
$$\align&c_{2_5}c_{1_{k'}}=c_2c_1c_{2_3}c_{1_{k'}}-\ze c_{2_3}c_{1_{k'}}
-c_{2_1}c_{1_{k'}}\\&=c_2c_1c_{2_{k'+3}}+\de_{k'>1}\ze c_2c_1c_{2_{k'+1}}+
(\de_{k'>1}+\de_{k'>3})c_2c_1c_{2_{k'-1}}+\de_{k'>3}\ze c_2c_1c_{2_{k'-3}}\\&
+\de_{k'>5}c_2c_1c_{2_{k'-5}}-\ze c_{2_{k'+3}}-\de_{k'>1}\ze^2 c_{2_{k'+1}}-
(\de_{k'>1}+\de_{k'>3})\ze c_{2_{k'-1}}\\&-\de_{k'>3}\ze^2 c_{2_{k'-3}}-
\de_{k'>5}\ze c_{2_{k'-5}}-c_{2_{k'+1}}-\de_{k'>1}\ze c_{2_{k'-1}}-\de_{k'>3}
c_{2_{k'-3}}\\&=c_2c_{1_{k'+4}}+\de_{k'>1}\ze c_2c_{1_{k'+2}}+(\de_{k'>1}+
\de_{k'>3})c_2c_{1_{k'}}+\de_{k'>3}\ze c_2c_{1_{k'-2}}\\&+\de_{k'>5}c_2
c_{1_{k'-4}}-\ze c_{2_{k'+3}}-\de_{k'>1}\ze^2 c_{2_{k'+1}}-(\de_{k'>1}+
\de_{k'>3})\ze c_{2_{k'-1}}\\&-\de_{k'>3}\ze^2 c_{2_{k'-3}}-\de_{k'>5}\ze 
c_{2_{k'-5}}-c_{2_{k'+1}}-\de_{k'>1}\ze c_{2_{k'-1}}-\de_{k'>3}c_{2_{k'-3}}\\&=
c_{2_{k'+5}}+\ze c_{2_{k'+3}}+c_{2_{k'+1}}+\de_{k'>1}\ze c_{2_{k'+3}}+
\de_{k'>1}\ze^2 c_{2_{k'+1}}+\de_{k'>1}\ze c_{2_{k'-1}}\\&+(\de_{k'>1}+
\de_{k'>3})c_{2_{k'+1}}+(\de_{k'>1}+\de_{k'>3})\ze c_{2_{k'-1}}+2\de_{k'>3}
c_{2_{k'-3}}+\de_{k'>3}\ze c_{2_{k'-1}}\\&+\de_{k'>3}\ze^2c_{2_{k'-3}}+
\de_{k'>5}\ze c_{2_{k'-5}}+\de_{k'>5}c_{2_{k'-3}}+\de_{k'>5}\ze c_{2_{k'-5}}+
\de_{k'>7}c_{2_{k'-7}}\\&-\ze c_{2_{k'+3}}-\de_{k'>1}\ze^2c_{2_{k'+1}}-
(\de_{k'>1}+\de_{k'>3})\ze c_{2_{k'-1}}-\de_{k'>3}\ze^2c_{2_{k'-3}}-\de_{k'>5}
\ze c_{2_{k'-5}}\\&-c_{2_{k'+1}}-\de_{k'>1}\ze c_{2_{k'-1}}-\de_{k'>3}
c_{2_{k'-3}}\\&=c_{2_{k'+5}}+\de_{k'>1}\ze c_{2_{k'+3}}+(\de_{k'>1}+\de_{k'>3})
c_{2_{k'+1}}+\de_{k'>3} c_{2_{k'-3}}+\de_{k'>3}\ze c_{2_{k'-1}}\\&+\de_{k'>5}
c_{2_{k'-3}}+\de_{k'>5}\ze c_{2_{k'-5}}+\de_{k'>7}c_{2_{k'-7}}\\&=c_{2_{k'+5}}+
\de_{k'>1}\ze c_{2_{k'+3}}+(\de_{k'>1}+\de_{k'>3})c_{2_{k'+1}}\\&+\de_{k'>3}\ze
c_{2_{k'-1}}+(\de_{k'>3}+\de_{k'>5})c_{2_{k'-3}}+\de_{k'>5}\ze c_{2_{k'-5}}+
\de_{k'>7}c_{2_{k'-7}}.\endalign$$
A similar argument applies for $k\ge 4$. This proves (a).

(b) is obtained by multiplying both sides of (a) by $c_1$ on the left. (c),(d)
are obtained by applying the involution in 3.4 to both sides of (a),(b). We 
prove (e). We have

$c_{2_{2k+2}}c_{1_{k'}}=c_{2_{2k+1}}c_1c_{1_{k'}}=f_1c_{2_{2k+1}}c_{1_{k'}}$
\nl
and the last expression can be computed from (a). This proves (e). Similarly,
(f) follows from (b); (g) is a special case of 6.6. The proposition is proved.

\subhead 17.5\endsubhead
From 7.4,7.6 we see that the function $\De:W@>>>\bn$ has the following values:

$\De(2_{2k})=kL_1+kL_2$,

$\De(2_{2k+1})=-kL_1+(k+1)L_2$,

$\De(1_1)=L_1$,

$\De(1_{2k+1})=(k-1)L_1+kL_2$, if $k\ge 1$,

$\De(1_{2k})=kL_1+kL_2$.

It follows that P1 holds and that $\cd$ consists of the involutions
$2_0=1_0,2_1,1_1,1_3$. Thus, P6 holds.

The formulas in 17.3, 17.4 determine $h_{x,y,z}$ for all $x,y,z$ except when 
$x=1$ or $y=1$, in which case $h_{1,y,z}=\de_{y,z}$, $h_{x,1,z}=\de_{x,z}$.
From these formulas we see that the triples $(x,y,d)$ with $d\in\cd$, 
$\ga_{x,y,d}\ne 0$ are:

$(2_{2k+1},2_{2k+1},2_1)$, $(1_{2k+2},2_{2k+2},1_3)$, $(1_1,1_1,1_1)$,

$(1,1,1)$, $(2_{2k+2},1_{2k+2},2_1)$,$(1_{2k+3},1_{2k+3},1_3)$, 
\nl
where $k\ge 0$. This implies that P2,P3 hold. From the results in 8.8 we see 
that P4,P9,P13 hold. From 14.5 we see that P5 holds. From 14.7 we see that P7 
holds. From 14.8 we see that P8 holds. From 14.10 we see that P10 holds. From 
14.11 we see that P11 holds. From 14.12 we see that P12 holds. From 14.14 we 
see that P14 holds. 

We now verify P15 in our case. With the notation in 14.15, it is enough to show
that, if $a,b\in\{1,2\}$, $w\in W$, $s_aw>w,ws_b>w$, then  
$$(c_ae_w)c'_b-c_a(e_wc'_b)\in\ti{\ch}_{\ge \aa(w)+1}.$$
Here $c_a=c_{s_a},c'_b=c'_{s_b}$. If $a$ or $b$ is $1$, then from 17.2 we have
$(c_ae_w)c'_b-c_a(e_wc'_b)=0$. Hence we may assume that $a=b=2$ and 
$w=1_{2k+1}$. Using 17.2 we have
$$\align&c_2(e_{1_{2k+1}}c'_2)=c_2(e_{1_{2k+2}}+\de_{k>0})\ze'e_{1_{2k}}
+\de_{k>1}e_{1_{2k-2}}\\&=e_{2_{2k+3}}+\ze e_{2_{2k+1}}+\de_{k>0}e_{2_{2k-1}}
+\de_{k>0}\ze'e_{2_{2k+1}}+\de_{k>0}\ze\ze'e_{2_{2k-1}}\\&+\de_{k>1}\ze'
e_{2_{2k-3}}+\de_{k>1}e_{2_{2k-1}}+\de_{k>1}\ze e_{2_{2k-3}}+\de_{k>2}
e_{2_{2k-5}}\\&=e_{2_{2k+3}}+\ze e_{2_{2k+1}}+\de_{k>0}\ze'e_{2_{2k+1}}+
\de_{k>0}e_{2_{2k-1}}+\de_{k>1}e_{2_{2k-1}}\\&+\de_{k>0}\ze\ze'e_{2_{2k-1}}+
\de_{k>1}(\ze+\ze')e_{2_{2k-3}}+\de_{k>2}e_{2_{2k-5}}.\endalign$$
Similarly,
$$\align&(c_2e_{1_{2k+1}})c'_2=e_{2_{2k+3}}+\ze' e_{2_{2k+1}}+\de_{k>0}\ze 
e_{2_{2k+1}}+\de_{k>0})e_{2_{2k-1}}+\de_{k>1}e_{2_{2k-1}}\\&+\de_{k>0}\ze\ze'
e_{2_{2k-1}}+\de_{k>1}(\ze+\ze')e_{2_{2k-3}}+\de_{k>2}e_{2_{2k-5}}.\endalign$$
Hence 

$c_2(e_{1_{2k+1}}c'_2)-(c_2e_{1_{2k+1}})c'_2=(\ze-\ze')(1-\de_{k>0})
e_{2_{2k+1}}$.
\nl
If $k>0$, the right hand side is zero. Thus we may assume that $k=0$. In this
case,

$c_2(e_{1_1}c'_2)-(c_2e_{1_1})c'_2=(\ze-\ze')e_{2_1}$.
\nl
We have $\aa(1_1)=L_1<L_2=\aa(2_1)$. 
This completes the verification of P15 in our
case.

\head 18. The ring $J$\endhead
\subhead 18.1\endsubhead
{\it In this section we assume that $W$ is tame and that P1-P15 in \S14 
are valid.} 

A reference for this section is \cite{\LCC}.

\proclaim{Theorem 18.2} (a) $W$ has only finitely many left cells. 

(b) $W$ has only finitely many right cells. 

(c) $W$ has only finitely many two-sided cells.

(d) $\cd$ is a finite set.
\endproclaim
We prove (a). Since $\aa(w)$ 
is bounded above it is enough to show that, for any
$a\in\bn$, $\aa\i(a)$ is a union of finitely many left cells. By P4, $\aa\i(a)$
is a union of left cells. Let $\ch^1$ be the $\bz$-algebra $\bz\ot_\ca\ch$ 
where $\bz$ is regarded as an $\ca$-algebra via $v\mto 1$. We write $c_w$ 
instead of $1\ot c_w$. For any $a'\ge 0$ let $\ch^1_{\ge a'}$ be the subgroup 
of $\ch^1$ spanned by $\{c_w|\aa(w)\ge a'\}$ (a two-sided ideal of $\ch^1$, by
P4). We have a direct sum decomposition 

(e) $\ch^1_{\ge a}/\ch^1_{\ge a+1}=\opl_\Ga E_\Ga$
\nl
where $\Ga$ runs over the left cells contained in $\aa\i(a)$ and $E_\Ga$ is 
generated as a group by the images of $c_w,w\in\Ga$; these images form a 
$\bz$-basis of $E_\Ga$. Now $\ch^1_{\ge a}/\ch^1_{\ge a+1}$ inherits a left 
$\ch^1$-module structure from $\ch^1$ and (by P9) each $E_\Ga$ is a 
$\ch^1$-submodule. Since $W$ is tame, there exists a finitely generated 
abelian subgroup $W_1$ of finite index of $W$.
Now $\ch^1=\bz[W]$ contains $\bz[W_1]$ as a subring. Since 
$\ch^1_{\ge a}/\ch^1_{\ge a+1}$ is a subquotient of $\ch^1$ (a finitely
generated $\bz[W_1]$-module) and $\bz[W_1]$ is a noetherian ring, it follows 
that $\ch^1_{\ge a}/\ch^1_{\ge a+1}$ is a finitely generated $\bz[W_1]$-module.
Hence in the direct sum decomposition (e) with only non-zero summands, the 
number of summands must be finite. This proves (a).

Since any right cell is of the form $\Ga\i$ where $\Ga$ is a left cell, we see
that (b) follows from (a). Since any two-sided cell is a union of left cells,
we see that (c) follows from (a). From P16 we see that (d) follows from (a). 
The theorem is proved.

\subhead 18.3\endsubhead
Let $J$ be the free abelian group with basis $(t_w)_{w\in W}$. We define
$$t_xt_y=\sum_{z\in W}\ga_{x,y,z\i}t_z.$$
The sum is finite since $\ga_{x,y,z\i}\ne 0\implies h_{x,y,z}\ne 0$ and this
implies that $z$ runs through a finite set (for fixed $x,y$). We show that this
defines an (associative) ring structure on $J$. We must check the identity
$$\sum_z\ga_{x,y,z\i}\ga_{z,u,u'{}\i}=\sum_w\ga_{y,u,w\i}\ga_{x,w,u'{}\i}\tag a
$$
for any $x,y,u,u'\in W$. From P8,P4 we see that both sides of (a) are $0$ 
unless 

(b) $\aa(x)=\aa(y)=\aa(u)=\aa(u')=a$
\nl
for some $a\in\bn$. Hence we may assume that (b) holds. By P8,P4, in the first
sum in (a) we may assume that $\aa(z)=a$ and in the second sum in (a) we may 
assume that $\aa(w)=a$. The equation $(c_xc_y)c_u=c_x(c_yc_u)$ in $\ch$ implies
$$\sum_zh_{x,y,z}h_{z,u,u'}=\sum_wh_{y,u,w}h_{x,w,u'}.\tag c$$
If $h_{x,y,z}h_{z,u,u'}\ne 0$ then $u'\le_{\Cal R}z\le_{\Cal R}x$ hence, by P4,
$\aa(u')\ge\aa(z)\ge\aa(x)$ and $\aa(z)=a$. Hence in the first sum in (c) we 
may assume that $\aa(z)=a$. Similarly in the second sum in (c) we may assume 
that $\aa(w)=a$. Taking the coefficient of $v^{2\aa(z)}$ in both sides of (c) 
we find (a).

The ring $J$ has a unit element $\sum_{d\in\cd}n_dt_d$. Here $n_d=\pm 1$ is as
in 14.1(a), see P5. (The sum is well defined by 18.2(d).) Let us check that 
$t_x\sum_dt_d=t_x$ for $x\in W$. This is equivalent to the identity 
$\sum_dn_d\ga_{x,d,z\i}=\de_{z,x}$. By P7 this is equivalent to 
$\sum_dn_d\ga_{z\i,x,d}=\de_{z,x}$. This follows from P2,P3,P5. The equality 
$(\sum_dt_d)t_x=t_x$ is checked in a similar way.

For any subset $X$ of $W$, let $J^X$ be the subgroup of $J$ generated by 
$\{t_x|x\in X\}$. If $\boc$ is a two-sided cell of $W,L$ then, by P8, 
$J^{\boc}$ is a subring of $J$ and $J=\opl_\boc J^\boc$ is a direct sum 
decomposition of $J$ as a ring. The unit element of $J^\boc$ is 
$\sum_{d\in\cd\cap\boc}t_d$. Similarlym if $\Ga$ is a left cell of $W,L$ then
$J^{\Ga\cap\Ga\i}$ is a subring of $J$ with unit element $t_d$ where 
$d\in\cd\cap\Ga$. 

\proclaim{Proposition 18.4} Assume that we are in the setup of 15.1. Let 
$x,y\in W$.

(a) The condition $x\sim_\cl y$ is equivalent to the condition that
$t_xt_{y\i}\ne 0$ and to the condition that, for some $u$, $t_y$ appears with 
$\ne 0$ coefficient in $t_ut_x$.

(b) The condition $x\sim_{\Cal R}y$ is equivalent to the condition that
$t_{x\i}t_y\ne 0$ and to the condition that, for some $u$, $t_y$ appears with 
$\ne 0$ coefficient in $t_xt_u$.

(c) The condition $x\sim_{\cl\Cal R}y$ is equivalent to the condition that
$t_xt_ut_y\ne 0$ for some $u$ and to the condition that, for some $u,u'$, $t_y$
appears with $\ne 0$ coefficient in $t_{u'}t_xt_u$.
\endproclaim
Let $J^+=\sum_z\bn t_z$. By 15.1(a) we have $J^+J^+\sub J^+$.

We prove (a). The second condition is equivalent to $\ga_{x,y\i,u}\ne 0$ for 
some $u$; the third condition is equivalent to $\ga_{x,u,y\i}\ne 0$ for some
$u$. These conditions are equivalent by P7. 

Assume that $\ga_{x,y\i,u}\ne 0$ for some $u$. Using P8 we deduce that 
$x\sim_\cl y$. 

Assume now that $x\sim_\cl y$. Let $d\in\cd$ be such that $x\sim_\cl d$. Then 
we have also $y\sim_\cl d$. By P13 we have $\ga_{x\i,x,d}\ne 0$,
$\ga_{y\i,y,d}\ne 0$. Hence $\ga_{x\i,x,d}=1,\ga_{y\i,y,d}=1$. Hence 
$t_{x\i}t_x\in t_d+J^+$, $t_{y\i}t_y\in t_d+J^+$. Since $t_dt_d=t_d$, it 
follows that $t_{x\i}t_xt_{y\i}t_y\in t_dt_d+J^+=t_d+J^+$. In particular,
$t_xt_{y\i}\ne 0$. This proves (a).

The proof of (b) is entirely similar.

We prove (c). Using the associativity of $J$ we see that the third condition on
$x,y$ is a transitive relation on $W$. Hence to prove that the first condition
implies the third condition we may assume that either $x\sim_\cl y$ or
$x\sim_{\Cal R}y$, in which case this follows from (a) or (b). The fact that 
the third condition implies the first condition also follows from (a),(b). Thus
the first and third condition are equivalent. 

Assume that $t_xt_ut_y\ne 0$ for some $u$. By (a),(b) we then have 
$x\sim_\cl u\i, u\i\sim_{\Cal R}y$. Hence $x\sim_{\cl\Cal R}y$. 

Conversely, assume that $x\sim_{\cl\Cal R}y$. Using P14 we deduce that
$x\sim_{\cl\Cal R}y\i$. By the earlier part of the proof, $t_{y\i}$ appears 
with $\ne 0$ coefficient in $t_{u'}t_xt_u$ for some $u,u'$. We have 
$t_{u'}t_xt_u\in at_{y\i}+J^+$ where $a>0$. Hence
$t_{u'}t_xt_ut_y\in at_{y\i}t_y+J^+$. Since $t_{y\i}t_y$ has a coefficient $1$
and the other coefficients are $\ge 0$, it folows that $t_{u'}t_xt_ut_y\ne 0$.
Thus, $t_xt_ut_y\ne 0$. We see that the first and second conditions are 
equivalent. The proposition is proved.

\subhead 18.5\endsubhead
Assume now that we are in the setup of 7.1 with $m=\infty$ and $L_2>L_1$. From
the formulas in 17.3,17.4 we can determine the multiplication table of $J$. We
find

$t_{2_{2k+1}}t_{2_{2k'+1}}=\sum_{u\in[0,\ti k]}t_{2_{2k+2k'+1-4u}}$,

$t_{1_{2k+3}}t_{1_{2k'+3}}=\sum_{u\in[0,\ti k]}t_{1_{2k+2k'+3-4u}}$,

$t_{2_{2k+1}}t_{2_{2k'+2}}=\sum_{u\in[0,\ti k]}t_{2_{2k+2k'+2-4u}}$,

$t_{1_{2k+3}}t_{1_{2k'+2}}=\sum_{u\in[0,\ti k]}t_{1_{2k+2k'+2-4u}}$,

$t_{2_{2k+2}}t_{1_{2k'+3}}=\sum_{u\in[0,\ti k]}t_{2_{2k+2k'+2-4u}}$,

$t_{2_{2k+2}}t_{1_{2k'+2}}=\sum_{u\in[0,\ti k]}t_{2_{2k+2k'+1-4u}}$,

$t_{1_{2k+2}}t_{2_{2k'+1}}=\sum_{u\in[0,\ti k]}t_{1_{2k+2k'+2-4u}}$,

$t_{1_{2k+2}}t_{2_{2k'+2}}=\sum_{u\in[0,\ti k]}t_{1_{2k+2k'+3-4u}}$,

$t_{1_1}t_{1_1}=t_{1_1}$,

$t_1t_1=t_1$;
\nl
here $k,k'\ge 0$ and $\ti k=\min(k,k')$. All other products are $0$. 

Let $R$ be the free abelian group with basis $(b_k)_{k\in\bn}$. We regard $R$ 
as a commutative ring with multiplication 

$b_kb_{k'}=\sum_{u\in[0,\min(k,k')]}b_{k+k'-2u}$.
\nl
Let $J_0=\sum_{w\in W-\{1,1_1\}}\bz t_w$. The formulas above show that 
$J=J_0\opl\bz t_1\opl\bz t_{1_1}$ (direct sum of rings) and that the ring $J_0$
is isomorphic to the ring of $2\tim 2$ matrices with entries in $R$, via the 
isomorphism defined by:
$$t_{2_{2k+1}}\mto\left(\sm b_k&0\\0&0\esm\right),\quad
t_{1_{2k+3}}\mto\left(\sm 0&0\\0&b_k\esm\right),\quad    
t_{2_{2k+2}}\mto\left(\sm 0&b_k\\0&0\esm\right),\quad  
t_{1_{2k+2}}\mto\left(\sm 0&0\\b_k&0\esm\right).$$
Note that $R$ is canonically isomorphic to the representation ring of 
$SL_2(\bc)$ with its canonical basis consisting of irreducible representations.

\subhead 18.6\endsubhead
Assume that we are in the setup of 7.1 with $m=\infty$ and $L_2=L_1$. By 
methods similar (but simpler) to those of \S17 and 18.5, we find

$t_{2_{2k+1}}t_{2_{2k'+1}}=\sum_{u\in[0,2\min(k,k')]}t_{2_{2k+2k'+1-2u}}$.
\nl
Let $J^1$ be the subring of $J$ generated by $t_{2_{2k+1}},k\in\bn$. While, in
18.5, the analogue of $J^1$ was isomorphic to $R$ as a ring with basis, in the
present case, $J^1$ is canonically isomorphic to $R'$, the subgroup of $R$ 
generated by $b_k$ with $k$ even. (Note that $R'$ is a subring of $R$,
naturally isomorphic to the representation ring of $PGL_2(\bc)$.) 

\subhead 18.7\endsubhead
In the setup of 7.1 with $m=4$ and $L_2=2,L_1=1$ (a special case of the
situation in \S15), we have

$J=\bz t_1\opl \bz t_{1_1}\opl J_0\opl \bz t_{2_3}\opl \bz t_{2_4}$
\nl
(direct sum of rings) where $J_0$ is the subgroup of $J$ generated by
$t_{2_1},t_{2_2},t_{1_2},t_{1_3}$. The ring $J_0$ is isomorphic to the ring of
$2\tim 2$ matrices with entries in $\bz$, via the isomorphism defined by:
$$t_{2_1}\mto\left(\sm 1&0\\0&0\esm\right),\quad
t_{1_3}\mto\left(\sm 0&0\\0&1\esm\right),\quad    
t_{2_2}\mto\left(\sm 0&1\\0&0\esm\right),\quad  
t_{1_2}\mto\left(\sm 0&0\\1&0\esm\right).$$
Moreover, $t_1,t_{1_1},t_{2_4}$ are idempotent. On the other hand,

$t_{2_3}t_{2_3}=-t_{2_3}$.
\nl
Notice the minus sign! (It is a special case of the computation in ....)

\subhead 18.8 \endsubhead
For any $z\in W$ we set $\hat n_z=n_d$ where $d$ is the unique element of $\cd$
such that $d\sim_{\cl}z\i$ and $n_d=\pm 1$ is as in 14.1(a), see P5. Note that
$z\mto\hat n_z$ is constant on right cells.

\proclaim{Theorem 18.9} Let $J_\ca=\ca\ot J$ and let $\phi:\ch@>>>J_\ca$ be the
$\ca$-linear map given by 

$\phi(c_x^\dag)=\sum_{z\in W,d\in\cd;\aa(d)=\aa(z)}h_{x,d,z}\hat n_zt_z$ 
\nl
for all $x\in W$. Then $\phi$ is a homomorphism of $\ca$-algebras with $1$.
\endproclaim
Consider the equality

(a) $\sum_wh_{x_1,x_2,w}h'_{w,x_3,y}=\sum_wh_{x_1,w,y}h'_{x_2,x_3,w}$
\nl
(see P15) with $\aa(x_2)=\aa(y)=a$. In the left hand side we may assume that 
$y\le_{\Cal R}w\le_\cl x_2$ hence (by P4) $\aa(y)\ge\aa(w)\ge\aa(x_2)$, hence 
$\aa(w)=a$. Similarly in the right hand side we may assume that $\aa(w)=a$. 
Picking the coefficient of $v'{}^a$ in both sides of (a) gives

(b) $\sum_wh_{x_1,x_2,w}\ga_{w,x_3,y\i}=\sum_wh_{x_1,w,y}\ga_{x_2,x_3,w\i}$.
\nl
Let $x,x'\in W$. The desired identity 
$\phi(c_x^\dag c_{x'}^\dag)=\phi(c_x^\dag)\phi(c_{x'}^\dag)$ 
is equivalent to 
$$\sum\Sb w\in W,d\in\cd\\\aa(d)=a'\endSb h_{x,x',w}h_{w,d,u}\hat n_u=
\sum\Sb z,z'\in W,d,d'\in\cd\\\aa(d)=\aa(z)\\\aa(d')=\aa(z')\endSb h_{x,d,z}
h_{x',d',z'}\ga_{z,z',u\i}\hat n_z\hat n_{z'}$$
for any $u\in W$ such that $\aa(u)=a'$. In the right hand we may assume that 

$\aa(d)=\aa(z)=\aa(d')=\aa(z')=a'$ and $\hat n_z=\hat n_u$
\nl
(by P8,P4). Hence the right hand side can be rewritten (using (b)):
$$\align&\sum\Sb
z'\in W,d,d'\in\cd\\\aa(d)=\aa(d')=\aa(z')=a'\endSb h_{x',d',z'}
\sum_{z; \aa(z)=a'}h_{x,d,z}\ga_{z,z',u\i}\hat n_u\hat n_{z'}\\&
=\sum\Sb z'\in W,d,d'\in\cd\\ \aa(d)=\aa(d')=\aa(z')=a'\endSb h_{x',d',z'}
\sum_{w; \aa(w)=a'}h_{x,w,u}\ga_{d,z',w\i}\hat n_u\hat n_{z'}.\endalign$$
By P2,P3,P5, this equals
$$\sum_{z'\in W,d'\in\cd;\aa(d')=\aa(z')=a'}h_{x',d',z'}h_{x,z',u}\hat n_u$$
which by the identity $(c_xc_{x'})c_{d'}=c_x(c_{x'}c_{d'})$ equals
$$\sum_{w\in W,d'\in\cd;\aa(d')=a'}h_{x,x',w}h_{w,d',u}\hat n_u.$$
Thus $\phi$ is compatible with multiplication.

Next we show that $\phi$ is compatible with the unit elements of the two
algebras. An equivalent statement is that for any $z\in W$ such that 
$\aa(z)=a$, the sum $\sum_{d\in\cd;\aa(d)=a}h_{1,d,z}\hat n_z$ equals $n_z$ if
$z\in\cd$ and is $0$ if $z\notin\cd$. This is clear since 
$h_{1,d,z}=\de_{z,d}$.

\subhead 18.10\endsubhead
If we identify the $\ca$-modules $\ch$ and $J_\ca$ via 
$c_w^\dag\mto\hat n_wt_w$, the obvious left $J_\ca$-module structure on $J_\ca$
becomes the left $J_\ca$-module structure on $\ch$ given by

$t_x*c_w^\dag=\sum_{z\in W}\ga_{x,w,z\i}\hat n_w\hat n_zc_z^\dag$
\nl
Let $\ch_a=\opl_{w;\aa(w)=a}\ca c_w^\dag,
\ch_{\ge a}=\opl_{w;\aa(w)\ge a}\ca c_w^\dag$. 
Note that $t_x*c_w^\dag\in\ch_{\aa(w)}$ for all $x,w$. For any $h\in\ch,w\in W$
we have

(a) $hc_w^\dag=\phi(h)*c_w^\dag\mod \ch_{\ge \aa(w)+1}$.
\nl
Indeed, we may assume that $h=c_x^\dag$. Using 18.9(b), we have
$$\align&\phi(c_x^\dag)*c_w^\dag=\sum\Sb d\in\cd,z\\ \aa(d)=\aa(z)\endSb 
h_{x,d,z}\hat n_zt_z*c_w^\dag\\&=\sum\Sb d\in\cd,z,u\\\aa(d)=\aa(z)
\endSb h_{x,d,z}\ga_{z,w,u\i}\hat n_z
\hat n_w\hat n_uc_u^\dag=\sum\Sb d\in\cd,z,u\\\aa(d)=\aa(w)=\aa(u)\endSb 
h_{x,d,z}
\ga_{z,w,u\i}\hat n_wc_u^\dag\\&=\sum\Sb d\in\cd,t,u\\ \aa(d)=\aa(w)=\aa(u)
\endSb h_{x,t,u}
\ga_{d,w,t\i}\hat n_wc_u^\dag=\sum\Sb d\in\cd,u\\ \aa(d)=\aa(w)=\aa(u)\endSb 
h_{x,w,u}\ga_{d,w,w\i}\hat n_wc_u^\dag\\&=\sum\Sb u\\ \aa(w)=\aa(u)
\endSb h_{x,w,u}c_u^\dag=c_x^\dag c_w^\dag
\mod\ch_{\ge \aa(w)+1},\endalign$$
as required.

\subhead 18.11\endsubhead
Let $\ca@>>>R$ be a ring homomorphism of $\ca$ into a commutative ring $R$ with
$1$. Let $\ch_R=R\ot_\ca\ch,J_R=R\ot_\ca(J_\ca)=R\ot J$,
$\ch_{R,\ge a}=R\ot_\ca\ch_{\ge a}$. Then $\phi$ extends to a homomorphism of 
$R$-algebras $\phi_R:\ch_R@>>>J_R$. The $J_\ca$-module in 18.10 extends to a 
$J_A$-module structure on $\ch_R$ denoted again by $*$. From 18.10(a) we deduce

(a) $hc_w^\dag=\phi_R(h)*c_w^\dag\mod\ch_{R,\ge \aa(w)+1}$ for any 
$h\in\ch_R,w\in W$.

\proclaim{Proposition 18.12} (a) If $N$ is a bound for $W,L$, then 
$(\ker\phi_R)^{N+1}=0$.

(b) If $R=R_0[v,v\i]$ where $R_0$ is a commutative ring with $1$, $v$ is an
indeterminate and $\ca@>>>R$ is the obvious ring homomorphism, then 
$\ker\phi_R=0$.
\endproclaim
We prove (a). If $h\in\ker\phi_R$ then by 18.11(a), we have
$h\ch_{R,\ge a}\sub\ch_{R,\ge a+1}$ for any $a\ge 0$. Applying this repeatedly,
we see that, if $h_1,h_2,\dots,h_{N+1}\in\ch$, we have
$h_1h_2\dots h_{N+1}\in\ch_{R,\ge N+1}=0$. This proves (a).

We prove (b). Let $h=\sum_xp_xc_x^\dag\in\ker\phi_R$ where $p_x\in R$. Assume 
that $h\ne 0$. Then $p_x\ne 0$ for some $x$. We can find $a\ge 0$ such that 
$p_x\ne 0\implies\aa(x)\ge a$ and $X=\{x\in W|p_x\ne 0,\aa(x)=a\}$ is 
non-empty.
We can find $b\in\bz$ such that $p_x\in v^b\bz[v\i]$ for all $x\in X$ and such
that $X'=\{x\in X|\pi_b(p_x)\ne 0\}$ is non-empty. Let $x_0\in X'$. We can find
$d\in\cd$ such that $\ga_{x_0,d,x_0\i}=\ga_{x_0\i,x_0,d}\ne 0$. We have 
$hc_d^\dag=\sum_xp_xc_x^\dag c_d^\dag$. If $\aa(x)>a$, then 
$c_x^\dag c_d^\dag\in\ch_{R,\ge a+1}$. Hence
$hc_d^\dag=\sum_{x\in X}p_xc_x^\dag c_d^\dag\mod\ch_{R,\ge a+1}$. Since 
$\phi_R(h)=0$, from
18.11(a) we have $hc_d^\dag\in\mod\ch_{R,\ge a+1}$. It follows that
$\sum_{x\in X}p_xc_x^\dag c_d^\dag\in\ch_{R,\ge a+1}$. In particular the 
coefficient of $c_{x_0}^\dag$ in $\sum_{x\in X}p_xc_x^\dag c_d^\dag$ is $0$. 
In other words, $\sum_{x\in X}p_xh_{x,d,x_0}=0$. The coefficient of $v^{a+b}$ in the last sum is

$\sum_{x\in X}\pi_b(p_x)\ga_{x,d,x_0\i}=\pi_b(p_{x_0})\ga_{x_0,d,x_0\i}$
\nl
and this is on the one hand $0$ and on the other hand is non-zero since
$\pi_b(p_{x_0})\ne 0$ and $\ga_{x_0,d,x_0\i}\ne 0$, by the choice of $x_0,d$.
This contradiction completes the proof.

\head 19. Algebras with trace form\endhead  
\subhead 19.1\endsubhead
Let $R$ be a field and let $A$ be an associative $R$-algebra with $1$ of finite
dimension over $R$. We assume that $A$ is semisimple and split over $R$ and
that we are given a {\it trace form} on $A$ that is, an $R$-linear map 
$\tau:A@>>>R$ such that $(a,a')=\tau(aa')=\tau(a'a)$ is a non-degenerate 
(symmetric) $R$-bilinear form (,):$A\tim A@>>>R$. Note that 
$(aa',a'')=(a,a'a'')$ for all $a,a',a''$ in $A$. Let $\Mod A$ be the category 
whose objects are left $A$-modules of finite dimension over $R$. We write 
$E\in\Irr A$ for "$E$ is a simple object of $\Mod A$".

Let $(a_i)_{i\in I}$ be an $R$-basis of $A$ and let $(a'_i)_{i\in I}$ be the 
$R$-basis defined by $(a_i,a'_j)=\de_{ij}$. Then 

(a) $\sum_ia_i\ot a'_i\in A\ot A$ is independent of the choice of $(a_i)$.

\proclaim{Proposition 19.2} (a) We have $\sum_i\tau(a_i)a'_i=1$. 

(b) If $E\in\Irr A$, then $\sum_i\tr(a_i,E)a'_i$ is in the centre of $A$. It 
acts on $E$ as a scalar $f_E\in R$ times the identity and on $E'\in\Irr A$, not
isomorphic to $E$, as zero. Moreover, $f_E$ does not depend on the choice of 
$(a_i)$.

(c) One can attach uniquely to each $E\in\Irr A$ a scalar $g_E\in R$ (depending
only on the isomorphism class of $E$), so that

$\sum_Eg_E\tr(a,E)=\tau(a)$ for all $a\in A$,
\nl
where the sum is taken over all $E\in\Irr A$ up to isomorphism. 

(d) For any $E\in\Irr A$ we have $f_Eg_E=1$. In particular, 
$f_E\ne 0, g_E\ne 0$.

(e) If $E,E'\in\Irr A$, then $\sum_i\tr(a_i,E)\tr(a'_i,E')$ is $f_E\dim E$ if 
$E,E'$ are isomorphic and is $0$, otherwise.
\endproclaim
Let $A=\opl_{n=1}^tA_n$ be the decomposition of $A$ as a sum of simple 
algebras. Let $\tau_n:A_n@>>>R$ be the restriction of $\tau$. Then $\tau_n$ is
a trace form for $A_n$, whose associated form is the restriction of $(,)$ and 
$(A_n,A_{n'})=0$ for $n\ne n'$. Hence we can choose $(a_i)$ so that each $a_i$
is contained in some $A_n$ and then $a'_i$ will be contained in the same $A_n$
as $a'_i$.

We prove (a). From 19.1(a) we see that $\sum_i\tau(a_i)a'_i$ is independent of
the choice of $(a_i)$. Hence we may choose $(a_i)$ as in the first paragraph of
the proof. We are thus reduced to the case where $A$ is simple. In that case   
the assertion is easily verified.

We prove (b).  From 19.1(a) we see that $\sum_i\tr(a_i,E)a'_i$ is independent
of the choice of $(a_i)$. Hence we may choose $(a_i)$ as in the first paragraph
of the proof. We are thus reduced to the case where $A$ is simple. In that case
the assertion is easily verified.

We prove (c). It is enough to note that $a\mto\tr(a,E)$ form a basis of the 
space of $R$-linear functions $A@>>>R$ which vanish on all $aa'-a'a$ and 
$\tau$ is such a function. 

We prove (d). We consider the equation in (c) for $a=a_i$ and we multiply both
sides by $a'_i$ and sum over $i$. Using (a), we obtain

$\sum_i\sum_Eg_E\tr(a_i,E)a'_i=\sum_i\tau(a_i)a'_i=1$.
\nl
Hence $\sum_Eg_E\sum_i\tr(a_i,E)a'_i=1$. By (b), the left hand side acts on a
$E'\in\Irr A$ as a scalar $g_{E'}f_{E'}$ times the identity. This proves (d).

(e) follows immediately from (b). The proposition is proved.

\subhead 19.3\endsubhead
Now let $A'$ be a semisimple subalgebra of $A$ such that $\tau'$, the 
restriction of $\tau$ to $A'$ is a trace form of $A'$. (We do not assume that
the unit element $1_{A'}$ of $A'$ coincides to the unit element $1$ of $A$.) If
$E\in\Mod A$ then $1_{A'}E$ is naturally an object of $\Mod A'$. Hence if 
$E'\in\Irr A'$, then the multiplicity $[E':1_{A'}E]$ of $E'$ in $1_{A'}E'$ is 
well defined.

Note that, if $a'\in A'$, then $\tr(a',1_{A'}E)=\tr(a',E)$.

\proclaim{Lemma 19.4} Let $E'\in\Irr A'$. We have

$g_{E'}=\sum_E[E':1_{A'}E]g_E$,
\nl
sum over all $E\in\Irr A$ (up to isomorphism).
\endproclaim
By the definition of $g_{E'}$, it is enough to show that

(a) $\sum_{E'}\sum_E[E':1'E]g_E\tr(a',E')=\tau(a')$
\nl
for any $a'\in A'$. Here $E'$ (resp. $E$) runs over the isomorphism classes of
simple objects of $\Mod A'$ (resp. $\Mod A$). The left hand of (a) is

$\sum_Eg_E\sum_{E'}[E':1'E]\tr(a',E')=\sum_Eg_E\tr(a',1'E)=\sum_Eg_E\tr(a',E)$
\nl
which, by the definition of $g_E$ is equal to $\tau(a')$. This completes the
proof.

\head 20. The function $\aa_E$\endhead
\subhead 20.1\endsubhead
In this section we assume that the assumptions of 18.1 hold and that $W$ is 
finite. 

The results of \S19 will be applied in the following cases.

(a) $A=\ch_\bc, R=\bc$. Here $\ca@>>>\bc$ takes $v$ to $1$. We identify
$\ch_\bc$ with the group algebra $\bc[W]$ by $w\mto T_w$ for all $w$. It is 
well known that $\bc[W]$ is a semisimple split algebra. We take $\tau$ so that
$\tau(x)=\de_{x,1}$ for $x\in W$. Then the bases $(x)$ and $(x\i)$ are dual
with respect to $(,)$. 

We will say $W$-module instead of $\bc[W]$-module.
We will write $\Mod W,\Irr W$ instead of $\Mod\bc[W],\Irr\bc[W]$.

(b) $A=J_\bc$, $R=\bc$. Since $\bc[W]$ is semisimple, we see from 18.12(a) that
the kernel of $\phi_\bc:\bc[W]@>>>J_\bc$ is $0$ so that $\phi_\bc$ is 
injective. Since $\dim\bc[W]=\dim J_\bc=\sha W$ it follows that $\phi_\bc$ is
an isomorphism. In particular $J_\bc$ is a semisimple split algebra. We take
$\tau:J_\bc@>>>\bc$ so that $\tau(t_z)$ is $n_z$ if $z\in\cd$ and $0$, 
otherwise. Then $(t_x,t_y)=\de_{xy,1}$. The bases $(t_x)$ and $(t_{x\i})$ are 
dual with respect to $(,)$.

(c) $A=\ch_{\bc(v)}, R=\bc(v)$. Here $\ca@>>>\bc$ takes $v$ to $v$. The 
homomorphism $\phi_{\bc(v)}:\ch_{\bc(v)}@>>>J_{\bc(v)}$ is injective. This 
follows from 18.12(b), using the fact that injectivity is preserved by 
tensoring with a field of fractions. Since $\ch_{\bc(v)},J_{\bc(v)}$ have the 
same dimension, it follows that $\phi_{\bc(v)}$ is an isomorphism. Since
$J_{\bc(v)}=\bc(v)\ot J_\bc$, and $J_\bc$ is semisimple, split, it follows that
$J_{\bc(v)}$ is semisimple, split, hence $\ch_{\bc(v)}$ is semisimple, split.
We take $\tau:\ch_{\bc(v)}$ so that $\tau(T_w)=\de_{w,1}$. The bases $(T_x)$ 
and $(T_{x\i})$ are dual with respect to $(,)$.

{\it Remark.} The argument above shows also that,

(d) if $R=R_0(v)$, with $R_0$ an arbitrary field and $\ca@>>>R$ carries $v$ to
$v$, then $\phi_R:\ch_R@>>>J_R$ is an isomorphism;

(e) if $R$ in 18.11 is a field of characteristic $0$ then $\phi_R:\ch_R@>>>J_R$
is an isomorphism if and only if $\ch_R$ is a semisimple $R$-algebra.

\subhead 20.2\endsubhead
For any $E\in\Mod W$ we denote by $E_\spa$ the corresponding $J_\bc$-module.
Thus, $E_\spa$ coincides with $E$ as a $\bc$-vector space and the action of 
$j\in J_\bc$ on $E_\spa$ is the same as the action of $\phi_\bc\i(j)$ on $E$.
The $J_\bc$-module structure on $E_\spa$ extends in a natural way to a 
$J_{\bc(v)}$-module structure on $E_v=\bc(v)\ot_\bc E_\spa$. We will also 
regard $E_v$ as an $\ch_{\bc(v)}$-module via the algebra isomorphism 
$\phi_{\bc(v)}:\ch_{\bc(v)}@>\sim>>J_{\bc(v)}$. If $E$ is simple, then 
$E_\spa$ and $E_v$ are simple.

If $E\in\Irr W$. Then $E_\spa$ is a simple $J^\boc_\bc$-module for a unique 
two-sided cell $\boc$ of $W$. Then for any $x\in\boc$, we write 
$E\sim_{\cl\Cal R}x$. If $E,E'\in\Irr W$, we write $E\sim_{\cl\Cal R}E'$ if for
some $x\in W$ we have $E\sim_{\cl\Cal R}x$,$E'\sim_{\cl\Cal R}x$.

\subhead 20.3\endsubhead
There is the following direct relationship between $E$ and $E_v$ (without
going through $J$):

(a) $\tr(x,E)=\tr(T_x,E_v)|_{v=1}$ for all $x\in W$.
\nl
Indeed, it is enough to show that $\tr(c_x^\dag,E)=\tr(c_x^\dag,E_v)|_{v=1}$. 
Both sides are equal to 
$\sum_{z\in W,d\in\cd}\ga_{x,d,z\i}\hat n_z\tr(t_z,E_\spa)$.

\subhead 20.4\endsubhead
Assume that $E\in\Irr W$. We have 

(a) $(f_{E_v})_{v=1}\dim(E)=|W|$.
\nl
Indeed, setting $v=1$ in 
$\sum_{x\in W}\tr(T_x,E_v)\tr(T_{x\i},E_v)=f_{E_v}\dim(E)$ gives

$\sum_{x\in W}\tr(x,E)\tr(x\i,E)=(f_{E_v})_{v=1}\dim(E)$.
\nl
The left hand side equals $|W|$; (a) follows.

\subhead 20.5\endsubhead
Let $I\sub S$, let $E'\in\Irr W_I$ and let $E\in\Irr W$. We have 

(a) $[E'_v:E_v]=[E':E]$.
\nl
The right hand side is $|W_I|\i\sum_{x\in W_I}\tr(x,E')\tr(x\i,E)$. The left 
hand side is 
$$f_{E'_v}\i\dim(E')\i\sum_{x\in W_I}\tr(T_x,E'_v)\tr(T_{x\i},E_v).$$ 
Since this is a constant, it is equal to its value for $v=1$. Hence it is equal
to

$(f_{E'_v}\i)_{v=1}\dim(E')\i\sum_{x\in W_I}\tr(x,E'_v)\tr(x\i,E_v)$.
\nl
Thus it is enough to show that $(f_{E'_v})_{v=1}\dim(E')=|W_I|$. But this is a
special case of 20.4(a). 

\proclaim{Proposition 20.6} Let $E\in\Irr W$.

(a) There exists a unique integer $\aa_E\ge 0$ such that 
$\tr(T_x,E_v)\in v^{-\aa_E}\bc[v]$ for all $x\in W$ and 
$\tr(T_x,E_v)\notin v^{-\aa_E+1}\bc[v]$ for some $x\in W$.

(b) For $x\in W$ we have 
$\tr(T_x,E_v)=\sgn(x)v^{-\aa_E}\tr(t_x,E_\spa)\mod v^{-\aa_E+1}\bc[v]$.

(c) Let $\boc$ be the two-sided cell such that $E_\spa\in\Irr J^\boc_\bc$. Then
$\aa_E=\aa(x)$ for any $x\in\boc$.
\endproclaim
Let $a=\aa(x)$ for any $x\in\boc$. By definition,

$\tr(c_x^\dag,E_v)=\sum_{z\in W,d\in\cd;\aa(d)=\aa(z)}h_{x,d,z}\hat n_z
\tr(t_z,E_\spa)$.
\nl
In the last sum we have $\tr(t_z,E_\spa)=0$ unless $z\in\boc$ in which case 
$\aa(z)=a$. For such $z$ we have 
$h_{x,d,z}=\bar h_{x,d,z}=\ga_{x,d,z\i}v^{-a}\mod v^{-a+1}\bz[v]$, hence we 
have

$\tr(c_x^\dag,E_v)
=\sum_{z\in W,d\in\cd}\ga_{x,d,z\i}\hat n_zv^{-a}\tr(t_z,E_\spa)
\mod v^{-a+1}\bc[v]$.
\nl
For each $z$ in the last sum we have 
$\sum_{d\in\cd}\ga_{x,d,z\i}\hat n_z=\de_{x,z}n_d\hat n_z=\de_{x,z}$. From this
we deduce that

(d) $\tr(c_x^\dag,E_v)=v^{-a}\tr(t_x,E_\spa)\mod v^{-a+1}\bc[v]$.
\nl
We have $T_x=\sum_{y;y\le x}q'_{y,x}c_y$. Hence
$\sgn(x)\bar T_x=T_x^\dag=\sum_{y;y\le x}q'_{y,x}c_y^\dag$. Applying $\bar{}$
we obtain $\sgn(x)T_x=\sum_{y;y\le x}\bar q'_{y,x}c_y^\dag$. Hence

$\tr(T_x,E_v)=\sgn(x)\sum_{y;y\le x}\bar q'_{y,x}\tr(c_y^\dag,E_v)$.
\nl
Using (d) together with $\bar q'_{x,x}=1$, $\bar q'_{y,x}\in v\bz[v]$ (see 
10.1), we deduce

$\tr(T_x,E_v)=\sgn(x)v^{-a}\tr(t_x,E_\spa)\mod v^{-a+1}\bc[v]$.
\nl
Since $E_\spa\in\Irr J_\bc$, we have $\tr(t_x,E_\spa)\ne 0$ 
for some $x\in W$. The proposition follows.

\proclaim{Corollary 20.7} $f_{E_v}=f_{E_\spa}v^{-2\aa_E}
+\text{strictly higher powers of } v$. 
\endproclaim  
Using 19.2(e) for $\ch_{\bc(v)}$ and $J_\bc$, we obtain
$$\align &f_{E_v}\dim E=\sum_x\tr(T_x,E_v)\tr(T_{x\i},E_v) \\&\in
v^{-2\aa_E}\sum_x\tr(t_x,E_\spa)\tr(t_{x\i},E_\spa)+v^{-2\aa_E+1}\bc[v]\\&
=v^{-2\aa_E}f_{E_\spa}\dim E+v^{-2\aa_E+1}\bc[v].\endalign$$
The corollary follows.

\medpagebreak

Let $\bar{}:\bc[v,v\i]@>>>\bc[v,v\i]$ be the $\bc$-algebra homomorphism given 
by $v^n\mto v^{-n}$ for all $n$.

\proclaim{Corollary 20.8} For any $x\in W$ we have 
$\Tr(T_{x\i}\i,E_v)=\ov{\tr(T_x,E_v)}$.
\endproclaim
We have $T_x=\sgn(x)\sum_{y;y\le x}\bar q'_{y,x}c_y^\dag$,
$T_{x\i}\i=\sgn(x)\sum_{y;y\le x}q'_{y,x}c_y^\dag$, hence

$\tr(T_x,E_v)=\sgn(x)\sum_{y;y\le x}\bar q'_{y,x}\tr(c_y^\dag,E_v)$,

$\tr(T_{x\i}\i,E_v)=\sgn(x)\sum_{y;y\le x}q'_{y,x}\tr(c_y^\dag,E_v)$.
\nl  
Thus, it suffices to show that $\Tr(c_y^\dag,E_v)=\ov{\tr(c_y^\dag,E_v)}$ for 
any $y\in W$. As in the proof of 20.6 we have

$\tr(c_y^\dag,E_v)=\sum_{z\in W,d\in\cd;\aa(d)=\aa(z)}h_{y,d,z}\hat n_z
\tr(t_z,E_\spa)$.
\nl
Hence it suffices to show that $\ov{h_{y,d,z}}=h_{y,d,z}$ for all $d,z$ in the
last sum. But this is clearly true for any $y,d,z$ in $W$.

\medpagebreak

For any $E\in\Mod W$ we write $E^\dag$ instead of $E\ot\sgn$. We write 
$E_v^\dag$ instead of $(E^\dag)_v$

\proclaim{Lemma 20.9} Let $E\in\Irr W$. For any $x\in W$ we have

$\tr(T_x,(E^\dag)_v)=(-1)^{l(x)}\ov{\tr(T_x,E_v)}$
\endproclaim
${}^\dag:\ch@>>>\ch$ (see 3.5) extends uniquely to a $\bc(v)$-algebra 
involution ${}^\dag:\ch_{\bc(v)}@>>>\ch_{\bc(v)}$. Let 
$(E_v)^\dag$ be the $\ch_{\bc(v)}$-module with underlying vector space $E_v$ 
such that the action of $h$ on $E_v^\dag$ is the same as the action of $h^\dag$
on $E_v$. Clearly, $(E_v)^\dag\in\Irr\ch_{\bc(v)}$. For $x\in W$ we have 

$\tr(T_x,(E_v)^\dag)=(-1)^{l(x)}\tr(T_{x\i}\i,E_v)=
(-1)^{l(x)}\ov{\tr(T_x,E_v)}$.
\nl
(The last equation follows from 20.8.) Setting $v=1$ we obtain

$\tr(T_x,(E_v)^\dag)|_{v=1}=(-1)^{l(x)}\tr(x,E)=\tr(x,E^\dag)$.
\nl
Using 20.3, we deduce that $(E_v)^\dag\cong E^\dag_v$ in $\Mod\ch_{\bc(v)}$.
The lemma follows.

\proclaim{Proposition 20.10} For any $x\in W$ we have

$\tr(T_x,E_v)=v^{\aa_{E^\dag}}\tr(t_x,E^\dag_\spa)+\text{strictly lower
powers of } v$.
\endproclaim
Using 20.9 and 20.6 we have
$$\align&\tr(T_x,E_v)=\sgn(x)\ov{\tr(T_x,E^\dag_v)}
=\ov{v^{-\aa_{E^\dag}}\tr(t_x,E^\dag_\spa)
+\text{strictly higher powers of } v}\\&
=v^{\aa_{E^\dag}}\tr(t_x,E^\dag_\spa)+\text{strictly lower powers of } v. 
\endalign$$
The proposition is proved.

\proclaim{Corollary 20.11}
$f_{E_v}=f_{E^\dag_\spa}v^{2\aa_{E^\dag}}+\text{strictly lower powers of } v$.
\endproclaim
Using 20.10 we have
$$\align &f_{E_v}\dim E=\sum_x\tr(T_x,E_v)\tr(T_{x\i},E_v) \\&\in
v^{2\aa_{E^\dag}}\sum_x\tr(t_x,E^\dag_\spa)\tr(t_{x\i},E^\dag_\spa)
+v^{2\aa_{E^\dag}-1}\bc[v\i]\\&
=v^{2\aa_{E^\dag}}f_{E^\dag_\spa}\dim E+v^{2\aa_{E^\dag}-1}\bc[v\i].
\endalign$$

\proclaim{Lemma 20.12} Let $E'\in\Irr W_I$. We have
$g_{E'_v}=\sum_{E;E\in\Irr W}[E':E]g_{E_v}$.
\endproclaim
We apply 19.4 with $A=\ch_{\bc(v)}$ and $A'$ the analogous algebra for $W_I$
instead of $W$, identified naturally with a subspace of $A$. (In this case the
unit elements of the two algebras are compatible hence $1_{A'}E_v=E_v$.) It
remains to use 20.5(a).

\proclaim{Lemma 20.13} Let $E\in\Irr W$.

(a) For any $x\in W$, $\tr(t_{x\i},E_\spa)$ is the complex conjugate of 
$\tr(t_x,E_\spa)$.

(b) $f_{E_\spa}$ is a strictly positive real number.
\endproclaim
We prove (a). Let $\lan,\ran$ be a positive definite hermitian form on $E$. We
define $\lan,\ran':E_\spa\tim E_\spa@>>>\bc$ by 

$\lan e,e'\ran'=\sum_{z\in W}\lan t_ze,t_ze'\ran$.
\nl
This is again a positive definite hermitian form on $E_\spa$. We show that

$\lan t_xe,e'\ran'=\lan e,t_{x\i}e'\ran'$
\nl
for all $e,e'$. This is equivalent to

$\sum_{y,z}\ga_{z,x,y\i}\lan t_ye,t_ze'\ran=
\sum_{y,z}\ga_{y,x\i,z\i}\lan t_ye,t_ze'\ran$
\nl
which follows from $\ga_{z,x,y\i}=\ga_{y,x\i,z\i}$. We see that $t_{x\i}$ is 
the adjoint of $t_x$ with respect to a positive definite hermitian form. (a) 
follows.

We prove (b). By 19.2(e) we have 
$f_{E_\spa}\dim(E)=\sum_x\tr(t_x,E_\spa)\tr(t_{x\i},E_\spa)$. The right 
hand side of this equality is a real number $\ge 0$, by (a). Hence so is the 
left hand side. Now $f_{E_\spa}\ne 0$ by 19.2(d) and (b) follows.

\proclaim{Proposition 20.14} Let $E'\in\Irr W_I$.

(a) For any $E\in\Irr W$ such that $[E':E]\ne 0$ we have
$\aa_{E'}\le\aa_E$.

(b) We have $g_{E'_\spa}=\sum [E':E]g_{E_\spa}$, sum over all 
$E\in\Irr W$ (up to isomorphism) such that $\aa_E=\aa_{E'}$.
\endproclaim
Let $X$ be the set of all $E$ (up to isomorphism) such that $[E':E]\ne 0$ and
such that $\aa_E$ is minimum, say equal to $a$. Assume first that $a<\aa_{E'}$.

Using 19.2(d) we rewrite 20.12 in the form 

(c) $v^{-2a}f_{E'_v}\i=\sum_E[E':E]v^{-2a}f_{E_v}\i$.
\nl
By 20.7, we have 

(d) $(v^{-2\aa_E}f_{E_v}\i)|_{v=0}=f_{E_\spa}\i$,
$(v^{-2\aa_{E'}}f_{E'_v}\i)|_{v=0}=f_{E'_\spa}\i$, 
\nl
hence by setting $v=0$ in (c) we obtain

$0=\sum_{E\in X}[E':E]f_{E_\spa}\i$.
\nl
The right hand side is a real number $>0$ by 20.8(b). This is a contradiction.
Thus we must have $a\ge\aa_{E'}$ and (a) is proved.

We now rewrite (c) in the form

(e) $v^{-2\aa_{E'}}f_{E'_v}\i=\sum_E[E':E]v^{-2\aa_{E'}}f_{E_v}\i$.
\nl
Using (d) and (a) we see that, setting $v=0$ in (e) gives

$f_{E'_\spa}\i=\sum_{E;\aa_E=\aa_{E'}}[E':E]f_{E_\spa}\i$.
\nl
This proves (b).

\subhead 20.15\endsubhead
Let $K(W)$ be the $\bc$-vector space with basis indexed by the $E\in\Irr W$ 
(up to isomorphism). If $\ti E\in\Mod W$ we identify $\ti E$ with the element
$\sum_E[E:\ti E]E\in K(W)$ ($E$ as above).

We define a $\bc$-linear map $\boj_{W_I}^W:K(W_I)@>>>K(W)$ by

$\boj_{W_I}^W(E')=\sum_E[E':E]E$,
\nl
sum over all $E\in\Irr W$ (up to isomorphism) such that $\aa_E=\aa_{E'}$;
here $E'\in\Irr W_I$. We call this {\it truncated induction}.

Let $I''\sub I'\sub S$. We show that the following transitivity formula holds:

(a) $\boj_{W_{I'}}^W \boj_{W_{I''}}^{W_{I'}}=\boj_{W_{I''}}^W:K(W_{I''})
@>>>K(W)$.
\nl
Let $E''\in\Irr W_{I''}$. We must show that

$[E'':E]=\sum_{E';\aa_{E'}=\aa_{E''}}[E'':E'][E':E]$
\nl
for any $E''\in\Irr W_{I''},E\in\Irr W$ such that $\aa_{E''}=\aa_E$;
in the sum we have $E'\in\Irr W_{I'}$. It is clear that

$[E'':E]=\sum_{E'}[E'':E'][E':E]$.
\nl
Hence it is enough to show that, if $[E'':E'][E':E]\ne 0$, then we have
automatically $\aa_{E'}=\aa_{E''}$. By 2.10(a) we have
$\aa_{E''}\le\aa_{E'}\le\aa_E$. Since $\aa_{E''}=\aa_E$, the desired 
conclusion follows.

\subhead 20.16\endsubhead
For any $x\in W$ define 

$\ga_x=\sum_{E;E\in\Irr W}\tr(t_x,E_\spa)E\in K(W)$.
\nl
We sometimes write $\ga_x^W$ instead of $\ga_x$, to emphasize dependence on 
$W$.

Note that $\ga_x$ is a $\bc$-linear combination of $E$ such that
$E\sim_{\cl\Cal R}x$. Hence, if $E,E'$ appear with $\ne 0$ coefficient in
$\ga_x$ then $E\sim_{\cl\Cal R}E'$. 

\proclaim{Proposition 20.17} If $x\in W_I$, then 
$\ga_x^W=\boj_{W_I}^W(\ga_x^{W_I})$.
\endproclaim
An equivalent statement is

(a) $\tr(t_x,E_\spa)=\sum_{E';\aa_E=\aa_{E'}}\tr(t_x,E'_\spa)[E':E]$
\nl
for any $E\in\Irr W$; in the sum we have $E'\in\Irr W_I$. Clearly, we have

(b) $v^{\aa_E}\tr(T_x,E_v)=\sum_{E';E'\in\Irr W_I}
v^{\aa_E}\tr(T_x,E'_v)[E':E]$.
\nl
In the right hand side we may assume that $\aa_{E'}\le\aa_E$. Using this and 
20.6, we see that setting $v=0$ in (b) gives (a). The proposition is proved. 

\proclaim{Lemma 20.18} (a) We have $\aa_{\sgn}=L(w_0)$.

(b) We have $f_{\sgn_\spa}=1$.

(c) We have $\ga_{w_0}=\sgn$. 
\endproclaim
$\sgn_v$ is the one dimensional $\ch_{\bc(v)}$-module on which $T_x$ acts as 
$\sgn(x)v^{-L(x)}$. (This follows from 20.3.) From 20.6(b) we see that 
$\aa_{\sgn}=L(w_0)$ and that $\tr(t_{w_0},\sgn_\spa)=1$. This proves (a). 

To prove (c) it remains to show that, if $\tr(t_{w_0},E_\spa)\ne 0$ ($E$ 
simple) then $E\cong\sgn$. This assumption shows, by 20.6(c), that 
$E_\spa\in\Irr J^\boc_\bc$ where $\boc$ is the two-sided cell such that 
$\sgn_\spa\in\Irr J^\boc_\bc$. Since 
$\tr(t_{w_0},\sgn_\spa)=1$, we have $w_0\in\boc$. From 13.8 it follows that 
$\{w_0\}$ is a two-sided cell. Thus $\boc=\{w_0\}$ and $J^\boc_\bc$ is one 
dimensional. Hence it cannot have more than one simple module. Thus, 
$E\cong\sgn$. This yields (c) and also (b). The lemma is proved.

\subhead 20.19\endsubhead
Assume that $I,I'$ form a partition of $S$ such that $W=W_I\tim W_{I'}$ (direct
product). If $E\in\Irr W_I$ and $E'\in\Irr W_{I'}$, then $E\boxt E'\in\Irr W$.
From the definitions one checks easily that

$\aa_{E\boxt E'}=\aa_E+\aa_{E'}, 
f_{(E\boxt E')_\spa}=f_{E_\spa}f_{E'_\spa}$.
\nl
Moreover, if $x\in W_I,x'\in W_{I'}$, then

$\ga_{xx'}^W=\ga_x^{W_I}\boxt\ga_{x'}^{W_{I'}}$.

\subhead 20.20\endsubhead
In the remainder of this section we assume that $w_0$ is in the centre of $W$.
Then, for any $E\in\Irr W$, $w_0$ acts on $E$ as $\ep_E$ times identity 
where $\ep_E=\pm 1$. Now $E\mto\ep_EE$ extends to a $\bc$-linear involution 
$\ze:K(W)@>>>K(W)$.

\proclaim{Lemma 20.21} Let $E\in\Irr W$. For any $x\in W$ we have

$\tr(T_{w_0x},E_v)=\ep_E v^{-\aa_E+\aa_{E^\dag}}\ov{\tr(T_x,E_v)}$.
\endproclaim
Since $w_0$ is in the centre of $W$, $T_{w_0}$ is in the centre of 
$\ch_{\bc(v)}$ hence it acts on $E_v$ as a scalar $\la\in\bc(v)$ times the
identity. Now $\tr(T_x,E_v)\in\bc[v,v\i]$ and $\tr(T_x\i,E_v)\in\bc[v,v\i]$. In
particular, $\la\in\bc[v,v\i]$ and $\la\i\in\bc[v,v\i]$. This implies 
$\la=cv^n$ where $c\in\bc$. For $v=1$, $\la$ becomes $\ep_E$. Hence 
$\la=\ep_Ev^n$ for some $n$. We have 

$\tr(T_{w_0x},E_v)=\tr(T_{w_0}T_{x\i}\i,E_v)=\la\tr(T_{x\i}\i,E_v)=
\la\ov{\tr(T_x,E_v)}$.
\nl
We have

$\sum_x\tr(T_{w_0x},E_v)\tr(T_{x\i w_0},E_v)=
\la^2\sum_x\ov{\tr(T_x,E_v)}\ov{\tr(T_{x\i},E_v)}$
\nl
hence $f_{E_v}\dim(E)=\la^2\ov{f_{E_v}}\dim(E)$ so that
$f_{E_v}=v^{2n}\ov{f_{E_v}}$. By 20.9, we have

$\sum_x\ov{\tr(T_x,E_v)}\ov{\tr(T_{x\i},E_v)}=
\sum_x\tr(T_x,E^\dag_v)\tr(T_{x\i},E^\dag_v)$
\nl
hence $\ov{f_{E_v}}=f_{E^\dag_v}$. We see that $f_{E_v}=v^{2n}f_{(E^\dag)_v}$. 
Comparing the lowest terms we see that

$-2\aa_E=2n-2\aa_{E^\dag}$ hence $n=-\aa_E+\aa_{E^\dag}$
\nl
and that 

(a) $f_{E_\spa}=f_{E^\dag_\spa}$.

\proclaim{Lemma 20.22} $v^{\aa_E}\tr(T_{w_0x},E_v)
=\ep_E(-1)^{l(x)}v^{\aa_{E^\dag}}\tr(T_x,E^\dag_v)$.
\endproclaim
We combine 20.8, 20.21.

\proclaim{Lemma 20.23} For any $x\in W$ we have 
$\ga_{xw_0}=\sgn(x)\ze(\ga_x)\ot\sgn$.
\endproclaim
An equivalent statement is

$\tr(t_{xw_0},E_\spa)=\sgn(x)\tr(t_x,E^\dag_\spa)\ep_{E^\dag}$
\nl
for any $E\in\Irr W$. Setting $v=0$ in the identity in 20.22 gives 

$\sgn(xw_0)\tr(t_{w_0x},E_\spa)=\ep_E\tr(t_x,E^\dag_\spa)$.
\nl
It remains to show that $\ep_{E^\dag}=\ep_E\sgn(w_0)$. This is clear.

\subhead 20.24\endsubhead
By the Cayley-Hamilton theorem, any element $r\in J$ satisfies an equation of 
the form $r^n+a_1r^{n-1}+\dots+a_n=0$ where $a_i\in\bz$. (We use that the 
structure constants of $J$ are integers.) This holds in particular for $r=t_x$
where $x\in W$. Hence for any $\ce\in\Irr J_\bc$, $\tr(t_x,\ce)$ is an 
algebraic integer. If $R$ is a subfield of $\bc$ such that the group algebra
$R[W]$ is split over $R$, then $J_R$ is split over $R$ and it follows that for
$x,\ce$ as above,  $\tr(t_x,\ce)$ is an algebraic integer in $R$. In 
particular, if we can take $R=\bq$, then $\tr(t_x,\ce)\in\bz$.

\head 21. Study of a left cell\endhead
\subhead 21.1\endsubhead
In this section we preserve the setup of 20.1. 
Let $\Ga$ be a left cell of $W,L$. Let $d$ be the unique element in 
$\Ga\cap\cd$. The $\ca$-submodule 
$\sum_{y\in\Ga}\ca c_y^\dag$ of $\ch$ can be regarded as an $\ch$-module by the
rule $c_x^\dag\cdot c_w^\dag=\sum_{z\in\Ga}h_{x,y,z}c_z^\dag$ with 
$x\in W,y\in W$. By change of scalars ($v\mto 1$) this gives rise to an 
$\ch_\bc=\bc[W]$-module $[\Ga]$. On 
the other hand, $J^\Ga_\bc=\opl_{y\in\Ga}\bc t_y$ is a left ideal in $J_\bc$ by
14.2(P8).

\proclaim{Lemma 21.2} The $\bc$-linear isomorphism $t_y\mto\hat n_yc_y^\dag$ 
for $y\in\Ga$ is an isomorphism of $J_\bc$-modules 
$J^\Ga_\bc@>\sim>>[\Ga]_\spa$.
\endproclaim
We have $\Ga\sub X_a=\{w\in W|\aa(x)=a\}$ for some $a\in\bn$.
The $\ca$-submodule $\sum_{y\in X_a}\ca c_y^\dag$ of $\ch$ can be regarded as 
an $\ch$-module by the rule 

$c_x^\dag\cdot c_w^\dag=
\sum_{z\in X_a}h_{x,y,z}c_z^\dag$
\nl
with $x\in W,y\in W$. By change of scalars ($v\mto 1$) this gives rise to an 
$\ch_\bc=\bc[W]$-module $[X_a]$. On the other hand, 
$J^{X_a}_\bc=\opl_{y\in X_a}\bc t_y$ is a left (even two-sided) ideal in 
$J_\bc$. The $\bc$-linear map in the lemma extends by the same formula to a
$\bc$-linear isomorphism $J^{X_a}_\bc@>\sim>>[X_a]_\spa$. It is enough to
show that this is $J_\bc$-linear. This follows from the computation in 18.10.
The lemma is proved.

\proclaim{Lemma 21.3} Let $\ce\in\Irr J_\bc$. The $\bc$-linear map 
$u:\Hom_{J_\bc}(J^\Ga_\bc,\ce)@>>>t_d\ce$ given by $\xi\mto\xi(n_dt_d)$ is an 
isomorphism.
\endproclaim
$u$ is well defined since $\xi(n_dt_d)=t_d\xi(t_d)\in t_d\ce$. We define a 
linear map in the opposite direction by $e\mto[j\mto je]$. It is clear that 
this is the inverse of $u$. (We use that $jn_dt_d=j$ for $j\in J^\Ga_\bc$.) The
lemma is proved.

\proclaim{Proposition 21.4} We have $\ga_d=n_d\sum_E[E:[\Ga]]E$ (sum over 
all $E\in\Irr W$ up to isomorphism).
\endproclaim
An equivalent statement is that $\tr(n_dt_d,E_\spa)=[E:[\Ga]]$, for $E$ as 
above. By 21.2, we have $[E:[\Ga]]=[E_\spa:J^\Ga_\bc]$. Hence it remains to 
show that $\tr(n_dt_d,\ce)=[\ce:J^\Ga_\bc]$ for any $\ce\in\Irr J_\bc$. Since 
$\ce=\opl_{d'\in\cd}n_{d'}t_{d'}\ce$ and $n_dt_d:\ce@>>>\ce$ is the projection
to the summand $n_dt_d\ce$, we see that $\tr(n_dt_d,\ce)=\dim(t_d\ce)$. It
remains to show that $\dim(t_d\ce)=[\ce:J^\Ga_\bc]$. This follows from 21.3.

\proclaim{Proposition 21.5} $[\Ga]^\dag,[\Ga w_0]$ are isomorphic in $\Mod W$.
\endproclaim
We may identify $[\Ga]^\dag$ with the $W$-module with $\bc$-basis 
$e_y (y\in\Ga)$ where $s\in S$ acts by 
$e_y\mto-e_y+\sum_{z\in\Ga}h_{s,y,z}e_z$. 

On the other hand we may identify $[\Ga w_0]$ with the $W$-module with 
$\bc$-basis $e'_{yw_0} (y\in\Ga)$ where $s\in S$ acts by 
$e'_{yw_0}\mto e'_{yw_0}-\sum_{z\in\Ga}h_{s,yw_0,zw_0}e'_{zw_0}$.

The $W$-module dual to $[\Ga]^\dag$ has a $\bc$-basis $e''_y (y\in\Ga)$ 
(dual to $(e_y)$) in which the action of $s\in S$ is given by
$e''_y\mto-e''_y+\sum_{z\in\Ga}h_{s,z,y}e''_z$. We define a $\bc$-isomorphism 
between this last space and $[\Ga w_0]$ by $e''_y\mto\sgn(y)e'_{yw_0}$ for all 
$y$. We show that this comutes with the action of $W$. It suffices to show that
for any $s\in S$, we have

(a) $-h_{s,z,y}=\sgn(y)\sgn(z)h_{s,yw_0,zw_0}$ for all $z\ne y$ and

(b) $1-h_{s,y,y}=-1+h_{s,yw_0,yw_0}$ for all $y$.
\nl
We use 6.6. Assume first that $sz>z$. If $sy>y$ and $y\ne z$, both sides of (a)
are $0$. If $sy<y<z$ then (a) follows from 11.6. If $y=sz$ then  both sides of
(a) are $-1$. If $sy<y$ but $y\not<z$ or $y\ne sz$ then both sides of (a) are 
$0$.

Assume next that $sz<z$. If $z\ne y$ then both sides of (a) are $0$.

If $sy>y$, both sides of (b) are $1$. If $sy<y$, both sides of (b) are $-1$.
Thus (a),(b) are verified. Since $[\Ga]^\dag$ and its dual are isomorphic 
in $\Mod W$ (they are defined over $\bq$), the lemma follows.

\proclaim{Corollary 21.6} Let $E\in\Irr W$ and let $\boc$ be the two-sided
cell of $W$ such that $E_\spa\in\Irr J^{\boc}_\bc$. Then 
$E^\dag_\spa\in\Irr J^{\boc w_0}_\bc$.
\endproclaim
Replacing $\Ga$ by $\boc$ in the definition of $[\Ga]$ we obtain a 
$W$-module $[\boc]$. Then 21.2, 21.5 hold with $\Ga$ replaced by $\boc$
with the same proof. Our assumption implies (by 21.2 for $\boc$) that $E$ 
appears in the $W$-module $[\boc]$. Using 21.5 for $\boc$ we deduce that 
$E^\dag$ 
appears in the $W$-module $[\boc w_0]$. Using 21.2 for $\boc w_0$,
we deduce that $E^\dag_\spa$ appears in the $J_\bc$-module $J^{\boc w_0}_\bc$.
The corollary follows. 

\proclaim{Corollary 21.7} Let $E,E'\in\Irr W$ be such that
$E\sim_{\cl\Cal R}E'$. Then $E^\dag\sim_{\cl\Cal R}E'{}^\dag$.
\endproclaim
By assumption there exists a two-sided cell $\boc$ such that $E_\spa,
E'_\spa\in\Irr J^{\boc}_\bc$. By 21.6, 
$E^\dag_\spa,E'{}^\dag_\spa\in\Irr J^{\boc w_0}_\bc$. The corollary follows.

\subhead 21.8\endsubhead
The results of \S19 are applicable to $A$, the $\bc$-subspace 
$J^{\Ga\cap\Ga\i}_\bc$ of $J_\bc$ spanned by $\Ga\cap\Ga\i$ and $R=\bc$. This 
is a $\bc$-subalgebra of $J_\bc$ with unit element $t_d$. In 21.9 we will show
that $J^{\Ga\cap\Ga\i}_\bc$ is semisimple. It is then clearly split. We take 
$\tau:J^{\Ga\cap\Ga\i}_\bc@>>>\bc$ so that $\tau(t_x)=\de_{x,d}$. (This is the
restriction of $\tau:J_\bc@>>>\bc$.) We have $(t_x,t_y)=\de_{xy,1}$. The bases
$(t_x)$ and $(t_{x\i})$ (where $x$ runs through $\Ga\cap\Ga\i$) are dual with 
respect to $(,)$. 

\subhead 21.9\endsubhead
We show that the $\bc$-algebra $J^{\Ga\cap\Ga\i}_\bc$ is semisimple. It is 
enough to prove the analogous statement for the $\bq$-algebra $A'$, the
$\bq$-span of $\Ga\cap\Ga\i$ in $J_\bq$. We define a $\bq$-bilinear pairing 
$(|):A'\tim A'\to\bq$ by $(t_x|t_y)=\de_{x,y}$  for $x,y\in\Ga\cap\Ga\i$. Let 
$j\mto\ti j$ be the $\bq$-linear map $A'\to A'$ given by $\ti t_x=t_{x\i}$ for
all $x$. We show that 

(a) $(j_1j_2|j_3)=(j_2|\ti j_1j_3)$
\nl
for all $j_1,j_2,j_3$ in our ring. We may assume that 
$j_1=t_x,j_2=t_y,j_3=t_z$. Then (a) follows from 

$\ga_{x,y,z\i}=\ga_{x\i,z,y\i}$.
\nl
Now let $I$ be a left ideal of $A'$. Let
$I^\perp=\{a\in A'| (a|I)=0\}$. Since $(|)$ is positive definite, we have
$A'=I\opl I^\perp$. From (c) we see that $I^\perp$ is a left ideal. This proves
that $A'$ is semisimple.

The same proof could be used to show directly that $J_\bc$ is semisimple.

\proclaim{Proposition 21.10} Let $E,E'\in\Irr W$. Let
$N=\sum_{x\in\Ga\cap\Ga\i}\tr(t_x,E_\spa)\tr(t_{x\i},E'_\spa)$. Then 
$N=f_{E_\spa}[E:[\Ga]]$ if $E,E'$ are isomorphic and $N=0$, otherwise.
\endproclaim
If $\ce\in\Irr J_\bc$, then $t_d\ce$ is either $0$ or in 
$\Irr J^{\Ga\cap\Ga\i}_\bc$. Moreover, $\ce\mto t_d\ce$ defines a bijection 
between the set of simple $J_\bc$-modules (up to isomorphism) which appear in 
the $J_\bc$-module $J^\Ga_\bc$ and the set of simple 
$J^{\Ga\cap\Ga\i}_\bc$-modules (up to isomorphism). We then have
$\dim(t_d\ce)=[\ce:J_\bc^\Ga]$. Note that, for $j\in J^{\Ga\cap\Ga\i}_\bc$ we 
have $\tr(j,\ce)=\tr(j,t_d\ce)$. If $t_dE_\spa=0$ or $t_dE'_\spa=0$, then
$N=0$ and the result is clear. If $t_dE_\spa\ne 0$ and $t'_dE_\spa\ne 0$ 
then, by 19.2(e), we see that $N=f_{t_dE_\spa}[E_\spa:J^\Ga_\bc]$ if $E,E'$
are isomorphic and to $0$, otherwise. It remains to show that 
$f_{t_dE_\spa}=f_{E_\spa}$, $[E:[\Ga]]=[E_\spa:J^\Ga_\bc]$ and the 
analogous equalities for $E'$. Now $f_{t_dE_\spa}=f_{E_\spa}$ follows from
19.4 applied to $(A',A)=(J^{\Ga\cap\Ga\i}_\bc,J_\bc)$; the equality 
$[E:[\Ga]]=[E_\spa:J^\Ga_\bc]$ follows from 21.2. The proposition is proved.

\head 22. Constructible representations\endhead
\subhead 22.1\endsubhead
In this section we preserve the setup of 20.1. 

We define a class $Con(W)$ of $W$-modules (relative to our fixed 
$L:W@>>>\bn$) by induction on $|S|$. If $|S|=0$ so that $W=\{1\}$, $Con(W)$
consists of the unit representation. Assume now that $|S|>0$. Then $Con(W)$ 
consists of the $W$-modules of the form $\boj_{W_I}^W(E')$ or 
$\boj_{W_I}^W(E')\ot\sgn$ for various subsets $I\sub S, I\ne S$ and various
$E'\in Con(W_I)$. (If we restrict ourselves to $I$ such that $|I|=|S|-1$ we get
the same class of $W$-modules, by the transitivity of truncated 
induction.) The $W$-modules in $Con(W)$ are said to be the 
{\it constructible representations} of $W$.

Note that the unit representation of $W$ is constructible (it is obtained by
truncated induction from the unit representation of the subgroup with one
element). Hence $\sgn\in Con(W)$.

\proclaim{Lemma 22.2} If $E\in Con(W)$, then there exists a left cell $\Ga$ of
$W$ such that $E=[\Ga]$.
\endproclaim
We argue by induction on $|S|$. If $|S|=0$ the result is obvious. Assume now 
that $|S|>0$. Let $E\in Con(W)$.

{\it Case 1}. $E=\boj_{W_I}^W(E')$ where $I\sub S,I\ne S$ and $E'\in Con(W_I)$.
By the induction hypothesis there exists a left cell $\Ga'$ of $W_I$ such that
$E'=[\Ga']$. Let $d\in\Ga'\cap\cd$. By 21.4 we have $\ga_d^{W_I}=[\Ga']=E'$. By
20.17 we have $E=\boj_{W_I}^W(E')=\boj_{W_I}^W(\ga_d^{W_I})=\ga_d^W$. Let $\Ga$
be the left cell of $W$ that contains $d$. By 21.4 we have $\ga_d^W=[\Ga]$. 
Hence $E=[\Ga]$.

{\it Case 2}. $E=\boj_{W_I}^W(E')\ot\sgn$ where $I\sub S,I\ne S$ and 
$E'\in Con(W_I)$. Then by Case 1, $E\ot\sgn=[\Ga]$ for some left cell $\Ga$ of
$W$. By 21.5 we have $E=[\Ga w_0]$. The lemma is proved.

\proclaim{Proposition 22.3} For any $E\in\Irr W$ there exists a constructible 
representation of $W$ which contains a simple component isomorphic to $E$.
\endproclaim
The general case can be easily reduced to the case where $W$ is irreducible.
Assume now that $W$ is irreducible. If $L=al$ for some $a>0$, the constructible
representations of $W$ are listed in \cite{...} and the proposition is easily 
checked.
(See also the discussion of types $A,D$ in 22.5, 22.26.) In the cases where $W$
is irreducible but $L$ is not of the form $al$, the constructible 
representations are described later in this section and this yields the 
proposition in all cases.

\subhead 22.4\endsubhead
Let $W=\fS_n$ be the group of permutations of $1,2,\dots,n$. We regard $W$ as 
a Coxeter group with generators

$s_1=(1,2),s_2=(2,3),\dots,s_{n-1}=(n-1,n)$,
\nl
(transpositions). We take $L=al$ where $a>0$.

The simple $W$-modules (up to isomorphism) are in 1-1 correspondence with
the partitions $\al=(\al_1\ge\al_2\ge\dots)$ such that $\al_N=0$ for large $N$
and $\sum_i\al_i=n$. The correspondence (denoted by $\al\mto\pi_\al$) is 
defined as follows. Let $\al$ be as above, let $(\al'_1\ge\al'_2\ge\dots)$ be 
the partition dual to $\al$. Let $\pi_\al$ be the simple $W$-module whose
restriction to $\fS_{\al_1}\tim\fS_{\al_2}\dots$ contains $1$ and whose 
restriction to $\fS_{\al'_1}\tim\fS_{\al'_2}\dots$ contains the sign 
representation. We have (a consequence of results of Steinberg):
$$f_{(\pi_\al)_v}=v^{-\sum_i2\binom{al'_i}{2}}+\text{strictly higher powers of}
v.$$
It follows that 

(a) $\aa_{\pi_\al}=\sum_ia\binom{\al'_i}{2}$ and $f_{(\pi_\al)_\spa}=1$.

\proclaim{Lemma 22.5} In the setup of 22.4, a $W$-module is constructible
if and only if it is simple.
\endproclaim
For any sequence $\beta=(\beta_1,\beta_2,\dots)$ in $\bn$ such that $\beta_N=0$
for large $N$ and $\sum_i\beta_i=n$, we set

$I_{\beta}=\{s_i|i\in[1,n-1],i\ne\beta_1,i\ne\beta_1+\beta_2,\dots\}$.
\nl
From 22.4(a) we see easily that, if $\beta$ is the same as $\al'$ up to order, 
then

(a) $\boj_{W_{I_\beta}}^W(\sgn)=\pi_\al$. 
\nl
Since the $\sgn\in Con(W_{I_\beta})$, it follows that $\pi_\al\in Con(W)$. Thus
any simple $W$-module is constructible.

We now show that any constructible representation $E$ of $W=\fS_n$ is simple. 
We may assume that $n\ge 1$ and that the analogous result is true for any 
$W_{I'}\ne W$. We may assume that $E=\boj_{W_{I_\beta}}^W(C)$ where $\beta$ is
as above, $W_{I_\beta}\ne W$ and $C\in Con(W_{I_\beta})$. By the induction
hypothesis, $C$ is simple. Since the analogue of (a) holds for $W_{I_\beta}$
(instead of $W$) we have $C=\boj_{W_{I_{\beta'}}}^{W_{I_\beta}}(\sgn)$ for some
$\beta'$ such that $W_{I_{\beta'}}\sub W_{I_\beta}$. By the transitivity of
truncated induction we have $E=\boj_{W_{I_{\beta'}}}^W(\sgn)$. Hence, by (a),
for $\beta'$ instead of $\beta$, $E$ is simple. The lemma is proved.

\subhead 22.6\endsubhead
We now develop some combinatorics which is useful for the verification of 22.3
for W of classical type.

Let $a>0,b\ge 0$ be integers. We can write uniquely $b=ar+b'$ where 
$r,b'\in\bn$ and $b'<a$. Let $N\in\bn$. Let $\cm_{a,b}^N$ be the set of 
multisets $\ti Z=\{\ti z_1\le\ti z_2\le\dots\le\ti z_{2N+r}\}$ of integers 
$\ge 0$ such that 

(a) if $b'=0$, there are at least $N+r$ distinct entries in $\ti Z$, no entry
is repeated more than twice and all entries of $\ti Z$ are divisible by $a$;

(b) if $b'>0$, all inequalities in $\ti Z$ are strict and $N$ entries of 
$\ti Z$ are divisible by $a$ and $N+r$ entries of $\ti Z$ are equal to $b'$
modulo $a$.
\nl
The entries which appear in $\ti Z$ exactly once are called the {\it singles}
of $\ti Z$; they form a set $Z$. The other entries of $\ti Z$ are called the
{\it doubles} of $\ti Z$.

For example, the multiset $\ti Z^0$ whose entries are (up to order)

$0,a,2a,\dots,(N-1)a,b',a+b',2a+b',\dots,(N+r-1)a+b'$
\nl
belongs to $\cm_{a,b}^N$. Clearly, the sum of entries of $\ti Z$ minus the sum
of entries of $\ti Z^0$ is $\ge 0$ and divisible by $a$, hence it is equal to 
$an$ for a well defined $n\in\bn$ said to be the {\it rank} of $\ti Z$. We have
$$\sum_{k=1}^{2N+r}\ti z_k=an+aN^2+N(b-a)+a\binom{r}{2}+b'r.$$
Note that $\ti Z^0$ has rank $0$. Let $\cm_{a,b;n}^N$ be the set of multisets
of rank $n$ in $\cm_{a,b}^N$. We define an (injective) map 
$\cm_{a,b}^N@>>>\cm_{a,b}^{N+1}$ by 

$\{\ti z_1\le\ti z_2\le\dots\ti z_{2N+r}\}\mto
\{0,b',\ti z_1+a\le\ti z_2+a\le\dots\ti z_{2N+r}+a\}$.
\nl
This restricts for any $n\in\bn$ to an (injective) map 

(c) $\cm_{a,b;n}^N@>>>\cm_{a,b;n}^{N+1}$.
\nl
It is easy to see that, for fixed $n$, $|\cm_{a,b;n}^N|$ is bounded as 
$N\to\infty$, hence the maps (c) are bijections for large $N$. Let 
$\cm_{a,b;n}$ be the inductive limit of $\cm_{a,b;n}^N$ as $N\to\infty$ (with 
respect to the maps (c)).

\subhead 22.7\endsubhead
Let $\Sy_{a,b;n}^N$ be the set consisting of all tableaux (or {\it symbols})
$$\align\la_1&,\la_2,\dots,\la_{N+r}\\&\mu_1,\mu_2,\dots,\mu_N\tag a\endalign$$
where $\la_1<\la_2<\dots<\la_{N+r}$ are integers $\ge 0$, congruent to $b'$
modulo $a$, $\mu_1,\mu_2,\dots,\mu_N$ are integers $\ge 0$, divisible by $a$
and
$$\sum_i\la_i+\sum_j\mu_j=an+aN^2+N(b-a)+a\binom{r}{2}+b'r.$$
If we arrange the entries of $\La$ in a single row, we obtain a multiset
$\ti Z\in\cm_{a,b;n}^N$. This defines a (surjective) map 
$\pi_N:\Sy_{a,b;n}^N@>>>\cm_{a,b;n}^N$.

We define an (injective) map 

(b) $\Sy_{a,b;n}^N@>>>\Sy_{a,b;n}^{N+1}$
\nl
by associating to (a) the symbol
$$\align b',\la_1+a&,\la_2+a,\dots,\la_{N+r}+a\\&
0,\mu_1+a,\mu_2+a,\dots,\mu_N+a.\endalign$$
This is compatible with the map $\cm_{a,b}^N@>>>\cm_{a,b}^{N+1}$ in 22.6 (via 
$\pi_N,\pi_{N+1}$).

Since for fixed $n$, $|\Sy_{a,b;n}^N|$ is bounded as $N\to\infty$, the maps (b)
are bijections for large $N$. Let $\Sy_{a,b;n}$ be the inductive limit of 
$\Sy_{a,b;n}^N$ as $N\to\infty$ (with respect to the maps (b)).

\subhead 22.8\endsubhead
Let $\ti Z=\{\ti z_1\le\ti z_2\le\dots\ti z_{2N+r}\}\in\cm_{a,b;n}^N$. Let $t$
be an integer which is large enough so that the multiset

(a) $\{at+b'-\ti z_1,at+b'-\ti z_2,\dots,at+b'-\ti z_{2N+r}\}$
\nl
is contained in the multiset

(b) $\{0,a,2a,\dots,ta, b',a+b',2a+b',\dots,ta+b'\}$
\nl
and let $\bar{\ti Z}$ be the complement of (a) in (b). Then
$\bar{\ti Z}\in\cm_{a,b}^{t+1-N-r}$. The sum of entries of $\bar{\ti Z}$ is
$$\align&\sum_{k\in[0,t]}(2ka+b')-(at+b')(2N+r)+\sum_h\ti z_h\\&=
at(t+1)+(t+1)b'-(at+b')(2N+r)+an+aN^2+N(b-a)+a\binom{r}{2}+b'r\\&
=an+a(t+1-N-r)^2+(t+1-N-r)(b-a)+a\binom{r}{2}+b'r.\endalign$$
Thus, $\bar{\ti Z}$ has rank $n$.

We define a bijection $\pi_N\i(\ti Z)@>\sim>>\pi_{t+1-N-r}\i(\bar{\ti Z})$ by
$\La\mto\bar\La$ where $\La$ is as in 22.7(a) and $\bar\La$ is
$$\align&\{b',a+b',2a+b',3a+b',\dots,ta+b'\}
-\{at+b'-\mu_1,at+b'-\mu_2,\dots,at+b'-\mu_N\}\\&\{0,a,2a,3a,\dots,ta\}
-\{at+b'-\la_1,at+b'-\la_2,\dots,at+b'-\la_{N+r}\}.\endalign$$ 

\subhead 22.9\endsubhead
Let $W=W_n$ be the group of permutations of $1,2,\dots,n,n',\dots,2',1'$ which
commute with the involution $i\mto i', i'\mto i$. We regard $W_n$ as a Coxeter
group with generators $s_1,s_2,\dots,s_n$ given as products of transpositions 
by 

$s_1=(1,2)(1',2'),s_2=(2,3)(2',3'),\dots,s_{n-1}=(n-1,n)((n-1)',n')$,

$s_n=(n,n')$.

\subhead 22.10\endsubhead
A permutation in $W$ defines a permutation of the $n$ element set consisting of
the pairs $(1,1'),(2,2'),\dots,(n,n')$. Thus we have a natural homomorphism of
$W_n$ onto $\fS_n$, the symmetric group in $n$ letters. Let 
$\chi_n:W_n@>>>\pm 1$ be the homomorphism defined by

$\chi_n(\si)=1$ if $\{\si(1),\si(2),\dots,\si(n)\}\cap\{1',2',\dots,n'\}$ has 
even cardinality,

$\chi_n(\si)=-1$, otherwise.
\nl
The simple $W$-modules (up to isomorphism) are in 1-1 correspondence with
the ordered pairs $\al,\beta$ where $\al=(\al_1\ge\al_2\ge\dots)$ and
$\beta=(\beta_1\ge\beta_2\ge\dots)$ are partitions such that $\al_N=\be_N=0$ 
for large $N$ and $\sum_i\al_i+\sum_j\beta_j=n$. The correspondence (denoted by
$\al,\beta\mto E^{\al,\beta}$) is defined as follows. Let $\al,\beta$ be as 
above, let $(\al'_1\ge\al'_2\ge\dots)$ be the partition dual to $\al$ and let 
$(\beta'_1\ge\beta'_2\ge\dots)$ be the partition dual to $\beta$. Let 
$k=\sum_i\al_i,l=\sum_j\beta_j$. Let $\pi_\al$ be the simple 
$\fS_k$-module defined as in 22.4 and let $\pi_\beta$ be the analogously 
defined simple $\fS_l$-module. We regard $\pi_\al,\pi_\beta$ as simple
modules of $W_k,W_l$ via the natural homomorphisms $W_k@>>>\fS_k,W_l@>>>\fS_l$
as above. We identify $W_k\tim W_l$ with the subgroup of $W$ consisting of all 
permutations in $W$ which map $1,2,\dots,k,k',\dots,2',1'$ into itself hence 
also map $k+1,k+2,\dots,n,n',\dots,(k+2)',(k+1)'$ into itself. Consider the
representation $\pi_\al\ot(\pi_\beta\ot\chi_l)$ of $W_k\tim W_l$. We induce it
to $W$; the resulting representation of $W$ is irreducible; we denote it by
$E^{\al,\beta}$. 

We fix $a>0,b\ge 0$ and we write $b=ar+b'$ as in 22.6.

Let $\al,\beta$ be as in 22.10. Let $N$ be an integer such that 
$\al_{N+r+1}=0,\beta_{N+1}=0$. (Any large enough integer satisfies these 
conditions.) We set
$$\la_i=a(\al_{N+r-i+1}+i-1)+b',(i\in[1,N+r]),\quad
\mu_j=a(\beta_{N-j+1}+j-1),(j\in[1,N]).$$    
We have $0\le\la_1<\la_2<\dots<\la_{N+r}$, $0\le\mu_1<\mu_2<\dots<\mu_N$. Let 
$\La$ denote the tableau 22.7(a).
It is easy to see that $\La\in\Sy_{a,b;n}^N$. Moreover, if $N$ is replaced by
$N+1$, then $\La$ is replaced by its image under
$\Sy_{a,b;n}^N@>>>\Sy_{a,b;n}^{N+1}$ (see 22.7). Let $[\La]=E^{\al,\beta}$. 
Note that $[\La]$ depends only on the image of $\La$ under the canonical map 
$\Sy_{a,b;n}^N@>>>\Sy_{a,b;n}$. In this way, we see that 

{\it the simple $W$-modules are naturally in bijection with the set 
$\Sy_{a,b;n}$.}
\nl
For $i\in[1,N]$ we have
$a(\al_{N-i+1}+i-1)+b=a(\al_{N+r-i-r+1}+i+r-1)+b'=\la_{i+r}$. 

If $N$ is large we have $\la_i=a(i-1)+b'$ for $i\in[1,r]$ and $\mu_j=a(j-1)$ 
for $j\in[1,r]$.

\subhead 22.11\endsubhead
Let $q,y$ be indeterminates. With the notation in 22.15, let
$$H_\al(q)=q^{-\sum_i\binom{\al'_i}{2}}
\prod_{i,j}\fra{q^{\al_i+\al'_j-i-j+1}-1}{q-1},$$
$$G_{\al,\beta}(q,y)=q^{-\sum_i\al'_i\beta'_i/2}
\prod_{i,j}(q^{\al_i+\beta'_j-i-j+1}y+1);$$
both products are taken over all $i\ge 1,j\ge 1$ such that
$\al_i\ge j,\al'_j\ge i$. 

Let $L:W@>>>\bn$ be the weight function defined be defined by 
$L(s_1)=L(s_2)=\dots=L(s_{n-1})=a$, $L(s_n)=b$. We now assume that both $a,b$
are $>0$. We also assume that $a,b$ are such that $W,L$ satisfies the 
assumptions of 18.1. Then $f_{E^{\al,\beta}_v}$ is defined in terms of this 
$L$.

\proclaim{Lemma 22.12 (Hoefsmit \cite{\HO})} We have
$$f_{E^{\al,\beta}_v}=H_\al(v^{2a})H_\beta(v^{2a})G_{\al,\beta}(v^{2a},v^{2b})
G_{\beta,\al}(v^{2a},v^{-2b}).$$
\endproclaim

We will rewrite the expression above using the following result.

\proclaim{Lemma 22.13} Let $N$ be a large integer. We have
$$H_\al(q)=q^{\sum_{i\in[1,N-1]}\binom{i}{2}}
\fra{\prod_{i=1}^N \prod_{h\in[1,\al_{N-i+1}+i-1]}\fra{q^h-1}{q-1}}
{\prod_{1\le i<j\le N}\fra{q^{\al_{N-j+1}+j-1}-q^{\al_{N-i+1}+i-1}}{q-1}},$$
$$\align&G_{\al,\beta}(q,y)G_{\beta,\al}(q,y\i)=
q^{\sum_{i\in[1,N-1]}i^2}(\sqrt y+\sqrt y\i)^N\\&
\tim\fra{\prod_{i=1}^N\prod_{h\in[1,\al_{N-i+1}+i-1]}(q^hy+1)\prod_{j=1}^N
\prod_{h\in[1,\beta_{N-j+1}+j-1]}(q^hy\i+1)}{\prod_{i,j\in[1,N]}
(q^{\al_{N-i+1}+i-1}\sqrt y+q^{\beta_{N-j+1}+j-1}\sqrt y\i)}.\endalign$$
\endproclaim
The proof is by induction on $n$. We omit it.

\proclaim{Proposition 22.14} (a) If $b'=0$ then $f_{[\La]_\spa}$ is equal to 
$2^d$ where $2d+r$ is the number of singles in $\La$. If $b'>0$ then 
$f_{[\La]_\spa}=1$.

(b) We have $\aa_{[\La]}=A_N-B_N$ where
$$A_N=\sum_{i\in[1,N+r],j\in[1,N]}\inf(\la_i,\mu_j)
+\sum_{1\le i<j\le N+r}\inf(\la_i,\la_j)+\sum_{1\le i<j\le N}\inf(\mu_i,\mu_j),
$$
$$\align&B_N=\sum_{i\in[1,N+r],j\in[1,N]}\inf(a(i-1)+b',a(j-1))\\&
+\sum_{1\le i<j\le N+r}\inf(a(i-1)+b',a(j-1)+b')
+\sum_{1\le i<j\le N}\inf(a(i-1),a(j-1)).\endalign$$
\endproclaim
It is enough to prove (a) assuming that $N$ is large. Since

$A_{N+1}-A_N=a(N+r)N+a\binom{N}{2}+a\binom{N+r}{2}+b'N+b'(N+r)=B_{N+1}-B_N$,
\nl
we have $A_{N+1}-B_{N+1}=A_N-B_N$ hence it is enough to prove (b) assuming that
$N$ is large. In the remainder of the proof we assume that $N$ is large.

For $f,f'\in\bc(v)$ we write $f\cong f'$ if $f'=fg$ with $g\in\bc(v)$,
$g|_{v=0}=1$. Using 22.12, 22.13, we see that
$$\align&f_{[\La]_v}\cong\prod_{i\in[1,N]}(v^{2a-2b}+1)(v^{4a-2b}+1)\dots
(v^{2\mu_i-2b}+1)\\&(v^b+v^{-b})^Nv^{2a\sum_{i=1}^{N-1}
(2i^2-i)}\prod_{i,j\in[1,N]}(v^{2\la_{i+r}-b}
+v^{2\mu_j-b})\i\\&\prod_{1\le i<j\le N}
(v^{2\la_{j+r}-2b}-v^{2\la_{i+r}-2b})\i\prod_{1\le i<j\le N}
(v^{2\mu_j}-v^{2\mu_i})\i\endalign$$
hence 
$$f_{[\La]_v}=2^dv^{-K}+\text{higher powers of $v$}$$
where

$d=0$ if $b'>0$,
$$\align&d=\sha(j\in[1,N]: b\le\mu_j)-\sha(i,j\in[1,N]: \la_{i+r}=\mu_j)\\&
=N-\sha(i\in[1,r], j\in[1,N]: (i-1)a=\mu_j)-\sha(i,j\in[1,N]: \la_{i+r}=\mu_j)
\\&=N-\sha(i\in[1,r],j\in[1,N]:\la_i=\mu_j)-\sha(i,j\in[1,N]:\la_{i+r}=\mu_j)
\\&=N-\sha(i\in[1,N+r], j\in[1,N]: \la_i=\mu_j)=(\sha singles-r)/2,
\text{ if } b'=0,\endalign$$
$$\align&-K=-bN+2a\sum_{i\in[1,N-1]}(2i^2-i)
+\sum_{j\in[1,N]}\sum\Sb k\in [1,r]\\ak\le\mu_j\endSb(2ak-2b)\\&
-\sum_{i,j\in[1,N]}(-b+2\inf(\la_{i+r},\mu_j))
-\sum_{1\le i<j\le N}(-2b+2\inf(\la_{i+r},\la_{j+r}))\\&
-\sum_{1\le i<j\le N}2\inf(\mu_i,\mu_j)
=-bN+2a\sum_{i\in[1,N-1]}(2i^2-i)+2bN^2-bN\\&
+\sum_{j\in[1,N]}\sum_{k\in[1,r],ak\le\mu_j}(2ak-2b)
-\sum_{i,j\in[1,N]}2\inf(\la_{i+r},\mu_j)\\&
-\sum_{1\le i<j\le N}2\inf(\la_{i+r},\la_{j+r})
-\sum_{1\le i<j\le N}2\inf(\mu_i,\mu_j)\\&
=\sum_{j\in[1,N]}\sum_{k\in[1,r],ak\le\mu_j}(2ak-2b)
-\sum_{i,j\in[1,N]}2\inf(\la_{i+r},\mu_j)\\&
-\sum_{1\le i<j\le N}2\inf(\la_{i+r},\la_{j+r})
-\sum_{1\le i<j\le N}2\inf(\mu_i,\mu_j)+\bst.\endalign$$
(We will generally write $\bst$ for an expression which depends only on 
$a,b,N$.) We have
$$\align&\sum\Sb j\in[1,N]\\k\in[1,r]\\ak\le\mu_j\endSb(2ak-2b)=
\sum\Sb j\in[1,r]\\k\in[1,r]\\ak\le\mu_j\endSb(2ak-2b)
+\sum\Sb j\in[r+1,N]\\k\in[1,r]\\ak\le\mu_j\endSb(2ak-2b)\\&
=\sum\Sb j\in[1,r]\\k\in[1,r]\\ak\le a(j-1)\endSb(2ak-2b)
+\sum\Sb j\in[r+1,N]\\k\in[1,r]\endSb(2ak-2b)=\bst,\endalign$$
hence
$$-K=-2(\sum_{i,j\in[1,N]}\inf(\la_{i+r},\mu_j)
-\sum_{1\le i<j\le N}\inf(\la_{i+r},\la_{j+r})
-\sum_{1\le i<j\le N}\inf(\mu_i,\mu_j))+\bst.$$
We have
$$\align&\sum_{i\in[1,r],j\in[1,N]}\inf(\la_i,\mu_j)
=\sum_{i\in[1,r],j\in[1,N]}\inf(a(i-1)+b',\mu_j)\\&
=\sum_{i\in[1,r],j\in[1,r]}\inf(a(i-1)+b',a(j-1))
+\sum_{i\in[1,r],j\in[r+1,N]}\inf(a(i-1)+b',\mu_j)\\&
=\sum_{i\in[1,r],j\in[1,r]}\inf(a(i-1)+b',a(j-1))
+\sum_{i\in[1,r],j\in[r+1,N]}(a(i-1)+b')=\bst,\endalign$$
hence
$$\sum_{i,j\in[1,N]}\inf(\la_{i+r},\mu_j)=
\sum_{i\in[1,N+r],j\in[1,N]}\inf(\la_i,\mu_j)+\bst.$$
We have
$$\align&\sum_{1\le i<j\le N+r}\inf(\la_i,\la_j)=
\sum_{1\le i<j\le N}\inf(\la_{i+r},\la_{j+r})+\sum_{i\in[1,r]}\la_i(N+r-i)\\&
=\sum_{1\le i<j\le N}\inf(\la_{i+r},\la_{j+r})
+\sum_{i\in[1,r]}(a(i-1)+b')(N+r-i)\\&
=\sum_{1\le i<j\le N}\inf(\la_{i+r},\la_{j+r})+\bst.\endalign$$
We see that
$$-K=-2A_N+\bst.\tag c$$
In the special case where $\al=\beta=(0\ge 0\ge\dots)$ we have $K=0$. On the 
other hand, by (c), we have $0=-2B_N+\bst$ where $\bst$ is as in (c). Hence in 
general we have $-K=-2A_N+2B_N$. This proves the proposition, in view of 20.11
and 20.21(a).

\subhead 22.15\endsubhead
We identify $\fS_k\tim W_l$ ($k+l=n)$ with the subgroup of $W$ consisting of
all permutations in $W$ which map $\{1,2,\dots,k\}$ into itself (hence also map
$\{1,2,\dots,k'\}$ and  $\{k+1,\dots,n,n'\dots,(k+1)'\}$ into themselves. This
is a standard parabolic subgroup of $W$. We consider an irreducible
representation of $\fS_k\tim W_l$ of the form $\sgn_k\boxt[\La']$ where 
$\sgn_k$ is the sign representation of $\fS_k$ and $\La'\in\Sy_{a,b;l}^N$. We
may assume that $\La'$ has at least $k$ entries. We want to associate to 
$\La'$ a symbol in $\Sy_{a,b;n}^N$ by increasing each of the $k$ largest 
entries in $\La'$ by $a$. It may happen that the set of $r$ largest entries of
$\La'$ is not uniquely defined but there are two choices for it. (This can only
happen if $b'=0$.) Then the same procedure gives rise to two distinct symbols 
$\La^I,\La^{II}$ in $\Sy_{a,b;n}^N$.

\proclaim{Lemma 22.16} (a) $g_{(\sgn_k\ot[\La'])_\spa}=g_{[\La']_\spa}$ is 
equal to $g_{[\La]_\spa}$ or to $g_{[\La^I]_\spa}+g_{[\La^{II}]_\spa}$.

(b) $\aa_{\sgn_k\ot[\La']}=a\binom{k}{2}+\aa_{[\La']}$ is equal to 
$\aa_{[\La]}$ or to $\aa_{[\La^I]}=\aa_{[\La^{II}]}$.
\endproclaim
$\La$, if defined, has the same number of singles as $\La'$. Moreover, $\La^I$
(and $\La^{II}$), if defined, has one more single than $\La'$. Hence (a) 
follows from 22.14(a) using 20.18, 20.19.

By 22.14(b), the difference $\aa_{[\ti{\La}]}-\aa_{[\La]}$ (where $\ti{\La}$ is
either $\La$ or $\La^I$ or $\La^{II}$) is $a$ times the number of $i<j$ in 
$[1,k]$. Thus, it is $a\binom{k}{2}$. Hence (a) follows from 20.18, 20.19. The
lemma is proved.

\proclaim{Lemma 22.17} $\boj_{\fS_k\tim W_l}^W(\sgn_k\ot[\La'])$ equals $[\La]$
or $[\La^I]+[\La^{II}]$.
\endproclaim
By a direct computation (involving representations of symmetric groups) we see
that:

(a) if $\La$ is defined then $[[\La']:[\La]]\ge 1$;

(b) if $\La^I,\La^{II}$ are defined then $[[\La']:[\La^I]]\ge 1$ and
$[[\La']:[\La^{II}]]\ge 1$.
\nl
In the setup of (a) we have (by 20.14(b)):

$g_{[\La']_\spa}=\sum_{E;\aa_E=\aa_{E'}}[[\La']:E]g_{E_\spa}$ hence using 
22.16(a) we have

(c) $g_{[\La]_\spa}=\sum_{E;\aa_E=\aa_{E'}}[[\La']:E]g_{E_\spa}$.
\nl
By 22.16(b), $E=[\La]$ enters in the last sum and its contribution is
$\ge g_{[\La]_\spa}$; the contribution of the other $E$ is $\ge 0$ (see 
20.13(b)). Hence (c) forces $[[\La']:[\La]]=1$ and $[[\La']:E]=0$ for all other
$E$ in the sum. In this case the lemma follows.

In the setup of (b) we have (by 20.14(b)):

$g_{[\La']_\spa}=\sum_{E;\aa_E=\aa_{E'}}[[\La']:E]g_{E_\spa}$ hence, using
22.16(a), we have

(d) $g_{[\La^I]_\spa}+g_{[\La^{II}]_\spa}
=\sum_{E;\aa_E=\aa_{E'}}[[\La']:E]g_{E_\spa}$.
\nl
By 22.16(b), $E=[\La^I]$ and $E=[\La^{II}]$ enter in the last sum and their
contribution is $\ge g_{[\La^I]_\spa}+g_{[\La^{II}]_\spa}$; the 
contribution of the other $E$ is $\ge 0$ (see 20.13(b)). Hence (d) forces 
$[[\La']:[\La^I]]=[[\La']:[\La^{II}]]=1$ and $[[\La']:E]=0$ for all other $E$ 
in the sum. The lemma follows.

\proclaim{Lemma 22.18}$[\La]\ot\sgn=[\bar\La]$. (Notation of 22.14.)
\endproclaim
This can be reduced to a known statement about the symmetric group. We omit the
details.

\subhead 22.19\endsubhead
Let $Z$ be a totally ordered finite set $z_1<z_2<\dots<z_M$. For any 
$r\in[0,M]$ such that $r=M\mod 2$ let $\und Z_r$ be the set of subsets of $Z$ 
of cardinal $(M-r)/2$. An involution $\io:Z@>>>Z$ is said to be $r$-admissible
if the following hold:

(a) $\io$ has exactly $r$ fixed points;

(b) if $M=r$, there is no further condition; if $M>r$, there exist two 
consecutive elements $z,z'$ of $Z$ such that $\io(z)=z',\io(z')=z$ and the
induced involution of $Z-\{z,z'\}$ is $r$-admissible.
\nl
Let $\Inv_r(Z)$ be the set of $r$-admissible involutions of $Z$. To 
$\io\in\Inv_r(Z)$ we associate a subset $\cs_\io$ of $\und Z_r$ as follows: a 
subset $Y\sub Z$ is in $\cs_\io$ if it contains exactly one element in each 
non-trivial $\io$-orbit. Clearly $|\cs_\io|=2^p$ where $p=(M-r)/2$. (In fact,
$\cs_\io$ is naturally an affine space over the field $F_2$.)

\proclaim{Lemma 22.20} Assume that $p>0$. Let $Y\in\und Z_r$.

(a) We can find two consecutive elements $z,z'$ of $Z$ such that exactly one of
$z,z'$ is in $Y$.

(b) There exists $\io\in\Inv_r(Z)$ such that $Y\in\cs_\io$.

(c) Assume that for some $k\in[0,p-1]$, $z_1,z_2,\dots,z_k$ belong to $Y$ but 
$z_{k+1}\notin Y$. Let $l$ be the smallest number such that $l>k$ and 
$z_l\in Y$. There exists $\io\in\Inv_r(Z)$ such that $Y\in\cs_\io$ and 
$\io(z_l)=z_{l-1}$.
\endproclaim
We prove (a). Let $z_k$ be the smallest element of $Y$. If $k>1$ then we can 
take $(z,z')=(z_{k-1},z_k)$. Hence we may assume that $z_1\in Y$. Let $z_{k'}$ 
be the next smallest elememt of $Y$. If $k'>2$ then we can take
$(z,z')=(z_{k'-1},z_{k'})$. Continuing like this we see that we may assume that
$Y=\{z_1,z_2,\dots,z_p\}$  Since $p<M$ we may take $(z,z')=(z_p,z_{p+1})$. 

We prove (b). Let $z,z'$ be as in (a). Let $Z'=Z-\{z,z'\}$ with the induced
order. Let $Y'=Y\cap Z'$. If $p\ge 2$ then by induction on $p$ we may assume
that there exists $\io'\in\Inv_r(Z')$ such that $Y'\in\cs_{\io'}$. Extend 
$\io'$ to an involution $\io$ of $Z$ by $z\mto z', z'\mto z$. Then 
$\io\in\Inv_r(Z)$ and $Y\in\cs_\io$. If $p=1$, define $\io:Z@>>>Z$ so that
$z\mto z', z'\mto z$ and $\io=1$ on $Z-\{z,z'\}$. Then $\io\in\Inv_r(Z)$ and 
$Y\in\cs_\io$.

We prove (c). We have $l\ge k+2$. Hence $z_{l-1}\notin Y$. Let 
$(z,z')=(z_{l-1},z_l)$. We continue as in the proof of (b), except that instead
of invoking an induction hypothesis, we invoke (b) itself.

\subhead 22.21\endsubhead
Assume that $M>r$. We consider the graph whose set of vertices is $\und Z_r$ 
and in which two vertices $Y\ne Y'$ are joined if there exists 
$\io\in\Inv_r(Z)$ such that $Y\in\cs_\io,Y'\in\cs_\io$.

\proclaim{Lemma 22.22} This graph is connected.
\endproclaim
We show that any vertex $Y=\{z_{i_1},z_{i_2},\dots,z_{i_p}\}$ is in the same 
connected component as $Y_0=\{z_1,z_2,\dots,z_p\}$. We argue by induction on
$m_Y=i_1+i_2+\dots+i_p$. If $m_Y=1+2+\dots+p$ then $Y=Y_0$ and there is nothing
to prove. Assume now that $m>1+2+\dots+p$ so that $Y\ne Y_0$. Then the 
assumption of Lemma 22.20(c) is satisfied. Hence we can find $l$ such that 
$z_l\in Y, z_{l-1}\notin Y$ and $\io\in\Inv_r(Z)$ such that $Y\in\cs_\io$ and
$\io(z_l)=z_{l-1}$. Let $Y'=(Y-\{z_l\})\cup\{z_{l-1}\}$. Then $Y'\in\cs_\io$
hence $Y,Y'$ are joined in our graph. We have $m_{Y'}=m_Y-1$ hence by the
induction hypothesis $Y',Y_0$ are in the same connected component. It follows
that $Y,Y_0$ are in the same connected component. The lemma is proved.

\subhead 22.23\endsubhead
Assume that $b'=0$. Let $\ti Z\in\cm_{a,b;n}^N$. Let $Z$ be the set of singles
of $\ti Z$. Each set $Y\in\und Z_r$ gives rise to a symbol $\La_Y$ in 
$\pi_N\i(\ti Z)$: the first row of $\La_Y$ consists of $Z-Y$ and one of each 
double of $\ti Z$; the second row consists of $Y$ and one of each double of 
$\ti Z$. For any $\io\in\Inv_r(Z)$ we define

$c(\ti Z,\io)=\opl_{Y\in\cs_\io}[\La_Y]\in\Mod W$. 

\proclaim{Proposition 22.24} (a) In the setup of 22.23, let $\io\in\Inv_r(Z)$.
Then $c(\ti Z,\io)\in Con(W)$.

(b) All constructible representations of $W$ are obtained as in (a).
\endproclaim
We prove (a) by induction on $n$. If $n=0$ the result is clear. Assume now that
$n\ge 1$. We may assume that $0$ is not a double of $\ti Z$. Let $at$ be the 
largest entry of $\ti Z$.

(A) Assume that there exists $i, 0\le i<t$, such that $ai$ does not appear in
$\ti Z$. Then $\ti Z$ is obtained from $\ti Z'\in\cm_{a,b;n-k}^N$ with $n-k<n$
by increasing each of the $k$ largest entries by $a$ and this set of largest 
entries is unambiguously defined. The set $Z'$ of singles of $\ti Z'$ is 
naturally in order preserving bijection with $Z$. Let $\io'$ correspond to 
$\io$ under this bijection. By the induction hypothesis,
$c(\ti Z',\io')\in Con(W_{n-k})$. Since, by 22.5, the sign representation 
$\sgn_k$ of $\fS_k$ is constructible, it follows that
$\sgn_k\boxt c(\ti Z',\io')\in Con(\fS_k\tim W_{n-k})$. Using 22.17, we have 
$$\boj_{\fS_k\tim W_{n-k}}^W(\sgn_k\boxt c(\ti Z',\io'))=c(\ti Z,\io)$$
hence $c(\ti Z,\io)\in Con(W)$.

(B) Assume that there exists $i, 0<i\le t$ such that $ai$ is a double of 
$\ti Z$. Let $\bar{\ti Z}$ be as in 22.8 (with respect to our $t$). Then $0$ is
not a double of $\bar{\ti Z}$ and the largest entry of $\bar{\ti Z}$ is $at$. 
Let $\bar Z$ be the set of singles of $\bar{\ti Z}$. We have 
$\bar Z=\{at-z|z\in Z\}$. Thus $\bar Z,Z$ are naturally in (order reversing) 
bijection under $j\mto at-j$. Let $\io'\in\Inv_r(\bar Z)$ correspond to $\io$ 
under this bijection. Since $at-ai$ does not appear in $\bar{\ti Z}$, (A) is 
applicable to $\bar{\ti Z}$. Hence $c(\bar{\ti Z},\io')\in Con(W)$. By 22.18 we
have $c(\bar{\ti Z},\io')\ot\sgn=c(\ti Z,\io)$ hence $c(\ti Z,\io)\in Con(W)$. 

(C) Assume that we are not in case (A) and not in case (B). Then 
$\ti Z=\{0,a,2a,\dots,ta\}=Z$. We can find $ia,(i+1)a$ in $Z$ such that $\io$ 
interchanges $ia,(i+1)a$ and induces on $Z-\{ia,(i+1)a\}$ an $r$-admissible
involution $\io_1$. We have

$\ti Z'=\{0,a,2a,\dots,ia,ia,(i+1)a,(i+2)a,\dots,(t-1)a\}\in\cm_{a,b;n-k}^N$
\nl
with $n-k<n$. The set of singles of $\ti Z'$ is 

$Z'=\{0,a,2a,\dots,(i-1)a,(i+1)a,\dots,(t-1)a\}$. 
\nl
It is in natural (order preserving) bijection with $Z-\{ia,(i+1)a\}$. Hence 
$\io_1$ induces $\io'\in\Inv_r(Z')$. By the induction hypothesis we have
$c(\ti Z',\io')\in Con(W_{n-k})$. Hence
$\sgn_k\boxt c(\ti Z',\io')\in Con(\fS_k\tim W_{n-k})$ where $\sgn_k$ is as in
(A). Using 22.17, we have 
$$\boj_{\fS_k\tim W_{n-k}}^W(\sgn_k\boxt c(\ti Z',\io'))=c(\ti Z,\io)$$
hence $c(\ti Z,\io)\in Con(W)$. This proves (a).

We prove (b) by induction on $n$. If $n=0$ the result is clear. Assume now that
$n\ge 1$. By an argument like the ones used in (B) we see that the class of 
representations of $W$ obtained in (a) is closed under $\ot\sgn$. Therefore, to
show that $C\in Con(W)$ is obtained in (a), we may assume that
$C=\boj_{\fS_k\tim W_{n-k}}^W(C')$ for some $k>0$ and some
$C'\in Con(\fS_k\tim W_{n-k})$. By 22.5 we have $C'=E\boxt C''$ where $E$ is a
simple $\fS_k$-module and $C''\in Con(W_{n-k})$. Using 22.5(a) we have
$E=\boj_{\fS_{k'}\tim\fS_{k''}}^{\fS_k}(\sgn\boxt E')$ where $k'+k''=k, k'>0$
and $E'$ is a simple $\fS_{k''}$-module. Let 
$\ti C=\boj_{\fS_{k''}\tim W_{n-k}}^{W_{n-k'}}(E'\ot C')\in Con(W_{n-k'})$.
Then $C=\boj_{\fS_{k'}\tim W_{n-k'}}^W(\sgn_{k'}\ot\ti C)$. By the induction 
hypothesis, $\ti C$ is of the form described in (a). Using an argument as in 
(A) or (C) we deduce that $C$ is of the form described in (a). 
The proposition is proved.

\proclaim{Proposition 22.25}Assume that $b'>0$. 

(a) Let $E\in\Irr W$. Then $E\in Con(W)$.

(b) All constructible representations of $W$ are obtained as in (a).
\endproclaim
We prove (a). We may assume that $E=[\La]$ where $\La\in\Sy_{a,b;n}^N$ does not
contain both $0$ and $b'$. We argue by induction on $n$. If $n=0$ the result is
clear. Assume now that $n\ge 1$.

(A) Assume that either (1) there exist two entries $z,z'$ of $\La$ such that 
$z'-z>a$ and there is no entry $z''$ of $\La$ such that $z<z''<z'$, or (2) 
there exists an entry $z'$ of $\La$ such that $z'\ge a$ and there is no entry 
$z''$ of $\La$ such that $z''<z'$. Let $\La'$ be the symbol obtained from $\La$
by substracting $a$ from each entry $\ti z$ of $\La$ such that $\ti z\ge z'$ 
and leaving the other entries of $\La$ unchanged. Then $\La'\in\Sy_{a,b;n-k}^N$
with $n-k<n$. By the induction hypothesis, $[\La']\in Con(W_{n-k})$. Since, by
22.5, the sign representation $\sgn_k$ of $\fS_k$ is constructible, it follows
that $\sgn_k\boxt [\La']\in Con(\fS_k\tim W_{n-k})$. Using 22.17, we have 
$\boj_{\fS_k\tim W_{n-k}}^W(\sgn_k\boxt [\La'])=[\La]$ hence $[\La]\in Con(W)$.

(B) Assume that there exist two entries $z,z'$ of $\La$ such that $0<z'-z<a$.
Let $t$ be the smallest integer such that $at+b'\ge\la_i$ for all $i\in[1,N+r]$
and $at\ge\mu_j$ for all $j\in[1,N]$. Let $\bar\La\in\Sy_{a,b;n}^{t+1-N-r}$ be
as in 22.8 with respect to this $t$. Then $\bar\La$ does not contain both $0$
and $b'$. Now (A) is applicable to $\bar\La$. Hence $[\bar\La]\in Con(W)$. By
22.18 we have $[\bar\La]\ot\sgn=[\La]$ hence $[\La]\in Con(W)$. 

(C) Assume that we are not in case (A) and not in case (B). Then the entries of
$\La$ are either $0,a,2a,\dots,ta$ or $b',a+b',2a+b',\dots,ta+b'$. This cannot
happen for $n\ge 1$. This proves (a).

The proof of (b) is entirely similar to that of 22.24(b). The proposition is 
proved.

\subhead 22.26\endsubhead 
We now assume that $n\ge 2$ and that $W'=W'_n$ is the kernel of
$\chi_n:W_n@>>>\pm 1$ in 22.10. We regard $W'_n$ as a Coxeter group with 
generators $s_1,s_2,\dots,s_{n-1}$ as in 22.9 and $s'_n=(n-1,n')((n-1)',n)$ 
(product of transpositions). Let $L:W'@>>>\bn$ be the weight function given by
$L(w)=al(w)$ for all $w$. Here $a>0$.

For $\La\in\Sy_{a,0}^N$ we denote by $\La^{tr}$ the symbol whose first (resp.
second) row is the second (resp. first) row of $\La$. We then have 
$\La^{tr}\in\Sy_{a,0}^N$. From the definitions we see that the simple 
$W_n$-modules $[\La],[\La^{tr}]$ have the same restriction to $W'$; 
this restriction is a simple $W'$-module $[\und\La]$ if $\La\ne\La^{\tr}$
and is a direct sum of two non-isomorphic simple $W'$-modules 
$[{}^I\und\La]$,$[{}^{II}\und\La]$ if $\La=\La^{\tr}$. In this way we see that

{\it the simple $W'$-modules are naturally in bijection with the set of
orbits of the involution of $\Sy_{a,0;n}$ induced by $\La\mto\La^{tr}$ except
that each fixed point of this involution corresponds to two simple 
$W'$-modules.}

Let $\ti Z\in\cm_{a,0;n}^N$. Let $Z$ be the set of singles of $\ti Z$. Assume 
first that $Z\ne\emptyset$. Each set $Y\in\und Z_0$ gives rise to a symbol 
$\La_Y$ in $\Sy_{a,0;n}^N$: the first row of $\La_Y$ consists of $Z-Y$ and one
of each double of $\ti Z$; the second row consists of $Y$ and one of each 
double of $\ti Z$. For any $\io\in\Inv_0(Z)$ we define $c(\ti Z,\io)\in\Mod W$
by 

$c(\ti Z,\io)\opl c(\ti Z,\io)=\opl_{Y\in\cs_\io}[\La_Y]\in\Mod W$.  
\nl
Note that $Y$ and $Z-Y$ have the same contribution to the sum. A proof entirely
similar to that of 22.24 shows that $c(\ti Z,\io)\in Con(W)$. Moreover, if 
$Z=\emptyset$ and $\La=\La^{tr}\in\Sy_{a,1;n}^N$ is defined by 
$\pi_N(\La)=\ti Z$, then $[{}^I\und\La^]\in Con(W)$ and 
$[{}^{II}\und\La]\in Con(W)$. All constructible representations of $W$ are 
obtained in this way.
 
\subhead 22.27\endsubhead 
Assume that $W$ is of type $F_4$ and that the values of $L:W@>>>\bn$ on $S$ are
$a,a,b,b$ where $a>b>0$.

{\it Case 1}. Assume that $a=2b$. There are four simple $W$-modules
$\rho_1,\rho_2,\rho_8,\rho_9$ (subscript equals dimension) with $\aa=3b$. Then

$\rho_1\opl\rho_2,\rho_1\opl\rho_8,\rho_2\opl\rho_9,\rho_8\opl\rho_9\in 
Con(W)$.
\nl
(They are obtained by $\boj$ from the $W_I$ of type $B_3$ with parameters 
$a,b,b$.) 

The simple $W$-modules $\rho_1^\dag,\rho_2^\dag,\rho_8^\dag,\rho_9^\dag$
have $\aa=15b$ and

$\rho_1^\dag\opl\rho_2^\dag,\rho_1^\dag\opl\rho_8^\dag,\rho_2^\dag\opl
\rho_9^\dag, \rho_8^\dag\opl\rho_9^\dag\in Con(W)$.
\nl
There are five simple $W$-modules $\rho_{12},\rho_{16},\rho_6,\rho'_6,\rho_4$
(subscript equals dimension) with $\aa=7b$. Then

$\rho_4\opl\rho_{16},\rho_{12}\opl\rho_{16}\opl\rho_6,
\rho_{12}\opl\rho_{16}\opl\rho'_6\in Con(W)$.
\nl
All $12$ simple $W$-modules other than the $13$ listed above, are 
constructible. All constructible representations of $W$ are thus obtained.

{\it Case 2}. Assume that $a\notin\{b,2b\}$. The simple $W$-modules 
$\rho_{12},\rho_{16},\rho_6,\rho'_6,\rho_4$ in Case 1 now have $\aa=3a+b$ and

$\rho_4\opl\rho_{16},\rho_{12}\opl\rho_{16}\opl\rho_6,
\rho_{12}\opl\rho_{16}\opl\rho'_6\in Con(W)$.
\nl
All $20$ simple $W$-modules other than the $5$ listed above, are 
constructible. All constructible representations of $W$ are thus obtained.

\subhead 22.28\endsubhead 
Assume that $W$ is of type $G_2$ and that the values of $L:W@>>>\bn$ on $S$ are
$a,b$ where $a>b>0$. Let $\rho_2,\rho'_2$ be the two $2$-dimensional simple 
$W$-modules. They have $\aa=a$ and $\rho_2\opl\rho'_2$ is constructible. 
All $4$ simple $W$-modules other than the $2$ listed above, are 
constructible. All constructible representations of $W$ are thus obtained.

\head 23. Two-sided cells\endhead
\subhead 23.1\endsubhead
We define a graph $\cg_W$ as follows. The vertices of $\cg_W$ are the simple
$W$-modules up to isomorphism. Two non-isomorphic simple $W$-modules 
are joined in $\cg_W$ if they both appear as components of some constructible
representation of $W$. Let $\cg_W/\sim$ be the set of connected components of 
$\cg_W$. The connected components of $\cg/\sim$ are determined explicitly by
the results in \S22 for $W$ irreducible.

For example, in the setup of 22.4,22.5 we have $\cg_W=\cg_W/\sim$. In the setup
of 22.24, $\cg_W/\sim$ is naturally in bijection with $\cm_{a,b;n}$. (Here, 
22.22 is used). In the setup of 22.25, we have $\cg_W=\cg_W/\sim$.

We show that:

(a) {\it if $E,E'$ are in the same connected component of $\cg_W$ then}
$E\sim_{\cl\Cal R}E'$.
\nl
We may assume that both $E,E'$ appear in some constructible representation of 
$W$. By 22.2, there exists a left cell $\Ga$ such that $[E:[\Ga]]\ne 0$, 
$[E':[\Ga]]\ne 0$. By 21.2, we have $[E_\spa:J^\Ga_\bc]\ne 0$, 
$[E'_\spa:J^\Ga_\bc]\ne 0$. Hence $E\sim_{\cl\Cal R}E'$, as desired.

\subhead 23.2\endsubhead
Let $c_W$ be the set of two-sided cells of $W,L$. Consider the (surjective) map
$\Irr W@>>>c_W$ which to $E$ associates the two-sided cell $\boc$ such that 
$E\sim_{\cl\Cal R}x$ for $x\in\boc$. From 23.1 we see that this map induces a 
(surjective) map 

(a) $\om_W:\cg_W/\sim@>>>c_W$.
\nl
We conjecture that $\om_W$ is a bijection. This is made plausible by:

\proclaim{Proposition 23.3} Assume that $W,L$ is split. Then $\om_W$ is a
bijection.
\endproclaim
Let $E,E'\in\Irr W$ be such that $E\sim_{\cl\Cal R}E'$. By 22.3, we can find 
constructible representations $C,C'$ such that $[E:C]\ne 0,[E':C']\ne 0$. By 
22.2, we can find left cells $\Ga,\Ga'$ such that $C=[\Ga],C'=[\Ga']$. Then 
$[E:[\Ga]]\ne 0,[E':[\Ga']]\ne 0$. Let $d\in\cd\cap\Ga,d'\in\cd\cap\Ga'$. Since
$\ga_d=[\Ga]$ and $[E:[\Ga]\ne 0$, we have 
$E\sim_{\cl\Cal R}d$. Similarly, $E'\sim_{\cl\Cal R}d'$. Hence 
$d'\sim_{\cl\Cal R}d'$. By
18.4(c), there exists $u\in W$ such that $t_dt_ut_{d'}\ne 0$. (Here we use the
splitness assumption.) Note that $j\mto jt_ut_{d'}$ is a $J_\bc$-linear map
$J^\Ga_\bc@>>>J^{\Ga'}_\bc$. This map is non-zero since it takes $t_d$ to 
$t_dt_ut_{d'}\ne 0$. Thus, $\Hom_{J_\bc}(J^\Ga_\bc,J^{\Ga'}_\bc)\ne 0$. Using 
21.2, we deduce that $\Hom_W([\Ga],[\Ga'])\ne 0$. Hence there exists 
$\ti E\in\Irr W$ such that $\ti E$ is a component of both $[\Ga]=C$ and 
$[\Ga']=C'$. Thus, both $E,\ti E$ appear in $C$ and both $\ti E,E'$ appear in 
$C'$. Hence $E,E'$ are in the same connected component of $\cg_W$. The 
proposition is proved.

\subhead 23.4\endsubhead
Assume now that $W,L$ is quasisplit, associated as in 16.2 to $\ti W$ and the 
automorphism $\si:\ti W@>>>\ti W$ with $\ti W$ finite, irreducible. 

Let $c_{\ti W}^!$ be the set of all $\si$-stable two-sided cells of $\ti W$.
Let $c_{\ti W}^{\sha}$ be the set of all two-sided cells of $\ti W$ which meet 
$W$. We have $c_{\ti W}^{\sha}\sub c_{\ti W}^!\sub c_{\ti W}$. Let 
$f:c_W@>>>c_{\ti W}^{\sha}$ be the map which attaches to a two-sided cell of 
$W$ the unique two-sided cell of $\ti W$ containing it; this map is well 
defined by 16.13(a) and is obviously surjective. 

\proclaim{Proposition 23.5} In the setup of 23.4, $\om_W$ is a bijection and
$f:c_W@>>>c_{\ti W}^{\sha}$ is a bijection.
\endproclaim
Since $\om_W,f$ are surjective, the composition 
$f\om_W:\cg_W/\sim@>>>c_{\ti W}^{\sha}$ is surjective. Hence it is enough to 
show that this composition is injective. For this it suffices to check one of 
the two statements below:

(a) $|\cg_W/\sim|=|c_{\ti W}^{\sha}|$;

(b) the composition $\cg_W/\sim@>f\om_W>>c_{\ti W}^{\sha}\sub c_{\ti W}@>f'>>
\bn\opl\bn$ (where $f'(\boc)=(\aa(x),\aa(xw_0))$ for $x\in\boc$) is injective.
\nl
Note that the value of the composition (b) at $E$ is $(\aa_E,\aa_{E^\dag})$.

{\it Case 1.} $W$ is of type $G_2$ and $\ti W$ is of type $D_4$. Then (b) 
holds: the composition (b) takes distinct values $(0,12),(1,7),(3,3),(7,1),
(12,0)$ on the $5$ elements of $\cg_W/\sim$. 

{\it Case 2.} $W$ is of type $F_4$ and  $\ti W$ is of type $E_6$. Then again
(b) holds.

{\it Case 3.} $W$ is of type $B_n$ with $n\ge 2$ and $\ti W$ is of type 
$A_{2n}$ or $A_{2n+1}$. Then $\si$ is conjugation by the longest element 
$\ti w_0$ of $\ti W$. We show that (a) holds.

Let $Y$ be the set of all $E\in\Irr\ti W$ (up to isomorphism) such that
$\tr(\ti w_0,E)\ne 0$. Let $Y'$ be the set of all $E'\in\Irr W$ (up to 
isomorphism). By 23.4 and 23.1 we have a natural bijection between $c_{\ti W}$
and the set of isomorphism classes of $E\in\Irr\ti W$. If $\boc\in c_{\ti W}$ 
corresponds to $E$, then the number of fixed points of $\si$ on $\boc$ is 
clearly $\pm\dim(E)\tr(\ti w_0,E)$. Hence $|c_{\ti W}^{\sha}|=|Y|$. From 23.1 
we have $|\cg_W/\sim|=|Y|$. Hence to show (a) it suffices to show that 
$|Y|=|Y'|$. But this is shown in \cite{\LU}.

{\it Case 4.} Assume that $\ti W$ is of type $D_n$ and $W$ is of type $B_{n-1}$
with $n\ge 3$. We will show that (a) holds. We change notation and write $W'$ 
instead of $\ti W$, $W'{}^\si$ instead of $W$. Then $W'$ is as in 22.26 and we 
may assume that $\si:W'@>>>W'$ is conjugation by $s_n$ (as in 22.26). Let 
$\cm_{1,0;n}^{N,!}$ be the set of all elements in $\cm_{1,0;n}^N$ whose set of
singles is non-empty. Let 

$\cm_{1,0;n}^!=\lim_{N\to\infty}\cm_{1,0;n}^{N,!}$.
\nl
By 22.26 and 23.3, $c_{W'}^!$ is naturally in bijection with $\cm_{1,0;n}^!$. 
By 23.1, $\cg_{W'{}^\si}/\sim$ is naturally in bijection with $\cm_{1,2;n-1}$. 
The identity map is clearly a bijection 
$\cm_{1,2;n-1}^N@>\sim>>\cm_{1,0;n}^{N+1,!}$. It induces a bijection
$\cm_{1,2;n-1}@>\sim>>\cm_{1,0;n}^!$. Hence to prove
$|\cg_{W'{}^\si}/\sim|\le|c_{W'}^{\sha}|$ it suffices to prove that
$|\cm_{1,0;n}^!|=|\cm_{1,0;n}^{\sha}|$. In other words, we must show that 

(c) {\it any $\si$-stable two-sided cell of $W'$ meets $W'{}^\si$.}
\nl
Now 26.26 and 23.3 provide an inductive procedure to obtain any $\si$-stable
two-sided cell of $W'$. Namely such a cell is obtained by one of two 
procedures:

(i) we consider a $\si$-stable two-sided cell in a parabolic subgroup of type 
$\fS_k\tim D_{n-k}$ (where $n-k\in[2,n-1]$) and we attach to it the unique
two-sided cell of $W'$ that contains it;

(ii) we take a two-sided cell obtained in (i) and multiply it on the right by
the longest element of $W'$.
\nl
Since we may assume that (c) holds when $n$ is replaced by $n-k\in[2,n-1]$, we
see that the procedures (i) and (ii) yield only two-sided cells that contain 
$\si$-fixed elements. This proves (c). The proposition is proved.

\head 24. Virtual cells\endhead
\subhead 24.1\endsubhead
In this section we preserve the setup of 20.1. 

A {\it virtual cell} of $W$ (with respect to $L:W@>>>\bn$) is an element of
$K(W)$ of the form $\ga_x$ for some $x\in W$.

\proclaim{Lemma 24.2} Let $x\in W$ and let $\Ga$ be the left cell containing
$x$. 

(a) If $\ga_x\ne 0$ then $x\in\Ga\cap\Ga\i$. 

(b) $\ga_x$ is a $\bc$-linear combination of $E\in\Irr W$ such that
$[E:[\Ga]]\ne 0$.
\endproclaim
Assume that $\ga_x\ne 0$. Then there exists $\ce\in\Irr J_\bc$ such that
$\tr(t_x,\ce)\ne 0$. We have $\ce=\opl_{d\in\cd}t_d\ce$ and $t_x:\ce@>>>\ce$ 
maps the summand $t_d\ce$ (where $x\sim_\cl d$) into $t_{d'}$, where 
$d'\sim_\cl x\i$ and all other summands to $0$. Since $\tr(t_x,\ce)\ne 0$, we 
must have $t_d\ce=t_{d'}\ce\ne 0$ hence $d=d'$ and $x\sim_\cl x\i$. This proves
(a).

We prove (b). Let $d\in\cd\cap\Ga$. Assume that $E\in\Irr W$ appears with 
$\ne 0$ coefficient in $\ga_x$. Then $\tr(t_x,E_\spa)\ne 0$. As we have 
seen in the proof of (a), we have $t_dE_\spa\ne 0$. Using 21.3,21.2, we
deduce $[E_\spa:J_\bc^\Ga]\ne 0$ and $[E_\spa:[\Ga]_\spa]\ne 0$.
Hence $[E:[\Ga]]\ne 0$. The lemma is proved.

\subhead 24.3\endsubhead
In the remainder of this section we will give a number of expicit computations
of virtual cells.

\proclaim{Lemma 24.4} In the setup of 22.10, $w_0$ acts on $[\La]$ as 
multiplication by

$\ep_{[\La]}=(-1)^{\sum_j(a\i\mu_j-j+1)}$.
\endproclaim
Using the definitions we are reduced to the case where $k=n$ or $l=n$. If $k=n$
we have $\ep_{[\La]}=1$ since $[\La]$ factors through $\fS_n$ and the longest 
element of $W_n$ is in the kernel of $W_n@>>>\fS_n$. Similarly, if $l=n$ we 
have $\ep_{[\La]}=\ep_{\chi_n}=(-1)^n$. The lemma is proved.

\proclaim{Proposition 24.5} Assume that we are in the setup of 22.23. Let 
$\io\in\Inv_r(Z)$ and let $\ka:\cs_\io@>>>F_2$ be an affine-linear function. 
Let
$c(\ti Z,\io,\ka)=\sum_{Y\in\cs_\io}(-1)^{\ka(Y)}[\La_Y]\in K(W)$. There exists
$x\in W$ such that $\ga_x=\pm c(\ti Z,\io,\ka)$.
\endproclaim
To some extent the proof is a repetition of the proof of 22.24(a), but we have
to keep track of $\ka$, a complicating factor.

We argue by induction on the rank $n$ of $\ti Z$. If $n=0$ the result is clear.
Assume now that $n\ge 1$. We may assume that $0$ is not a double of $\ti Z$. 
Let $at$ be the largest entry of $\ti Z$.

(A) Assume that there exists $i, 0\le i<t$, such that $ai$ does not appear in
$\ti Z$. Then $\ti Z$ is obtained from a multiset $\ti Z'$ of rank $n-k<n$ by
increasing each of the $k$ largest entries by $a$ and this set of largest 
entries is unambiguously defined. The set $Z'$ of singles of $\ti Z'$ is 
naturally in bijection with $Z$. 

Let $\io',\ka'$ correspond to $\io,\ka$ under this bijection. By the induction
hypothesis, there exists $x'\in W_{n-k}$ such that 
$\ga_{x'}^{W_{n-k}}=\pm c(\ti Z',\io',\ka')]\in K(W_{n-k})$. Let $w_{0,k}$ be 
the longest element of $\fS_k$. Then
$$\ga_{w_{0,k}x'}^{\fS_k\tim W_{n-k}}=
\ga_{w_{0,k}}^{\fS_k}\boxt\ga_{x'}^{W_{n-k}}=\sgn_k\boxt\ga_{x'}^{W_{n-k}}$$
and
$$\align&\ga_{w_{0,k}x'}^W=\boj_{\fS_k\tim W_{n-k}}^W
(\ga_{w_{0,k}x'}^{\fS_k\tim W_{n-k}})\\&=\boj_{\fS_k\tim W_{n-k}}^W
(\sgn_k\boxt\ga_{x'}^{W_{n-k}})=\pm\boj_{\fS_k\tim W_{n-k}}^W
(\sgn_k\boxt c(\ti Z',\io',\ka'))=\pm c(\ti Z,\io,\ka),\endalign$$
as required.

(B) Assume that there exists $i, 0<i\le t$ such that $ai$ is a double of 
$\ti Z$. Let $\bar{\ti Z}$ be as in 22.8 (with respect to our $t$). Then $0$ is
not a double of $\bar{\ti Z}$ and the largest entry of $\bar{\ti Z}$ is $at$. 
Let $\bar Z$ be the set of singles of $\bar{\ti Z}$. We have 
$\bar Z=\{at-z|z\in Z\}$. Thus $\bar Z,Z$ are naturally in (order reversing) 
bijection under $j\mto at-j$. Let $\io'\in\Inv_r(\bar Z)$ correspond to $\io$ 
this bijection and let $\ka':\cs_{\io'}@>>>F_2$ correspond to $\ka$ under this
bijection. Let $\ka'':\cs_{\io'}@>>>F_2$ be given by 
$\ka''(Y)=\ka'(Y)+\sum_{y\in Y}a\i y$ (an affine-linear function). Since 
$at-ai$ does not appear in $\bar{\ti Z}$, (A) is applicable to $\bar{\ti Z}$. 
Hence there exists $x'\in W$ such that 

$\ga_{x'}=\pm c(\bar{\ti Z},\io',\ka'')$. By 20.23, 22.18, 24.4, we have
$$\ga_{x'w_0}=(-1)^{l(x')}\ze(\ga_{x'})=\pm\ze(c(\bar{\ti Z},\io',\ka''))
\ot\sgn=\pm c(\bar{\ti Z},\io',\ka')\ot\sgn=\pm c(\ti Z,\io,\ka),$$
as desired.

(C) Assume that we are not in case (A) and not in case (B). Then 
$\ti Z=\{0,a,2a,\dots,ta\}=Z$. We can find $ia,(i+1)a$ in $Z$ such that $\io$ 
interchanges $ia,(i+1)a$ and induces on $Z-\{ia,(i+1)a\}$ an $r$-admissible
involution $\io_1$. 

(C1) Assume first that $\ka(Y)=\ka(Y*\{ia,(i+1)a\})$ for any $Y\in\cs_\io$. 
($*$ is symmetric difference.) Let 

$\ti Z'=\{0,a,2a,\dots,ia,ia,(i+1)a,(i+2)a,\dots,(t-1)a\}$.
\nl
This has rank $n-k<n$. The set of singles of $\ti Z'$ is 

$Z'=\{0,a,2a,\dots,(i-1)a,(i+1)a,\dots,(t-1)a\}$. 
\nl
It is in natural (order preserving) bijection with $Z-\{ia,(i+1)a\}$. Hence 
$\io_1$ induces $\io'\in\Inv_r(Z')$. We have an obvious surjective map of 
affine spaces $p:\cs_{\io}@>>>\cs_{\io'}$ and $\ka$ is constant on the fibres
of this map. Hence there is an affine-linear map $\ka':\cs_{\io'}@>>>F_2$ such
that $\ka=\ka'p$. By the induction hypothesis, there exists $x'\in W_{n-k}$ 
such that $\ga_{x'}^{W_{n-k}}=\pm c(\ti Z',\io',\ka')\in K(W_{n-k})$. Let 
$w_{0,k}$ be the longest element of $\fS_k$. Then
$$\ga_{w_{0,k}x'}^{\fS_k\tim W_{n-k}}=
\ga_{w_{0,k}}^{\fS_k}\boxt\ga_{x'}^{W_{n-k}}=\sgn_k\boxt\ga_{x'}^{W_{n-k}}$$
and
$$\align&\ga_{w_{0,k}x'}^W=\boj_{\fS_k\tim W_{n-k}}^W
(\ga_{w_{0,k}x'}^{\fS_k\tim W_{n-k}})\\&=\boj_{\fS_k\tim W_{n-k}}^W
(\sgn_k\boxt\ga_{x'}^{W_{n-k}})=\pm\boj_{\fS_k\tim W_{n-k}}^W
(\sgn_k\boxt c(\ti Z',\io',\ka'))=\pm c(\ti Z,\io,\ka),\endalign$$
as required.

(C2) Assume next that $\ka(Y)\ne\ka(Y*\{ia,(i+1)a\})$ for some (or equivalently
any) $Y\in\cs_\io$. We have
$$\bar{\ti Z}=\{0,0,a,a,2a,2a,\dots,ta,ta\}-\{at-0,at-a,\dots,at-at\}=\ti Z=Z.
$$
Let $\io'\in\Inv_r(Z)$ correspond to $\io$ under the bijection $z\mto ta-z$ of
$Z$ with itself; let $\ka':\cs_{\io'}@>>>F_2$ correspond to $\ka$ under this 
bijection. Let $\ka'':\cs_{\io'}@>>>F_2$ be given by 
$\ka''(Y)=\ka'(Y)+\sum_{y\in Y}a\i y$ (an affine-linear function). Note that 
$\io'$ interchanges $(t-i-1)a,(t-i)a$ and induces on $Z-\{(t-i-1)a,(t-i)a\}$ an
$r$-admissible involution. We show that for any $Y\in\cs_{\io'}$ we have 
$\ka''(Y)=\ka''(Y*\{(t-i-1)a,(t-i)a\})$ or equivalently 
$\ka'(Y)=\ka'(Y*\{(t-i-1)a,(t-i)a\})+1$. This follows from our assumption 
$\ka(Y)=\ka(Y*\{ia,(i+1)a\})+1$ for any $Y\in\cs_\io$. We see that case (C1) 
applies to $\io',\ka''$ so that there exists $x'\in W$ with
$\ga_{x'}=\pm c(\ti Z,\io',\ka'')$. By 20.23, 22.18, 24.4, we have 
$$\ga_{x'w_0}=(-1)^{l(x')}\ze(\ga_{x'})=\pm\ze(c(\bar{\ti Z},\io',\ka''))\ot
\sgn=\pm c(\bar{\ti Z},\io',\ka')\ot\sgn=\pm c(\ti Z,\io,\ka),$$
as desired. The proposition is proved.

\subhead 24.6\endsubhead
Assume that we are in the setup of 22.27. By 22.27,

$\rho_4+\rho_{16},\rho_{12}+\rho_{16}+\rho_6,\rho_{12}+\rho_{16}+\rho'_6$
\nl
are constructible, hence (by 22.2, 21.4) are of the form $n_d\ga_d$ for 
suitable
$d\in\cd$, hence are $\pm$ virtual cells. 

Let $d\in\cd$ be such that $n_d\ga_d=\rho_{12}+\rho_{16}+\rho_6$. Let $\Ga$ be
the left cell such that $d\in\Ga$. Recall (21.4) that $[\Ga]=A\opl B\opl C$ 
where $A=\rho_{12}, B=\rho_{16}, C=\rho_6$. By the discussion in
21.10 we see that $J_\bc^{\Ga\cap\Ga\i}$ has exactly three simple modules (up 
to isomorphism), namely $t_dA_\spa,t_dB_\spa,t_dC_\spa$, and these are 
$1$-dimensional. Since $J^{\Ga\cap\Ga\i}$ is a semisimple algebra (21.9), it 
follows that it is commutative of dimension $3$. Hence $\Ga\cap\Ga\i$ consists
of three elements $d,x,y$. Let $p_A,p_B,p_C$ denote the traces of $t_x$ on 
$A_\spa,B_\spa,C_\spa$ respectively. Let $q_A,q_B,q_C$ denote the traces
of $t_y$ on $A_\spa,B_\spa,C_\spa$ respectively. By 20.24, 
$p_A,p_B,p_C,q_A,q_B,q_C$ are integers. Recall that the traces of $n_dt_d$ on 
$A_\spa,B_\spa,C_\spa$ are $1,1,1$ respectively. Since 
$f_{A_\spa},f_{B_\spa},f_{C_\spa}$ are $6,2,3$ we see that the 
orthogonality formula 21.10 gives 

$1+p_A^2+q_A^2=6,1+p_B^2+q_B^2=2, 1+p_C^2+q_C^2=3$,

$1+p_Ap_B+q_Aq_B=0,1+p_Ap_C+q_Aq_C=0,1+p_Bp_C+q_Bq_C=0$.
\nl
Solving these equations with integer unknowns we see that there exist 
$\ep,\ep'\in\{1,-1\}$ so that (up to interchanging $x,y$) we have

$(p_A,q_A)=(2\ep,\ep'),(p_B,q_B)=(0,-\ep'),(p_C,q_C)=(-\ep,\ep')$.
\nl
Then $\ep\ga_x=2\rho_{12}-\rho_6$, $\ep'\ga_y=\rho_{12}-\rho_{16}+\rho_6$. 
Hence 

$2\rho_{12}-\rho_6, \rho_{12}-\rho_{16}+\rho_6$ are $\pm$ virtual cells.
\nl
Exactly the same argument shows that 

$2\rho_{12}-\rho_{6'},\rho_{12}-\rho_{16}+\rho_{6'}$ are $\pm$ virtual cells. 
\nl
A similar (but simpler) argument shows that

$\rho_4-\rho_{16}$ is $\pm$ a virtual cell. 
\nl
Assume now that we are in the setup of 22.27 (Case 1). By 22.27,

$\rho_1+\rho_2,\rho_1+\rho_8,\rho_2+\rho_9,\rho_8+\rho_9,
\rho_1^\dag+\rho_2^\dag,\rho_1^\dag+\rho_8^\dag,\rho_2^\dag+\rho_9^\dag,
\rho_8^\dag+\rho_9^\dag$,
\nl
are constructible, hence by 22.2, 21.4 are of the form $n_d\ga_d$ for suitable
$d\in\cd$, hence are $\pm$ virtual cells. By an argument similar to that above
(but simpler) we see that

$\rho_1-\rho_2,\rho_1-\rho_8,\rho_2-\rho_9,\rho_8-\rho_9,
\rho_1^\dag-\rho_2^\dag,\rho_1^\dag-\rho_8^\dag,\rho_2^\dag-\rho_9^\dag,
\rho_8^\dag-\rho_9^\dag$,
\nl
are $\pm$ virtual cells. 

\subhead 24.7\endsubhead
Assume that we are in the setup of 22.29. By 22.29, $\rho_2+\rho'_2$ is 
constructible, hence by 22.2, 21.4, is of the form $n_d\ga_d$ for some
$d\in\cd$, hence is $\pm$ a virtual cell. As in 24.6, we see that 
$\rho_2-\rho'_2$ is $\pm$ a virtual cell.

\enddocument